\documentclass[reqno,12pt]{amsart}
\usepackage[margin=1.1in, footskip=1cm]{geometry}
\usepackage{amsmath, amsthm, amssymb, booktabs}
\usepackage{mathrsfs}
\usepackage[usenames,dvipsnames]{color}
\pdfpagewidth=\paperwidth
\pdfpageheight=\paperheight

\newtheorem{theorem}{Theorem}[section]
\newtheorem{lemma}[theorem]{Lemma}
\newtheorem{corollary}[theorem]{Corollary}

\newtheorem{hypothesis}[theorem]{Hypothesis}

\theoremstyle{definition}

\theoremstyle{remark}

\numberwithin{equation}{section}

\newcommand{\mmod}[1]{\,\,({\rm{mod}}\,\,#1)}

\def\d{{\partial}}

\def\le{\leqslant} \def\ge{\geqslant}

\def\d{{\,{\rm d}}}

\makeatletter
\@namedef{subjclassname@2020}{\textup{2020} Mathematics Subject Classification}
\makeatother

\begin{document}
\title[Efficient congruencing]{Efficient congruencing in $\mathbb Z[t]$ modulo $t^2$}
\author[T. D. Wooley]{Trevor D. Wooley}
\address{Department of Mathematics, Purdue University, 150 N. University Street, West Lafayette, IN 47907-2067, USA}
\email{twooley@purdue.edu}
\subjclass[2020]{11L15, 11L07, 11P55}
\keywords{Exponential sums, Hardy-Littlewood method, congruences}
\date{}

\begin{abstract} We develop an analogue of nested efficient congruencing applicable in the ring 
$\mathbb Z[t]/t^2\mathbb Z[t]$ so as to estimate mean values related to that of Vinogradov. Let $k$ be a natural 
number, and suppose that $\varphi_j(u)\in \mathbb Z[t,u]$ $(1\le j\le k)$ is a system of polynomials having Wronskian 
$W(u;\boldsymbol \varphi)=\det ({\rm d}^i\varphi_j(u)/{\rm d}u^i)_{1\le i,j\le k}$ with respect to $u$. Consider the 
number $\Upsilon_{s,k}(X;\boldsymbol \varphi)$ of integral solutions of the system of equations
\[
\sum_{i=1}^s \bigl( \varphi_j(u_i+v_it)-\varphi_j(u_{s+i}+v_{s+i}t)\bigr) \equiv 0 \mmod{t^2}\quad (1\le j\le k),
\]
with $|u_i|\le X$ and $|v_i|\le X$ $(1\le i\le s)$. We show that when $W(u;\boldsymbol \varphi)\not\equiv 0\mmod{t}$ 
and $1\le s\le k(k+1)/2$, then $\Upsilon_{s,k}(X;\boldsymbol \varphi)\ll X^{2s+\varepsilon}$. Our conclusions extend 
to number fields and function fields, as we outline. We describe several consequences of our results, including 
forthcoming work establishing the main conjecture for certain systems of equations not possessing 
translation-dilation invariant properties.
\end{abstract}

\maketitle

\section{Introduction}
Our purpose in this memoir is to initiate the quantitative theory of congruences over $\mathbb Z[t]$ as they relate to 
relatives of Vinogradov's mean value theorem. In order to illustrate what we have in mind, consider a natural number 
$k$ and polynomials $\varphi_j(u)\in \mathbb Z [t,u]$ $(1\le j\le k)$. For the sake of expositional concision here, we 
have suppressed the dependence on $t$ in our notation for the polynomials $\varphi_j(u)$. Our focus lies on the 
system of congruences
\begin{equation}\label{1.1}
\sum_{i=1}^s \bigl( \varphi_j(x_i)-\varphi_j(x_{s+i})\bigr) \equiv 0\mmod{t^2}\quad (1\le j\le k),
\end{equation}
in which $s\in \mathbb N$, and for $1\le i\le 2s$ the variables $x_i$ lie in $\mathbb Z[t]$. There is no loss of 
generality here, of course, in supposing that $x_i=u_i+v_it$ for some $u_i,v_i\in \mathbb Z$ $(1\le i\le 2s)$. Our goal 
will be to obtain estimates for the number $\Upsilon_{s,k}(X;\boldsymbol \varphi)$ of solutions of the system 
\eqref{1.1} with $|u_i|\le X$ and $|v_i|\le X$ $(1\le i\le 2s)$. There is nothing particularly significant about the modulus 
$t^2$ chosen in \eqref{1.1}, aside from the simplicity of exposition. However, as we explain elsewhere 
\cite{Woo2026, Woo2026a}, it transpires that the estimates that we are able to obtain herein may be employed to 
obtain essentially optimal upper bounds for mean values associated with certain systems of equations failing to exhibit 
translation-dilation invariance. Since the latter systems are not accessible to techniques based on either decoupling or 
efficient congruencing methods, the circle of ideas presented here and in \cite{Woo2026, Woo2026a} have the merit of 
permitting precise investigation of problems hitherto inaccessible.\par

In order to announce the main conclusions of this memoir, we must record some additional notation. We consider 
polynomials $\varphi_j$ as introduced above, and we define the Wronskian of derivatives
\begin{equation}\label{1.2}
W(u;\boldsymbol \varphi)=\det \bigl( \varphi_j^{(i)}(u)\bigr)_{1\le i,j\le k},
\end{equation}
in which we write $\varphi_j^{(r)}(u)$ for the $r$-th derivative ${\rm d}^r\varphi_j(u)/{\rm d}u^r$. The system of 
polynomials $\boldsymbol \varphi=(\varphi_1,\ldots ,\varphi_k)$ may be considered independent modulo $t$ when 
$W(u;\boldsymbol \varphi)\not \equiv 0\mmod{t}$. We recall that the counting function 
$\Upsilon_{s,k}(X;\boldsymbol \varphi)$ associated with the system of equations \eqref{1.1} counts the number of 
integral solutions of the system of congruences
\begin{equation}\label{1.3}
\sum_{i=1}^s \bigl( \varphi_j(u_i+v_it)-\varphi_j(u_{s+i}+v_{s+i}t)\bigr) \equiv 0\mmod{t^2}\quad (1\le j\le k),
\end{equation}
with $|u_i|\le X$ and $|v_i|\le X$ $(1\le i\le 2s)$.

\begin{theorem}\label{theorem1.1}
Suppose that $\varphi_j(u)\in \mathbb Z[t,u]$ $(1\le j\le k)$ is a system of polynomials with 
$W(u;\boldsymbol \varphi)\not\equiv 0\mmod{t}$. Let $s$ be a positive integer with $s\le k(k+1)/2$. Then for each 
$\varepsilon>0$, one has
\begin{equation}\label{1.4}
\Upsilon_{s,k}(X;\boldsymbol \varphi)\ll X^{2s+\varepsilon}.
\end{equation}
\end{theorem}

Here and elsewhere, unless indicated otherwise, the implicit constants in Vinogradov's notation $\ll $ and $\gg $ 
depend at most on $\varepsilon$, $s$, $k$ and the coefficients of $\boldsymbol \varphi$. We note that the number 
of available choices for the variables $\mathbf u$ and $\mathbf v$ is essentially $(2X+1)^{4s}$, so the conclusion of  
Theorem \ref{theorem1.1} achieves an upper bound of size roughly equal to the square-root of this trivial estimate 
when $s\le k(k+1)/2$.\par

We illuminate the conclusion of Theorem \ref{theorem1.1} by highlighting a number of special cases as corollaries. 
We first distinguish the special case in which the polynomials $\varphi_j(u)$ are independent of $t$. In this scenario, 
it is convenient to define $\Upsilon^*_{s,k}(X;\boldsymbol \varphi)$ to be the number of integral solutions of the 
Diophantine system
\begin{equation}\label{1.5}
\left. \begin{aligned}
\sum_{i=1}^s\bigl( \varphi_j(u_i)-\varphi_j(u_{s+i})\bigr)&=0\\
\sum_{i=1}^s\bigl( v_i\varphi_j'(u_i)-v_{s+i}\varphi_j'(u_{s+i})&=0
\end{aligned} \right\} (1\le j\le k),
\end{equation}
with $|u_i|\le X$ and $|v_i|\le X$ $(1\le i\le 2s)$.

\begin{corollary}\label{corollary1.2}
Suppose that $\varphi_j\in \mathbb Z[u]$ $(1\le j\le k)$ is a system of polynomials with 
$W(u;\boldsymbol \varphi)\ne 0$. Let $s$ be a positive integer with $s\le k(k+1)/2$. Then for each $\varepsilon>0$, 
one has $\Upsilon^*_{s,k}(X;\boldsymbol \varphi)\ll X^{2s+\varepsilon}$.
\end{corollary}

In order to confirm that Corollary \ref{corollary1.2} follows from Theorem \ref{theorem1.1}, one has only to observe 
that since $\varphi_j$ is independent of $t$ for $1\le j\le k$, one has
\[
\varphi_j(u+vt)\equiv \varphi_j(u)+v\varphi'_j(u)t\mmod{t^2}.
\]
By equating coefficients of $1$ and $t$, we consequently see that the systems \eqref{1.3} and \eqref{1.5} are 
equivalent, whence $\Upsilon_{s,k}(X;\boldsymbol \varphi)=\Upsilon^*_{s,k}(X;\boldsymbol \varphi)$.\par

By orthogonality, the conclusion of Corollary \ref{corollary1.2} delivers upper bounds for mean values of exponential 
sums. As usual, we write $e(z)$ for ${\rm e}^{2\pi {\rm i}z}$. Then, given $k$-tuples $\boldsymbol \alpha$ and 
$\boldsymbol \beta$, we put
\[
\psi_k(u,v;\boldsymbol \alpha ,\boldsymbol \beta)=\sum_{j=1}^k\bigl( \alpha_j\varphi_j(u)+\beta_jv\varphi_j'(u)\bigr) ,
\]
and define
\[
F_k(\boldsymbol \alpha ,\boldsymbol \beta;X)=\sum_{|u|\le X}\sum_{|v|\le X}e(\psi_k(u,v;\boldsymbol \alpha ,
\boldsymbol \beta)).
\]
Thus, we see that under the hypotheses of Corollary \ref{corollary1.2}, one has
\begin{equation}\label{1.6}
\int_{[0,1)^k}\int_{[0,1)^k}|F_k(\boldsymbol \alpha ,\boldsymbol \beta;X)|^{2s}\d\boldsymbol \alpha \d\boldsymbol 
\beta \ll X^{2s+\varepsilon }.
\end{equation}

\par A further specialisation yields a concrete example more readily appreciated. We put $\varphi_j(u)=u^{d_j}$ 
$(1\le j\le k)$, and observe that $W(u;\boldsymbol \varphi)\ne 0$ whenever $1\le d_1<d_2<\ldots <d_k$. We now take 
$\Upsilon_{s,k}^+(X;\mathbf d)$ to be the number of integral solutions of the system of equations
\begin{equation}\label{1.7}
\left. \begin{aligned} \sum_{i=1}^s (u_i^{d_j}-u_{s+i}^{d_j})&=0\\
\sum_{i=1}^s(v_iu_i^{d_j-1}-v_{s+i}u_{s+i}^{d_j-1})&=0\end{aligned} \right\} (1\le j\le k),
\end{equation}
with $|u_i|\le X$ and $|v_i|\le X$ $(1\le i\le 2s)$.

\begin{corollary}\label{corollary1.3}
Let $d_1,\ldots ,d_k$ be distinct positive integers, and let $s$ be a positive integer with $s\le k(k+1)/2$. Then for 
each $\varepsilon>0$, one has $\Upsilon_{s,k}^+(X;\mathbf d)\ll X^{2s+\varepsilon}$.
\end{corollary}

Finally, we specialise further to the case that $d_j=j$ $(1\le j\le k)$. This special case is an extension of Vinogradov's 
mean value theorem that deserves special attention. We denote by $J_{s,k}^*(X)$ the number of integral solutions of 
the system
\begin{equation}\label{1.8}
\left. \begin{aligned} \sum_{i=1}^s(u_i^j-u_{s+i}^j)&=0\\
\sum_{i=1}^s (v_iu_i^{j-1}-v_{s+i}u_{s+i}^{j-1})&=0\end{aligned} \right\} (1\le j\le k),
\end{equation}
with $|u_i|\le X$ and $|v_i|\le X$ $(1\le i\le 2s)$. The quantity $J_{s,k}^*(X)$ counts the number of solutions of the 
Vinogradov system
\[
\sum_{i=1}^s \bigl( (u_i+v_it)^j-(u_{s+i}+v_{s+i}t)^j\bigr)\equiv 0\mmod{t^2}\quad (1\le j\le k),
\]
with $|u_i|\le X$ and $|v_i|\le X$ $(1\le i\le 2s)$. A consequence of Corollary \ref{corollary1.3} is the upper bound 
supplied in our final corollary.

\begin{corollary}\label{corollary1.4}
Let $s$ be a positive integer. Then for each $\varepsilon>0$, one has
\[
J_{s,k}^*(X)\ll X^\varepsilon (X^{2s}+X^{4s-k(k+1)}).
\]
\end{corollary} 
 
This corollary is sharp up to the factor of $X^\varepsilon$. On one hand, a consideration of those diagonal solutions 
of \eqref{1.8} with
\begin{equation}\label{1.9}
\{ u_1+v_1t,\ldots ,u_s+v_st\} =\{ u_{s+1}+v_{s+1}t,\ldots ,u_{2s}+v_{2s}t\}
\end{equation}
reveals that $J_{s,k}^*(X)\gg X^{2s}$ for $s\in \mathbb N$. Meanwhile, a crude major arc estimation of standard type 
reveals that
\[
J_{s,k}^*(X)\gg X^{4s-k(k+1)}.
\]
The significance of the exponent $k(k+1)$ is, of course, that it is equal to the sum of the degrees of the polynomials 
defining the system of equations \eqref{1.8}. We note that the conclusion of Corollary \ref{corollary1.4} is made 
available, in all essentials, by \cite[Theorem 1.1]{GZK2020}. The statement of the latter result is sufficiently heavy in 
notation that the reader will likely require some measure of determination to find the desired conclusion, although 
the reader will find \cite[Example 1.4]{GZK2020} very helpful. This work of Guo and Zorin-Kranich improves on earlier 
work of Parsell, Prendiville and Wooley \cite{PPW2013} that applies to quite general translation-dilation invariant 
systems of higher dimension (see also \cite{BDG2017, GZ2019}).\par

A key feature of the system \eqref{1.8} is that it is fully translation-dilation invariant. Thus, whenever 
$(u_i,v_i)_{1\le i\le 2s}$ gives a solution of the system \eqref{1.8}, then for any integers $a$ and $b$, the tuple 
$(u_i+a,v_i+b)_{1\le i\le 2s}$ also provides a solution of the system \eqref{1.8}. This feature is critical for the 
decoupling approach to bounding $J_{s,k}^*(X)$ in \cite{GZK2020}, as well as the efficient congruencing approach of 
\cite{PPW2013}. The approach of this paper implicitly makes use of the observation that the system \eqref{1.5} is 
approximated by the system \eqref{1.8}, so that translation-dilation invariance also plays a critical role in our 
present treatment. This translation-dilation invariance property, however, creates difficulties for us, in that additional 
special subvarieties are present in the solution set of the system \eqref{1.8} possessing a significant number of 
solutions. Thus, while the aforementioned diagonal solutions \eqref{1.9} contribute roughly $s!(2X+1)^{2s}$ 
solutions to the total counted by $J_{s,k}^*(X)$, there are additional solutions making a contribution of the same 
order of magnitude. Indeed, the solutions of \eqref{1.8} with
\[
u_1=u_2=\ldots =u_{2s}\quad \text{and}\quad v_1+\ldots +v_s=v_{s+1}+\ldots +v_{2s}
\]
contribute $\gg X\cdot X^{2s-1}=X^{2s}$ solutions to $J_{s,k}^*(X)$. We shall return to this observation when 
conditioning variables in \S6.\par

Our approach to proving Theorem \ref{theorem1.1} is not merely some kind of efficient congruencing analogue of the 
decoupling method of \cite{GZK2020}. The latter approach must accommodate the bivariate nature, involving $u_i$ 
and $v_i$, of the underlying monomials of the system \eqref{1.8}. Instead, we apply an efficient congruencing 
approach making intrinsic use of the arithmetic of $\mathbb Z[t]/t^2\mathbb Z[t]$, working in the univariate system 
\eqref{1.1}. This approach is flexible, permitting various generalisations beyond the test case examined in this memoir. 
We are confident that variants of our methods should apply to systems analogous to \eqref{1.1} in which the single 
variable polynomials $\varphi_j(x)$ are replaced by many variable variants, with $x_i$ a polynomial in $t$ of higher 
degree, and with the modulus $t^2$ replaced by a higher degree polynomial. We regard the present memoir as a 
proof-of-concept excursion on the theme.\par

A second important feature of the efficient congruencing approach, as opposed to real decoupling approaches, is 
that it is applicable in number fields, function fields, and indeed other settings in which one studies discrete sets of 
points in valued fields. We briefly sketch out both kinds of conclusions, and the ideas necessary for proofs, in \S\S11 
and 12, handling number field analogues in \S11 and function field analogues in \S12.\par

Before preparations begin for our proof of Theorem \ref{theorem1.1}, we pause to announce the main conclusion of a 
second memoir \cite{Woo2026} that makes critical use of Corollary \ref{corollary1.4}. This concerns the number 
$I_{s,k}(X)$ of integral solutions of the system of equations
\begin{equation}\label{1.10}
\sum_{i=1}^s(x_i^j-y_i^j)=0\quad (2\le j\le k),
\end{equation}
with $1\le x_i,y_i\le X$ $(1\le i\le s)$.

\begin{theorem}\label{theorem1.4}
When $k$ and $s$ are natural numbers with $k\ge 2$, and $\varepsilon >0$, one has
\[
I_{s,k}(X)\ll X^{s+\varepsilon}+X^{2s-(k^2+k-2)/2}.
\]
Thus, whenever $s$ and $\varepsilon$ are positive numbers, one has
\[
\int_{[0,1)^{k-1}}\biggl| \sum_{1\le x\le X}e(\alpha_2 x^2+\ldots +\alpha_kx^k)\biggr|^{2s}\d\boldsymbol \alpha \ll 
X^{s+\varepsilon}+X^{2s-(k^2+k-2)/2}.
\]
\end{theorem}

The conclusion of this theorem may be compared with a general consequence 
\cite[Corollary 1.2 and Theorem 14.1]{Woo2019} of efficient congruencing, which would achieve the same bound in 
the ranges $1\le s\le k(k-1)/2$ and $s\ge k(k+1)/2$, and \cite[Theorem 1.1]{BW2017}, which enhances the former 
range to $1\le s\le (k^2-1)/2$. Since elementary arguments yield the lower bound
\[
I_{s,k}(X)\gg X^s+X^{2s-(k^2+k-2)/2},
\]
it is evident that our upper bounds achieve the main conjecture associated with these mean value estimates. We 
stress that the system \eqref{1.10} is not translation-dilation invariant, and that, aside from the special case $k=2$ of 
the theorem, this is the first instance in which the main conjecture has been established for an incomplete 
Vinogradov system.\par

This memoir is organised as follows. We discuss Fourier analysis relevant to our work in 
$\mathbb Z[t]/t^r\mathbb Z[t]$ in \S2. This prepares the way for our detailed account of the proof of Theorem 
\ref{theorem3.2} in \S\S3-9, which provides the basic nested inductive step in the argument. The key infrastructure is 
described in \S3, with implications for translation-dilation invariant families in \S4. The nested structure is inductively 
built, with the case $k=1$ as the base of the induction. We discuss this base case in \S5. The general situation with 
$k\ge 2$ is analysed in \S\S6-9, with an initial conditioning step in \S6, the extraction of strong congruence 
constraints in \S7, and iterative procedures in \S8. The proof of the key Theorem \ref{theorem3.2} is then concluded 
in \S9. We advance in \S10 to derive consequences of Theorem \ref{theorem3.2} for mean values and equations. This 
work provides the conclusion embodied in Theorem \ref{theorem1.1} and its corollaries. In \S11 we describe how to 
adapt our basic methods so as to derive conclusions analogous to Theorem \ref{theorem1.1} for number fields, and 
in \S12 we provide a similar discussion for function fields. Finally, in \S13, we provide an appendix demonstrating 
how rather more general bivariate translation-dilation invariant systems may be interpreted in the framework 
provided by congruences over $\mathbb Z[t]$.\par   

Some notational conventions will streamline our exposition. Throughout, we make liberal use of vector notation in 
settings not always conventional in nature. Thus, for example, we may write $1\le \mathbf x\le X$ to denote that 
every coordinate $x_i$ of $\mathbf x$ satisfies $1\le x_i\le X$, and $\mathbf x\equiv \xi\mmod{p^h}$ to denote that 
$x_i\equiv \xi\mmod{p^h}$ for all indices $i$. Also, we adopt the convention that when 
$F:[0,1)^n\rightarrow \mathbb C$ is integrable, then
\[
\oint F(\boldsymbol \alpha)\d\boldsymbol \alpha =\int_{[0,1)^n}F(\boldsymbol \alpha)\d\boldsymbol \alpha .
\]
Finally, given a situation in which the parameter $\varepsilon$ has not already been fixed, in any statement involving 
the letter $\varepsilon$, it is implicitly asserted that the statement holds for each $\varepsilon>0$. In such 
circumstances, implicit constants in Vinogradov's notation $\ll $ and $\gg $ may depend on $\varepsilon$.\par

Work on this paper was conducted while the author was supported by NSF grants DMS-2001549 and DMS-2502625, 
and by Simons Fellowship in Mathematics SFM-00011955. The author is grateful to the Institute for Advanced Study, 
Princeton, for hosting his sabbatical, during which period this paper was completed.

\section{Fourier analysis in $\mathbb Z[t]/t^r\mathbb Z[t]$}
A feature of the work presented in this memoir is that we require the use of Fourier analysis applicable in the setting 
of $\mathbb Z[t]/t^2\mathbb Z[t]$. The basic foundations of this Fourier analysis are simple to formulate in the 
more general setting $\mathbb Z[t]/t^r\mathbb Z[t]$, with $r$ any fixed natural number. We discuss such matters in 
this section, though in the remainder of this memoir we restrict attention to the case $r=2$. Workers familiar with 
the application of the circle method in $\mathbb F_q[t]$ will experience no difficulty in deriving for themselves the 
necessary infrastructure for such endeavours. It is surely of assistance to the reader, however, to expend a modicum 
of effort reviewing such ideas in order to fix notation.\par

We shall be counting solutions of equations defined over $\mathbb Z[t]/t^r\mathbb Z[t]$, a ring that we abbreviate 
to $\mathfrak Z(t^r)$. A typical element of this ring $\mathfrak Z(t^r)$ takes the form 
\[
n_0+n_1t+\ldots +n_{r-1}t^{r-1},
\]
where $n_i\in \mathbb Z$ $(0\le i<r)$. With the possibility of confusion an omnipresent hazard, we shall write such a 
polynomial as $n$, relying on context to disambiguate this polynomial from a natural number. The need for myriad 
items of notation forces such expedience upon us.\par

As usual, we write $e(z)$ for ${\rm e}^{2\pi {\rm i}z}$ when $z\in \mathbb R$, noting the familiar relation
\[
\int_0^1 e(h\alpha)\d\alpha =\begin{cases}1,&\text{when $h=0$,}\\
0,&\text{when $h\in \mathbb Z\setminus \{ 0\}$,}\end{cases}
\]
that forms the foundation of Fourier analysis. This orthogonality relation detects vanishing in $\mathbb Z$. An 
analogous relation appropriate for $\mathfrak Z(t^r)$ makes use of the set
\[
\mathfrak T(t^r)=\Bigl\{ \frac{\alpha_1}{t}+\frac{\alpha_2}{t^2}+\ldots +\frac{\alpha_r}{t^r}: (\alpha_1,\ldots ,\alpha_r)
\in [0,1)^r\Bigr\} .
\] 
The set $\mathfrak T(t^r)$ forms a subset of the field of formal Laurent series
\[
\mathbb R((1/t))=\biggl\{ \sum_{m>-\infty}\frac{\alpha_m}{t^m}:\text{$\alpha_m\in \mathbb R$ for 
$m\in \mathbb Z$}\biggr\} .
\]
When $\xi\in \mathbb R((1/t))$, we define the {\it residue} of $\xi$, written $\text{res}(\xi)$, to be the coefficient 
$\xi_1$ of $1/t$ in the Laurent series expansion of $\xi$. We then define the exponential function 
$E:\mathbb R((1/t))\rightarrow \mathbb C$ by putting
\[
E(\xi)=e(\text{res}(\xi)).
\]
Thus, when 
\[
n=n_0+n_1t+\ldots +n_{r-1}t^{r-1}\in \mathfrak Z(t^r)\quad \text{and}\quad \alpha=\frac{\alpha_1}{t}+
\frac{\alpha_2}{t^2}+\ldots +\frac{\alpha_r}{t^r}\in \mathfrak T(t^r),
\]
one finds that
\[
\text{res}(n\alpha)=n_0\alpha_1+n_1\alpha_2+\ldots +n_{r-1}\alpha_r,
\]
and hence
\[
\int_0^1\!\!\ldots \int_0^1 E(n\alpha)\d\alpha_1\cdots \d\alpha_r=\begin{cases} 
1,&\text{when $n_i=0$ $(0\le i<r)$,}\\
0,&\text{otherwise.}
\end{cases}
\]

\par By adopting the convention that for integrable functions $F:\mathfrak T(t^r)\rightarrow \mathbb C$, we write
\[
\int_{\mathfrak T(t^r)}F(\alpha )\d\alpha =\int_0^1\!\!\ldots \int_0^1 F\Bigl( \frac{\alpha_1}{t}+\ldots 
+\frac{\alpha_r}{t^r}\Bigr) \d\alpha_1\cdots \d\alpha_r,
\]
we find that
\[
\int_{\mathfrak T(t^r)}E(n\alpha)\d\alpha =\begin{cases}1,&\text{when $n=0$,}\\
0,&\text{when $n\in \mathfrak Z(t^r)\setminus \{0\}$.}\end{cases}
\]
This discussion generalises at once to consider orthogonality in the ring $\mathbb Z[t]$ modulo $t^r$. By this we 
mean that when $m\in \mathbb Z[t]$ is a polynomial having degree $R$, in which $R$ might conceivably be larger 
than $r-1$, then a moment's reflection reveals that one has
\begin{equation}\label{2.1}
\int_{\mathfrak T(t^r)}E(m\alpha)\d\alpha =\begin{cases}1,&\text{when $m\equiv 0\mmod{t^r}$,}\\
0,&\text{when $m\not\equiv 0\mmod{t^r}$.}\end{cases}
\end{equation}

\par Now consider a prime number $p$ and a positive integer $b$. We write $(p^b,t^r)$ for the ideal of 
$\mathbb Z[t]$ generated by $p^b$ and $t^r$. In order to avoid ambiguity with the superficial appearance of a 
numerical value of a greatest common divisor, we note that, in concrete terms, one of course has
\[
(p^b,t^r)=\{ zp^b+wt^r: z,w\in \mathbb Z\}.
\]
Much of our discussion will occur in the module
\[
\mathbb Z[t]/(p^b,t^r)\cong (\mathbb Z/p^b\mathbb Z)\oplus \ldots \oplus (\mathbb Z/p^b\mathbb Z),
\]
with basis elements $1,t,\ldots ,t^{r-1}$. When $u\in \mathbb Z$ and $q\in \mathbb N$, write $e_q(u)=e(u/q)$. Then 
a little thought reveals that arguments very similar to those above lead from the fundamental orthogonality relation
\[
\frac{1}{p^b}\sum_{a=1}^{p^b}e_{p^b}(ua)=\begin{cases} 1,&\text{when $u\equiv 0\mmod{p^b}$,}\\
0,&\text{when $u\not\equiv 0\mmod{p^b}$,}\end{cases}
\]
to an orthogonality relation relevant for polynomials $m\in \mathbb Z[t]$, namely
\begin{equation}\label{2.2}
\frac{1}{p^{rb}}\sum_{a_1=1}^{p^b}\ldots \sum_{a_r=1}^{p^b} e_{p^b}\Bigl( 
\text{res}\Bigl( m\Bigl( \frac{a_1}{t}+\ldots +\frac{a_r}{t^r}\Bigr) \Bigr) \Bigr) =\begin{cases} 
1,&\text{when $m\equiv 0\mmod{(p^b,t^r)}$,}\\
0,&\text{when $m\not\equiv 0\mmod{(p^b,t^r)}$.}\end{cases}
\end{equation}
We seek to think of the relation \eqref{2.2} as an orthogonality relation streamlined as in \eqref{2.1}. Thus, we adopt 
the convention that given $F:\mathfrak T(t^r)\rightarrow \mathbb C$, we write
\[
\oint_{p^b} F(\alpha^{(r)})\d\alpha^{(r)}=\frac{1}{p^{rb}}\sum_{a_1=1}^{p^b}\ldots \sum_{a_r=1}^{p^b} 
F\biggl( \frac{1}{p^b}\Bigl(\frac{a_1}{t}+\ldots +\frac{a_r}{t^r}\Bigr)\biggr) .
\]
We then have
\begin{equation}\label{2.3}
\oint_{p^b}E(m\alpha^{(r)})\d\alpha^{(r)}=\begin{cases} 1,&\text{when $m\equiv 0\mmod{(p^b,t^r)}$,}\\
0,&\text{when $m\not\equiv 0\mmod{(p^b,t^r)}$.}\end{cases}
\end{equation}
In circumstances in which confusion is easily avoided, it may be expedient to suppress mention of the superscript 
$(r)$ associated with the implicit modulus $t^r$ from this notation. Indeed, for the remainder of this memoir, we shall 
take $r=2$ and drop mention of this index from our notation.

\section{The apparatus for nested efficient congruencing}
Before embarking on details of the nested efficient congruencing argument, we first assemble the apparatus required 
for our proof of Theorem \ref{theorem1.1}. This is inspired by our previous work on efficient congruencing and its 
nested variant (see, for example, \cite{Woo2012, Woo2015, Woo2019}). A number of modifications are necessary to 
accommodate the present application, so we will provide a complete account rather than cut corners in our 
exposition.\par

Let $k$ be an integer, and consider polynomials $\varphi_1(u),\ldots ,\varphi_k(u)\in \mathbb Z[t,u]$. We fix a prime 
number $p$ with $p>k$, and a large positive integer $B$. As in the discussion of \cite{Woo2019}, we require these 
polynomials $\varphi_j(u)$ to be sufficiently independent of one another for $1\le j\le k$. Thus, we impose the 
condition that the system of polynomials $\boldsymbol \varphi$ is {\it $p^c$-spaced} for an appropriate positive 
integer $c$, in the sense that
\begin{equation}\label{3.1}
\varphi_j(u)\equiv u^j\mmod{(p^c,t)}\quad (1\le j\le k).
\end{equation}
It follows that, given a $p^c$-spaced system of polynomials $\boldsymbol \varphi$, the Wronskian 
$W(u;\boldsymbol \varphi)$, defined via \eqref{1.2}, is non-zero modulo $t$. Indeed, one has
\[
W(u;\boldsymbol \varphi)\equiv \prod_{j=1}^kj!\not\equiv 0\mmod{(p,t)}.
\]

\par We next define the exponential sums and mean values that feature prominently in our discussion. We consider 
bounded complex weights $\mathfrak a_n$ with $n\in \mathfrak Z(t^2)$. When $q\in (0,\infty)$, we define the 
norms
\[
\|\boldsymbol{\mathfrak a}\|_q=
\lim_{T\rightarrow \infty}\Biggl( \sum_{|n_0|\le T}\sum_{|n_1|\le T}|\mathfrak a_{n_0+n_1t}|^q\Biggr)^{1/q}.
\]
We extend this definition in the conventional manner by putting
\[
\| \boldsymbol{\mathfrak a}\|_0=\lim_{q\rightarrow 0+}\|\boldsymbol{\mathfrak a}\|_q\quad \text{and}\quad 
\|\boldsymbol{\mathfrak a}\|_\infty =\lim_{q\rightarrow \infty}\|\boldsymbol{\mathfrak a}\|_q
=\sup_{n\in \mathfrak Z(t^2)}|\boldsymbol{\mathfrak a}_n|.
\]
We shall restrict attention to complex weights $\mathfrak a_n$ with $\|\boldsymbol{\mathfrak a}\|_0<\infty$. We 
note for future reference that we have in mind the situation in which the weights $\mathfrak a_n$ are supported 
only for $n=n_0+n_1t$, with $|n_0|\le N$ and $|n_1|\le N$, with roughly equal weight attached to each 
$n\in \mathfrak Z(t^2)$ thus prescribed. Next, when $h$ is a non-negative integer and $\xi\in \mathfrak Z(t^2)$, we 
introduce the subsequence of weights
\[
\boldsymbol{\mathfrak a}_h(\xi)=\left( \mathfrak a_n:\text{$n\equiv \xi\mmod{(p^h,t^2)}$}\right),
\]
and then define the normalising weights $\rho_h(\xi)=\rho_h(\xi;{\boldsymbol{\mathfrak a}})$ by putting
\begin{equation}\label{3.2}
\rho_h(\xi;{\boldsymbol{\mathfrak a}})=\|{\boldsymbol{\mathfrak a}_h(\xi)\|_0^{1/2}\cdot 
\|{\boldsymbol{\mathfrak a}}_h(\xi)}\|_\infty .
\end{equation}
It is convenient to associate the restricted subsequence of weights $\boldsymbol{\mathfrak b}$ to 
$\boldsymbol{\mathfrak a}_h(\xi)$ by putting
\[
\mathfrak b_n=\mathfrak a_{p^hn+\xi}\quad (n\in \mathfrak Z(t^2)).
\]
Temporarily, fix a choice of $\xi\in \mathfrak Z(t^2)$. Then we note that, in the scenario outlined above, there are 
roughly $(2N+1)^2$ non-zero values of $\mathfrak a_n$ each of size approximately 
$\|{\boldsymbol{{\boldsymbol{\mathfrak a}}}}\|_\infty$, so that $\rho_h(\xi;{\boldsymbol{\mathfrak a}})$ has size 
corresponding in rough terms to square-root cancellation amongst the weights $\mathfrak a_n$ constrained 
by the congruence condition $n\equiv \xi\mmod{(p^h,t^2)}$.\par

We impose an important hypothesis on the distribution of the weights $\mathfrak a_n$ for $n\in \mathfrak Z(t^2)$, 
which, fortunately, is easily checked in the situations of greatest interest for us.

\begin{hypothesis}\label{hypothesis3.1}
There exists a large positive integer $B_0$ and a small positive number $\omega$ for which the collection of weights 
$(\boldsymbol{\mathfrak a})$ satisfies the following properties:
\begin{enumerate}
\item[(i)] one has $p^{2B_0-2}<\|\boldsymbol{\mathfrak a}\|_0\le p^{2B_0}$;
\item[(ii)] when $\nu$ is a non-negative integer with $\nu \le \omega B_0$ and $\zeta_0\in \mathbb Z$, then for 
each $\varepsilon>0$, one has
\[
\sum_{\eta=1}^{p^\nu}\| {\boldsymbol{\mathfrak a}}_\nu(\zeta_0+\eta t)\|_0\ll (p^\nu)^{\varepsilon-1}
\sum_{\zeta=1}^{p^\nu}\sum_{\eta=1}^{p^\nu}\| {\boldsymbol{\mathfrak a}}_\nu(\zeta+\eta t)\|_0.
\]
\end{enumerate}
\end{hypothesis}

We shall need to make use of this hypothesis also when the sequence $(\boldsymbol{\mathfrak a})$ is restricted to 
particular congruence classes modulo $p^h$ that arise from our inductive procedures. We therefore extend 
Hypothesis \ref{hypothesis3.1} to an analogue with uniformity across such restrictions.

\begin{hypothesis}\label{hypothesis3.2}
Suppose that $B_1$ is a large positive integer having the property that the collection of weights 
$(\boldsymbol{\mathfrak a})$ satisfies the condition $p^{2B_1-2}<\| \boldsymbol{\mathfrak a}\|_0\le p^{2B_1}$. Then, 
there is a positive number $\delta$ having the property that whenever $h$ is an integer with 
$0\le h\le (1-\delta )B_1$ and $\xi\in \mathfrak Z(t^2)$, then the collection of weights $(\boldsymbol{\mathfrak b})=
(\boldsymbol{\mathfrak b}_h(\xi))$ satisfies Hypothesis \ref{hypothesis3.1}.
\end{hypothesis}

Observe that Hypothesis \ref{hypothesis3.2} has a nested property. Thus, if the collection of weights 
$(\boldsymbol{\mathfrak a})$ satisfies Hypothesis \ref{hypothesis3.2}, then the weights 
$(\boldsymbol{\mathfrak b})=(\boldsymbol{\mathfrak b}_h(\xi))$ occurring in the definition must also satisfy 
Hypothesis \ref{hypothesis3.2}, though with adjusted choices for $B_1$ and $\delta$. We shall make use of this nested 
property, with a large but finite number of nested levels, in our inductive argument.\par

These hypotheses assert that the support for the weights $\mathfrak a_n$ is not overly concentrated on particular 
residue classes $\zeta_0\mmod{(p^\nu,t)}$, for $1\le \nu \le \lfloor \omega B\rfloor$. We record for future reference 
that the conclusion of this hypothesis may be reintrepreted in the form
\[
\sum_{\substack{\xi\mmod{(p^\nu,t^2)}\\ \xi\equiv \zeta_0\mmod{(p^\nu, t)}}}\rho_\nu(\xi)^2\ll p^{\varepsilon-\nu}
\rho_0(1)^2.
\]
Consider temporarily the scenario in which $p^B\le N$ and
\[
\mathfrak a_n=\begin{cases} 1,&\text{when $n=n_0+n_1t$, with $|n_0|\le N$ and $|n_1|\le N$,}\\
0,&\text{otherwise.}\end{cases}
\]
Then one sees that Hypothesis \ref{hypothesis3.1} is transparently satisfied, for we have
\[
\sum_{\substack{\xi\mmod{(p^\nu,t^2)}\\ \xi\equiv \zeta_0\mmod{(p^\nu, t)}}}\rho_\nu(\xi)^2=
\sum_{\eta=1}^{p^\nu}\Bigl( \frac{2N}{p^\nu}+O(1)\Bigr)^2\ll p^{-\nu}(2N)^2,
\]
whilst
\[
p^{-\nu}\rho_0(1)^2=p^{-\nu}(2N+O(1))^2.
\]
A similar thought experiment applies in similar fashion for Hypothesis \ref{hypothesis3.2}. It is these and closely 
allied scenarios that will most likely find application, and so the Hypothesis \ref{hypothesis3.2} does not restrict our 
freedom of manoeuvre. The inclusion of the weights $\rho_h(\xi)$ assists us in correctly pursuing the nested 
iterative procedure as one investigates congruence classes modulo various powers of $p$. Thus, even though one 
could avoid Hypotheses \ref{hypothesis3.1} and \ref{hypothesis3.2} for the situations of greatest interest to us, the 
extra flexibility added would seem to be of value.\par

Next, when $n\in \mathfrak Z(t^2)$ and $\alpha_1,\ldots ,\alpha_k\in \mathfrak T(t^2)$, we write
\begin{equation}\label{3.3}
\psi(n;\boldsymbol \alpha)=\alpha_1\varphi_1(n)+\ldots +\alpha_k\varphi_k(n).
\end{equation}
When $h$ is a non-negative integer and $\xi \in \mathfrak Z(t^2)$, we define the exponential sum 
\[
\mathfrak f_h(\boldsymbol \alpha;\xi)=\mathfrak f_h(\boldsymbol \alpha ;\xi;\boldsymbol{\mathfrak a};
\boldsymbol \varphi)
\]
as follows. When $\|\boldsymbol{\mathfrak a}\|_\infty >0$, we put
\begin{equation}\label{3.4}
\mathfrak f_h(\boldsymbol \alpha ;\xi)=\rho_h(\xi)^{-1}\sum_{n\equiv \xi\mmod{(p^h,t^2)}}
\mathfrak a_nE(\psi(n;\boldsymbol \alpha)),
\end{equation}
where, here and elsewhere, the summation over $n$ is implicitly taken over all $n\in \mathfrak Z(t^2)$ constrained 
by the explicit congruence condition. For completeness, when $\|\boldsymbol{\mathfrak a}\|_\infty=0$, we instead 
put $\mathfrak f_h(\boldsymbol \alpha ;\xi)=0$. For brevity, we write $f_{\mathfrak a}(\boldsymbol \alpha)$ for 
$\mathfrak f_0(\boldsymbol \alpha;\xi)$. Note that in the scenario outlined in the previous paragraph, one has
\[
(2N+1)f_{\mathfrak a}(\boldsymbol \alpha)=\sum_{|n_0|\le N}\sum_{|n_1|\le N}E\bigl(\alpha_1\varphi_1(n_0+n_1t)+
\ldots +\alpha_k\varphi_k(n_0+n_1t)\bigr).
\]

\par Specialising the notation from \S2 to the case $r=2$, when $B$ is a positive integer and 
$F:\mathfrak T(t^2)^k\rightarrow \mathbb C$, we write
\begin{equation}\label{3.5}
\oint_{p^B}F(\boldsymbol \alpha)\d\boldsymbol \alpha =(p^2)^{-kB}\sum_{u_1\mmod{(p^B,t^2)}}\ldots 
\sum_{u_k\mmod{(p^B,t^2)}}F\biggl( \frac{u_1}{p^B},\ldots ,\frac{u_k}{p^B}\biggr) ,
\end{equation}
where the summations are taken over complete sets of elements of $\mathfrak Z(t^2)$ modulo $p^B$. Thus, for 
example, the summation over $u_1$ may be construed as running over all elements $v+wt$ with $1\le v,w\le p^B$. 
We then define the mean value 
$U_{s,k}^B(\boldsymbol{\mathfrak a})=U_{s,k}^{B,\boldsymbol \varphi}(\boldsymbol{\mathfrak a})$ by putting
\begin{equation}\label{3.6}
U_{s,k}^B(\boldsymbol{\mathfrak a})=\oint_{p^B}|f_{\boldsymbol{\mathfrak a}}(\boldsymbol \alpha)|^{2s}
\d\boldsymbol \alpha .
\end{equation}
By orthogonality, the mean value $U_{s,k}^B(\boldsymbol{\mathfrak a})$ counts the integral solutions of the 
simultaneous congruences
\begin{equation}\label{3.7}
\sum_{i=1}^s(\varphi_j(x_i)-\varphi_j(y_i))\equiv 0\mmod{(p^B,t^2)}\quad (1\le j\le k),
\end{equation}
with $\mathbf x,\mathbf y\in \mathfrak Z(t^2)$, and with each solution $\mathbf x,\mathbf y$ being counted with 
weight
\[
\rho_0(1)^{-2s}\prod_{i=1}^s \mathfrak a_{x_i}{\overline{\mathfrak a}}_{y_i}.
\]

\par A variant of $U_{s,k}^B(\boldsymbol{\mathfrak a})$, in which the underlying variables are restricted to a common 
congruence class $\xi$ modulo $(p^h,t^2)$, is defined by putting $U_{s,k}^{B,h}(\boldsymbol{\mathfrak a})=
U_{s,k}^{B,h,\boldsymbol \varphi}(\boldsymbol{\mathfrak a})$, where
\begin{equation}\label{3.8}
U_{s,k}^{B,h,\boldsymbol \varphi}(\boldsymbol{\mathfrak a})=\rho_0(1)^{-2}\sum_{\xi\mmod{(p^h,t^2)}}
\rho_h(\xi)^2\oint_{p^B}|\mathfrak f_h(\boldsymbol \alpha;\xi)|^{2s}\d\boldsymbol \alpha .
\end{equation}
By orthogonality, once again, the integral on the right hand side of \eqref{3.8} counts the solutions 
$\mathbf x,\mathbf y\in \mathfrak Z(t^2)^s$ of the system \eqref{3.7} satisfying 
$\mathbf x\equiv \mathbf y\equiv \xi\mmod{(p^h,t^2)}$ with weight
\[
\rho_h(\xi)^{-2s}\prod_{i=1}^s \mathfrak a_{x_i}{\overline{\mathfrak a}}_{y_i}.
\]

\par We briefly discuss a trivial estimate for the mean value $U_{s,k}^{B,h}(\boldsymbol{\mathfrak a})$ in order to 
motivate the analysis to come. Our starting point is a decomposition of the exponential sum 
$f_{\boldsymbol{\mathfrak a}}(\boldsymbol \alpha)$ into the residue classes modulo $(p^H,t^2)$ of the summands. 
Thus, from \eqref{3.4}, we find that
\begin{align*}
f_{\boldsymbol{\mathfrak a}}(\boldsymbol \alpha)&
=\rho_0(1)^{-1}\sum_{\xi\mmod{(p^H,t^2)}}\sum_{n\equiv \xi\mmod{(p^H,t^2)}}
\mathfrak a_nE(\psi(n;\boldsymbol \alpha))\\
&=\rho_0(1)^{-1}\sum_{\xi\mmod{(p^H,t^2)}}\rho_H(\xi)\mathfrak f_H(\boldsymbol \alpha ;\xi).
\end{align*}
Since
\[
\sum_{\xi\mmod{(p^H,t^2)}}\rho_H(\xi)^2=\sum_{\xi\mmod{(p^H,t^2)}}\|\boldsymbol{\mathfrak a}_H(\xi)\|_0\cdot 
\| \boldsymbol{\mathfrak a}_H(\xi)\|_\infty^2\le \|\boldsymbol{\mathfrak a}\|_0\cdot \|
\boldsymbol{\mathfrak a}\|_\infty^2=\rho_0(1)^2,
\]
an application of H\"older's inequality reveals that $|f_{\boldsymbol{\mathfrak a}}(\boldsymbol \alpha)|^{2s}$ is 
bounded above by
\begin{align}\rho_0(1)^{-2s}&
\Biggl( \,\sum_{\xi\mmod{(p^H,t^2)}}1\Biggr)^{s}\Biggl( \,\sum_{\xi\mmod{(p^H,t^2)}}\rho_H(\xi)^2\Biggr)^{s-1}
\sum_{\xi\mmod{(p^H,t^2)}}\rho_H(\xi)^2|\mathfrak f_H(\boldsymbol \alpha;\xi)|^{2s}
\notag \\
&\ll (p^{2H})^s\rho_0(1)^{-2}\sum_{\xi\mmod{(p^H,t^2)}}\rho_H(\xi)^2|\mathfrak f_H(\boldsymbol \alpha;\xi)|^{2s}.
\label{3.9}
\end{align}
Thus, on recalling \eqref{3.6} and \eqref{3.8}, we infer that
\begin{equation}\label{3.10}
U_{s,k}^B(\boldsymbol{\mathfrak a})\le (p^{2H})^sU_{s,k}^{B,H}(\boldsymbol{\mathfrak a}).
\end{equation}
As in the analogous discussion of \cite{Woo2019}, we are consequently justified in seeking an exponent $\lambda$ 
having the property that, for each $\varepsilon>0$, one has
\begin{equation}\label{3.11}
U_{s,k}^B(\boldsymbol{\mathfrak a})\ll (p^{2H})^{\lambda+\varepsilon}U_{s,k}^{B,H}(\boldsymbol{\mathfrak a}),
\end{equation}
with as much uniformity as is feasible to establish in the various parameters. A comparison of \eqref{3.10} and 
\eqref{3.11} makes it apparent already that such an exponent $\lambda$ exists with $\lambda \le s$.\par

We now discuss a number of helpful simplifications. We observe first that the definition \eqref{3.4} of 
$\mathfrak f_h(\boldsymbol \alpha ;\xi)$ ensures that we have scale invariance with repect to the coefficient weights 
$\boldsymbol{\mathfrak a}$. Thus, if $\gamma>0$ and $\boldsymbol{\mathfrak a}$ is replaced by 
$\gamma \boldsymbol{\mathfrak a}$, then since 
$\rho_h(\gamma \boldsymbol{\mathfrak a})=\gamma \rho_h(\boldsymbol{\mathfrak a})$, the exponential sum 
$\mathfrak f_h(\boldsymbol \alpha ;\xi)$ remains unchanged. The exponential sum 
$f_{\boldsymbol{\mathfrak a}}(\boldsymbol \alpha)$ and the mean values $U_{s,k}^B(\boldsymbol{\mathfrak a})$ 
and $U_{s,k}^{B,h}(\boldsymbol{\mathfrak a})$ all inherit this scale invariance property. Denote by $\mathbb D(B)$ 
the set of sequences $(\mathfrak a_n)_{n\in \mathfrak Z(t^2)}$ satisfying Hypothesis \ref{hypothesis3.2} with 
$|\mathfrak a_n|\le 1$ $(n\in \mathfrak Z(t^2))$ and
\[
p^{2B-2}<\| \boldsymbol{\mathfrak a}\|_0<p^{2B}.
\]
Also, write $\mathbb D_0(B)=\mathbb D(B)\cup \{\mathbf 0\}$. Then, in view of the scale invariance just discussed, 
we are at liberty to restrict our deliberations to sequences $(\mathfrak a_n)$ lying in $\mathbb D(B)$.\par

We next provide a precise definition of a real number tantamount to the exponent $\lambda$ discussed above. 
When $\tau>0$, we denote by $\Phi_\tau(B)$ the set of $p^c$-spaced $k$-tuples of polynomials 
$\boldsymbol \varphi(u)$, with $c\ge \tau B$. The upper bound \eqref{3.10} ensures that for each fixed 
$B\in \mathbb N$ and $\boldsymbol \varphi(u)\in \Phi_\tau(B)$, and for every non-negative integer $H$, one has
\[
\sup_{(\boldsymbol{\mathfrak a})\in \mathbb D(B)} 
\frac{\log \bigl( U_{s,k}^B(\boldsymbol{\mathfrak a})/U_{s,k}^{B,H}(\boldsymbol{\mathfrak a})\bigr)}{\log (p^{2H})}
\le s.
\]
Now, given $s>0$, $\theta\ge 1$ and $\tau>0$, we write $H=\lceil B/\theta\rceil$, and define
\begin{equation}\label{3.12}
\lambda^*(s,\theta;\tau)=\limsup_{B\rightarrow \infty}\sup_{\boldsymbol \varphi\in \Phi_\tau(B)}\sup_{(\boldsymbol 
{\mathfrak a})\in \mathbb D(B)}\frac{\log \bigl( U_{s,k}^B(\boldsymbol{\mathfrak a})/U_{s,k}^{B,H}
(\boldsymbol{\mathfrak a})\bigr)}{\log (p^{2H})}.
\end{equation}
We shall be focused on the limiting exponent
\begin{equation}\label{3.13}
\lambda (s,\theta)=\limsup_{\tau\rightarrow 0+}\lambda^*(s,\theta;\tau).
\end{equation}
Since it is apparent from the discussion above that $\lambda^*(s,\theta;\tau)\le s$, we deduce that
\begin{equation}\label{3.14}
\lambda(s,\theta)\le s.
\end{equation}

The main conclusion of this memoir is described in the following theorem.

\begin{theorem}\label{theorem3.2} Suppose that $k\in \mathbb N$ and that $p$ is a prime with $p>k$. Then one 
has $\lambda (k(k+1)/2,k)\le 0$.
\end{theorem}

The utility of this theorem is made transparent in the following corollary.

\begin{corollary}\label{corollary3.3} Suppose that $k\in \mathbb N$ and that $p$ is a prime number with $p>k$. In 
addition, let $\tau>0$ and $\varepsilon>0$. Finally, let $B$ be sufficiently large in terms of $k$, $\tau$ and 
$\varepsilon$. Put $s=k(k+1)/2$ and $H=\lceil B/k\rceil$. Then for every $\boldsymbol \varphi\in \Phi_\tau(B)$, and 
every sequence $(\boldsymbol{\mathfrak a})\in \mathbb D_0(B)$, one has
\begin{equation}\label{3.15}
U_{s,k}^B(\boldsymbol{\mathfrak a})\ll (p^B)^\varepsilon U_{s,k}^{B,H}(\boldsymbol{\mathfrak a}).
\end{equation}
\end{corollary}

It is apparent from the relation \eqref{3.15} that in a certain sense, one can assume an additional constraint within the 
underlying mean value \eqref{3.6} that effectively concentrates the underlying variables towards the same congruence 
class modulo $(p^H,t^2)$. In the model situation in which $\mathfrak a_n$ is supported only when $n=n_0+n_1t$ with 
$|n_0|\le N$ and $|n_1|\le N$, we shall see that when $B$ is large enough that $p^{B}>N^k$, then one has 
$U_{s,k}^{B,H}(\boldsymbol{\mathfrak a})\ll p^{B\varepsilon}$, and hence also 
$U_{s,k}^B(\boldsymbol{\mathfrak a})\ll p^{2B\varepsilon}$.\par

Our next step is to introduce auxiliary mean values that form the important cogs in our nested iterative machine. 
Throughout, we fix $s=k(k+1)/2$. Let $a$, $b$, $c$ and $\nu$ be non-negative integers. We consider a fixed 
$p^c$-spaced $k$-tuple of polynomials $\boldsymbol \varphi\in \mathbb Z[t,u]^k$, and we fix a complex 
sequence $(\boldsymbol{\mathfrak a})\in \mathbb D_0(B)$. When $0\le r\le k$, we define the mean value 
$K_{a,b}^r=K_{a,b,c}^{r,\boldsymbol \varphi,\nu}(\boldsymbol{\mathfrak a})$ by putting
\begin{equation}\label{3.16}
K_{a,b,c}^{r,\boldsymbol \varphi,\nu}(\boldsymbol{\mathfrak a})=\rho_0(1)^{-4}\sum_{\xi\mmod{(p^a,t^2)}}
\sum_{\substack{\eta\mmod{(p^b,t^2)}\\ \eta\not\equiv \xi\mmod{(p^\nu,t)}}}\rho_a(\xi)^2\rho_b(\eta)^2
K_{a,b,c}^{r,\boldsymbol \varphi,\nu}(\boldsymbol{\mathfrak a};\xi,\eta),
\end{equation}
in which
\begin{equation}\label{3.17}
K_{a,b,c}^{r,\boldsymbol \varphi,\nu}(\boldsymbol{\mathfrak a};\xi,\eta)=
\oint_{p^B}|\mathfrak f_a(\boldsymbol \alpha;\xi)^{2R}\mathfrak f_b(\boldsymbol \alpha;\eta)^{2s-2R}|\d\boldsymbol 
\alpha
\end{equation}
and $R=r(r+1)/2$. The inclusion of the superscript $\nu$ in this notation is intended to remind us of the assumption 
in the second summation of \eqref{3.16} that $p^\nu\nmid (\eta_0-\xi_0)$, in which we have written 
$\xi=\xi_0+\xi_1t$ and $\eta=\eta_0+\eta_1t$, with $\xi_i,\eta_i\in \mathbb Z$. The reader may care to note, for 
future reference, that this condition ensures that $\eta-\xi$ is approximately a unit in the ring 
$\mathbb Z[t]/(p^B,t^2)$. For future reference, we note that as a consequence of orthogonality, the mean value 
$K_{a,b,c}^{r,\boldsymbol \varphi,\nu}(\boldsymbol{\mathfrak a};\xi,\eta)$ counts the solutions in $\mathfrak Z(t^2)$ 
of the simultaneous congruences
\begin{equation}\label{3.18}
\sum_{i=1}^R \bigl( \varphi_j(x_i)-\varphi_j(y_i)\bigr) \equiv \sum_{l=1}^{s-R}\bigl( \varphi_j(v_l)-\varphi_j(w_l)\bigr) 
\mmod{(p^B,t^2)}\quad (1\le j\le k),
\end{equation} 
satisfying
\begin{equation}\label{3.19}
\mathbf x\equiv \mathbf y\equiv \xi\mmod{(p^a,t^2)}\quad \text{and}\quad \mathbf v\equiv \mathbf w\equiv \eta 
\mmod{(p^b,t^2)},
\end{equation}
with each solution being counted with weight
\begin{equation}\label{3.20}
\rho_a(\xi)^{-2R}\rho_b(\eta)^{2R-2s}\biggl( \prod_{i=1}^R \mathfrak a_{x_i}{\overline{\mathfrak a}}_{y_i}\biggr) 
\biggl( \prod_{l=1}^{s-R} \mathfrak a_{v_l}{\overline{\mathfrak a}}_{w_l}\biggr) .
\end{equation}
We note that the relations \eqref{3.18}, \eqref{3.19} and \eqref{3.20} require careful interpretation. Indeed, each 
variable $x_i,y_i,v_l,w_l$ is an element of $\mathfrak Z(t^2)$. Thus, for example, one should think of $x_i$ as being of 
the shape $u+vt$, with $u,v\in \mathbb Z$. Likewise, one should think of $\xi$ and $\eta$ as lying in 
$\mathfrak Z(t^2)$, with similar associated comments.\par

As in our previous work \cite{Woo2019} on nested efficient congruencing, we find it useful to introduce notation 
making the size of a mean value relative to its conjectured size more readily transparent. Following this previous 
path, when $1\le r\le k-1$ and $\Delta$ is a positive number, we define
\begin{equation}\label{3.21}
[[ K_{a,b,c}^{r,\boldsymbol \varphi,\nu}(\boldsymbol{\mathfrak a})]]_\Delta 
=\Biggl( \frac{K_{a,b,c}^{r,\boldsymbol \varphi,\nu}(\boldsymbol{\mathfrak a})}
{(p^{2H})^\Delta U_{s,k}^{B,H}(\boldsymbol{\mathfrak a})}\Biggr)^{\textstyle{\frac{k-1}{r(k-r)}}}.
\end{equation}
Much of our discussion in the coming sections is largely independent of the parameters $\nu$ and $c$, the 
sequence $\boldsymbol{\mathfrak a}$, and the polynomial system $\boldsymbol \varphi$. Hence, when ambiguity is 
easily avoided, we shall feel free to omit explicit mention of such notational decorations. We note also that the 
exponent 
\[
\frac{k-1}{r(k-r)},
\]
as in the analogous treatment of \cite{Woo2019}, serves to streamline the nested iteration by enabling a simplified 
single branch version to be utilised in place of the more complicated multi-branch process immediately to hand. 

\section{Uniform estimates for translation-dilation invariant families}
As explained in \cite[\S4]{Woo2019}, a $p^c$-spaced system of polynomials $\boldsymbol \varphi$ may be thought 
of as approximating the special system in which $\varphi_j(u)=u^j$ $(1\le j\le k)$. The system of congruences
\begin{equation}\label{4.1}
\sum_{i=1}^s \varphi_j(x_i)\equiv \sum_{i=1}^s\varphi_j(y_i)\mmod{(p^B,t^2)}\quad (1\le j\le k)
\end{equation}
is then translation-dilation invariant. Thus, whenever $\mathbf x,\mathbf y$ is a solution of \eqref{4.1}, then so too is 
$q\mathbf x+a,q\mathbf y+a$, for any $q,a\in \mathbb Z[t]$. This property permits a certain scaling property to be 
established for the mean value $U_{s,k}^{B,h}(\boldsymbol{\mathfrak a})$. In advance of the announcement and 
proof of this scaling property, we pause to establish two basic invertibility properties of elements in 
$\mathbb Z[t]/(p^B,t^d)$ and associated matrices.

\begin{lemma}\label{lemma4.1}
Let $d$ be a natural number, and suppose that $n\in \mathbb Z[t]/(p^B,t^d)$. Let $n_0$ be an integer having the 
property that $n\equiv n_0\mmod{(p^B,t)}$, and suppose that $(n_0,p)=1$. Then $n$ possesses a multiplicative 
inverse $n^{-1}$ lying in $\mathbb Z[t]/(p^B,t^d)$.
\end{lemma}

\begin{proof} Since $(n_0,p)=1$, there exists an integer $m_0$ having the property that one has 
$m_0n_0\equiv 1\mmod{p^B}$. On observing that $m_0n-1\equiv 0\mmod{(p^B,t)}$, we see that for suitable 
elements $\chi, \omega\in \mathbb Z[t]/(p^B,t)$, one has
\[
m_0n=1+p^B\chi+t\omega.
\]
We now define the polynomial
\[
m=1-(p^B\chi+t\omega)+(p^B\chi+t\omega)^2-\ldots +(-1)^{d-1}(p^B\chi+t\omega)^{d-1}.
\]
One sees that
\[
m(m_0n)=1+(-1)^{d-1}(p^B\chi+t\omega)^d\equiv 1+(-1)^{d-1}(t\omega)^d\equiv 1\mmod{(p^B,t^d)}.
\]
Hence, the element $n$ has the multiplicative inverse $n^{-1}=m_0m$ modulo $(p^B,t^d)$.
\end{proof}

In the next lemma, and later, we denote the $k\times k$ identity matrix over the ring $\mathbb Z[t]/(p^B,t^d)$ by 
$I_k$.

\begin{lemma}\label{lemma4.2}
Let $d$ be a natural number, and suppose that $A=(a_{ij})_{1\le i,j\le k}$ is a matrix with entries in 
$\mathbb Z[t]/(p^B,t^d)$ having the property that $A\equiv I_k\mmod{(p,t)}$. Then $A$ has an inverse matrix 
$A^{-1}$ in $\mathbb Z[t]/(p^B,t^d)$ satisfying $A^{-1}\equiv I_k\mmod{(p,t)}$.
\end{lemma}

\begin{proof} By appealing to the reader's knowledge of standard linear algebra, the existence of $A^{-1}$ is 
assured provided one is able to apply elementary row operations to obtain $I_k$ from $A$. We first show that a 
sequence of elementary row operations may be employed to transform $A$ into lower triangular form $A^{(0)}$ 
with $A^{(0)}\equiv I_k\mmod{(p,t)}$. This we establish in the familiar inductive manner. We suppose at stage $r$ 
that $A$ has been transformed into a matrix $A^{(r)}$, with $A^{(r)}\equiv I_k\mmod{(p,t)}$, having the property that 
the columns with index exceeding $r$ are already in lower triangular form. This of course holds trivially when $r=k$ 
with $A^{(k)}=A$. We then seek to establish that the same conclusion can be obtained for all columns with index 
exceeding $r-1$. For this purpose, it suffices to observe that, in view of Lemma \ref{lemma4.1}, the diagonal entry 
$a^{(r)}_{rr}$ possesses a multiplicative inverse in $\mathbb Z[t]/(p^B,t^d)$, an element we label as 
$(a^{(r)}_{rr})^{-1}$. We subtract $a^{(r)}_{ir}(a^{(r)}_{rr})^{-1}$ times the $r$-th row from the $i$-th row, for 
$1\le i<r$. In this way, we obtain a new matrix $A^{(r-1)}$ having lower triangular form amongst the columns with 
index exceeding $r-1$. Since $a^{(r)}_{ir}\equiv 0\mmod{(p,t)}$ for $1\le i<r$, moreover, we see that this new matrix 
remains congruent to $I_k$ modulo $(p,t)$. This completes the inductive step, so that we may now work with the 
lower triangular matrix $A^{(0)}$.\par

The next step is to convert the matrix $A^{(0)}$ from lower triangular form to diagonal form $D^{(k)}$ via elementary 
row operations, again with $D^{(k)}\equiv I_k\mmod{(p,t)}$. This we achieve by reversing course, again applying 
induction. We suppose at stage $r$ that $A^{(0)}$ has been transformed into a matrix $D^{(r-1)}$, with 
$D^{(r-1)}\equiv I_k\mmod{(p,t)}$, having the property that the columns with index smaller than $r$ are already in 
diagonal form. This holds trivially when $r=1$ with $D^{(0)}=A^{(0)}$. We then seek to establish that the same 
conclusion can be obtained for all columns with index smaller than $r+1$. In this instance, we subtract 
$d^{(r)}_{ir}(d^{(r)}_{rr})^{-1}$ times the $r$-th row from the $i$-th row, for $r<i\le k$. In this way, we obtain a new 
matrix $D^{(r)}$ having diagonal form amongst the columns with index smaller than $r+1$. Since 
$d^{(r)}_{ir}\equiv 0\mmod{(p,t)}$ for $r<i\le k$, moreover, we again see that this new matrix remains congruent to 
$I_k$ modulo $(p,t)$. This completes the inductive step in this second phase, and so we may now work with the 
diagonal matrix $D^{(k)}$.\par

Finally, by multiplying the $r$-th row through by $(d_{rr}^{(k)})^{-1}$ for $1\le r\le k$, we transform the diagonal 
matrix $D^{(k)}$ into $I_k$. We have thus shown that the matrix $A$ may be converted to $I_k$ through elementary 
row operations defined over $\mathbb Z[t]/(p^B,t^d)$, and thus $A$ possesses a multiplicative inverse $A^{-1}$ 
having entries in the same ring. Finally, since $A\equiv I_k\mmod{(p,t)}$, we must have 
$A^{-1}\equiv I_k\mmod{(p,t)}$. This completes the proof of the lemma.
\end{proof}

It is useful also to have an extension of this result applicable to triangular matrices.

\begin{lemma}\label{lemma4.3} Let $d$ be a natural number, and suppose that $A=(a_{ij})_{1\le i,j\le k}$ is a matrix 
with entries in $\mathbb Z[t]/(p^B,t^d)$ having the property that
\[
a_{ij}\equiv \begin{cases}0\mmod{(p,t)},&\text{when $i>j$,}\\
1\mmod{(p,t)},&\text{when $i=j$.}\end{cases}
\]
Then $A$ has an inverse matrix $A^{-1}$ in $\mathbb Z[t]/(p^B,t^d)$.
\end{lemma}

\begin{proof} We adopt the inductive approach of the proof of the previous lemma, applying elementary row 
operations to transform the matrix $A$ into a new matrix $C$ in which $C\equiv I_k\mmod{(p,t)}$. The matrix $A$ 
has an inverse matrix $A^{-1}$, with entries in $\mathbb Z[t]/(p^B,t)$, if and only if $C$ is invertible in the latter 
ring. However, the invertibility of $C$ is guaranteed by the conclusion of Lemma \ref{lemma4.2}, and thus the 
desired conclusion follows. This completes the proof of the lemma.
\end{proof}

With these preliminaries in hand, we are now equipped to establish the main conclusion of this section.

\begin{lemma}\label{lemma4.4}
Suppose that $k\in \mathbb N$ and that $p$ is a prime number with $p>k$. Let $\tau$, $\varepsilon$ and $\delta$ 
be positive numbers with $\varepsilon<\tau<\delta<1$, and let $B$ be sufficiently large in terms of $s$, $k$ and 
$\varepsilon$. Write $H=\lceil B/k\rceil$. Then for every $\boldsymbol \varphi\in \Phi_\tau(B)$, every sequence 
$(\boldsymbol{\mathfrak a})\in \mathbb D_0(B)$, and every non-negative integer $h$ with $h\le (1-\delta)H$, one has
\[
U_{s,k}^{B,h}(\boldsymbol{\mathfrak a})\ll (p^{2(H-h)})^{\lambda (s,k)+\varepsilon}U_{s,k}^{B,H}(\boldsymbol 
{\mathfrak a}).
\]
\end{lemma}

\begin{proof} The argument required in our proof is a mild adjustment of that applied to establish 
\cite[Lemma 4.1]{Woo2019}, with a small number of critical adjustments. Since the setting is somewhat alien, 
however, at the risk of being seemingly repetitive in places, we provide fairly complete details.\par

Consider a $p^c$-spaced $k$-tuple of polynomials $\boldsymbol \varphi$ with $c\ge \tau B$, and a complex 
sequence $(\boldsymbol{\mathfrak a})\in \mathbb D_0(B)$, with associated parameters satisfying the statement of 
the lemma. As in the proof of \cite[Lemma 4.1]{Woo2019}, we let $h$ be an integer with $0\le h\le (1-\delta )H$. 
Then it follows from \eqref{3.4} that
\begin{equation}\label{4.2}
\mathfrak f_h(\boldsymbol \alpha ;\xi)=\rho_h(\xi)^{-1}\sum_{y\in \mathfrak Z(t^2)}\mathfrak b_y(\xi)
E(\psi(p^hy+\xi;\boldsymbol \alpha)),
\end{equation}
with $\psi(n;\boldsymbol \alpha)$ defined by \eqref{3.3}, and where
\begin{equation}\label{4.3}
\mathfrak b_y(\xi)=\mathfrak a_{p^hy+\xi}\quad (y\in \mathfrak Z(t^2)).
\end{equation}
An appeal to orthogonality reveals that the mean value on the right hand side of \eqref{3.8} counts the solutions 
$\mathbf y, \mathbf z\in \mathfrak Z(t^2)$ of the simultaneous congruences
\begin{equation}\label{4.4}
\sum_{i=1}^s \varphi_j(p^hy_i+\xi)\equiv \sum_{i=1}^s \varphi_j(p^hz_i+\xi)\mmod{(p^B,t^2)}\quad (1\le j\le k),
\end{equation}
with each solution being counted with weight
\begin{equation}\label{4.5}
\rho_h(\xi)^{-2s}\prod_{i=1}^s{\mathfrak b}_{y_i}{\overline{\mathfrak b}}_{z_i}.
\end{equation}

\par Since the system $\boldsymbol \varphi$ is $p^c$-spaced, we find from \eqref{3.1} that for suitable polynomials 
$\psi_j(u),\chi_j(u)\in \mathbb Z[t,u]$, one may write
\[
\varphi_j(u)=\sum_{i=1}^k\omega_{ij}u^i+u^{k+1}\bigl( p^c\psi_j(u)+t\chi_j(u)\bigr) \quad (1\le j\le k),
\]
for some coefficient matrix $A_1=(\omega_{ij})_{1\le i,j\le k}$ congruent modulo $(p^c,t)$ to the $k\times k$ identity 
matrix $I_k$. Since $A_1\equiv I_k\mmod{(p^c,t)}$, it follows from Lemma \ref{lemma4.2} that $A_1$ possesses a 
multiplicative inverse $A_1^{-1}$ modulo $(p^B,t^2)$, having coefficients lying in $\mathbb Z[t]/(p^B,t^2)$, and 
satisfying $A_1^{-1}\equiv I_k\mmod{(p,t)}$. We replace $\boldsymbol \varphi$ by $A_1^{-1}\boldsymbol \varphi$, 
likewise $\boldsymbol \psi$ by $A_1^{-1}\boldsymbol \psi$, and also $\boldsymbol \chi$ by 
$A_1^{-1}\boldsymbol \chi$. In so doing, we are effectively taking suitable integral linear combinations of the 
congruences comprising \eqref{4.4}. In this way, we see that there is no loss of generality in supposing that for 
$1\le j\le k$, one has
\begin{equation}\label{4.6}
\varphi_j(u)=u^j+u^{k+1}\bigl( p^c\psi_j(u)+t\chi_j(u)\bigr) \quad (1\le j\le k).
\end{equation}

\par The polynomials $\varphi_j(u)$ defined as in \eqref{4.6} are now in a convenient form with which to apply the 
binomial theorem. Thus, for suitable polynomials $\Psi_j(u), {\rm X}_j(u)\in \mathbb Z[t,u]$, one has
\[
\varphi_j(p^hy+\xi)-\varphi_j(\xi)=\varPhi_j(p^hy),
\]
in which
\[
\Phi_j(w)=\sum_{i=1}^k\Omega_{ij}w^i+w^{k+1}\bigl( p^c\Psi_j(w)+t{\rm X}_j(w)\bigr) 
\quad (1\le j\le k),
\]
and the coefficients $\Omega_{ij}\in \mathbb Z[t]$ satisfy
\[
\Omega_{ij}\equiv \binom{j}{i}\xi^{j-i}\mmod{(p^c,t)}\quad (1\le i,j\le k).
\]
As usual here, we adopt the convention that the binomial coefficient $\binom{j}{i}$ is equal to $0$ for $i>j$. The matrix 
$A_2=(\Omega_{ij})_{1\le i,j\le k}$ satisfies the property that
\[
\Omega_{ij}\equiv \begin{cases}0\mmod{(p,t)},&\text{when $i>j$,}\\
1\mmod{(p,t)},&\text{when $i=j$.}\end{cases}
\]
Since we assume that $p>k$, we find from Lemma \ref{lemma4.3} that the matrix $A_2$ possesses a multiplicative 
inverse $A_2^{-1}$ in $\mathbb Z[t]/(p^B,t^2)$. We now replace $\boldsymbol \Phi$ by 
$A_2^{-1}\boldsymbol \Phi$, likewise $\boldsymbol \Psi$ by $A_2^{-1}\boldsymbol \Psi$, and also 
$\boldsymbol{\rm X}$ by $A_2^{-1}\boldsymbol{\rm X}$. This transformation again corresponds to taking linear 
combinations of the congruences comprising \eqref{4.4}. In this way, we see that there is no loss of generality in 
supposing that the coefficient matrix $A_2$ is equal to $I_k$. Consequently, there exist polynomials 
$\Upsilon_j(u),\Theta_j(u)\in \mathbb Z[t,u]$ having the property that whenever the system \eqref{4.4} is satisfied, 
then
\begin{equation}\label{4.7}
\sum_{i=1}^s (p^h)^j\bigl( \Phi_j(y_i)-\Phi_j(z_i)\bigr) \equiv 0\mmod{(p^B,t^2)}\quad (1\le j\le k),
\end{equation}
in which $\Phi_j(u)$ has been redefined by putting
\begin{equation}\label{4.8}
\Phi_j(u)=u^j+p^{h}u^{k+1}\bigl(p^c\Upsilon_j(u)+t\Theta_j(u)\bigr).
\end{equation} 
At the conclusion of this first phase of the argument, we have shown that the solutions 
$\mathbf y,\mathbf z\in \mathfrak Z(t^2)$ of the system of congruences \eqref{4.4}, counted with weight \eqref{4.5} 
associated with the definition \eqref{3.8} of $U_{s,k}^{B,h}(\boldsymbol{\mathfrak a})$, are constrained by an 
additional system of congruences. Thus, we have
\begin{equation}\label{4.9}
\sum_{i=1}^s \Phi_j(y_i)\equiv \sum_{i=1}^s \Phi_j(z_i)\mmod{(p^{B-kh},t^2)}\quad (1\le j\le k).
\end{equation}

\par The condition \eqref{4.9} can be incorporated into the mean value $U_{s,k}^{B,h}(\boldsymbol{\mathfrak a})$. 
When $\boldsymbol \beta\in \mathfrak T(t^2)^k$, we define
\begin{equation}\label{4.10}
\mathfrak g_h(\boldsymbol \alpha,\boldsymbol \beta;\xi)=\rho_h(\xi)^{-1}\sum_{y\in \mathfrak Z(t^2)}
\mathfrak c_y(\xi;\boldsymbol \alpha)E(\psi^*(y;\boldsymbol \beta)),
\end{equation}
where
\begin{equation}\label{4.11}
\mathfrak c_y(\xi;\boldsymbol \alpha)=\mathfrak b_y(\xi)E(\psi(p^hy+\xi;\boldsymbol \alpha)),
\end{equation}
with the coefficient $\mathfrak b_y(\xi)$ defined in \eqref{4.3}, and
\[
\psi^*(y;\boldsymbol \beta)=\beta_1\Phi_1(y)+\ldots +\beta_k\Phi_k(y).
\]
It then follows via orthogonality that
\begin{equation}\label{4.12}
\oint_{p^B}|\mathfrak f_h(\boldsymbol \alpha;\xi)|^{2s}\d\boldsymbol \alpha =\oint_{p^B}\oint_{p^{B-kh}}
|\mathfrak g_h(\boldsymbol \alpha,\boldsymbol \beta ;\xi)|^{2s}\d\boldsymbol \beta \d\boldsymbol \alpha .
\end{equation}

\par The polynomial system $\boldsymbol \Phi$ defined in \eqref{4.8} is $p^{c+h}$-spaced. Moreover, we assume 
that $h\le (1-\delta )H$, and hence $B-kh\ge \delta B-k$. We may therefore assume that $B-kh$ is sufficiently large 
in terms of $s$, $k$ and $\varepsilon$. In view of \eqref{4.3}, \eqref{4.11} and Hypothesis \ref{hypothesis3.2}, 
moreover, we may assume that the set of weights $(\boldsymbol{\mathfrak c}_y(\xi;\boldsymbol \alpha))$ satisfies 
Hypothesis \ref{hypothesis3.2}. Thus, we deduce from \eqref{3.12} and \eqref{3.13} that
\[
U_{s,k}^{B-kh}(\boldsymbol{\mathfrak c})\ll \bigl( p^{2(H-h)}\bigr)^{\lambda (s,k)+\varepsilon}U_{s,k}^{B-kh,H-h}
(\boldsymbol{\mathfrak c}).
\]
We therefore deduce from \eqref{3.8} and \eqref{4.10} that one has the bound
\begin{equation}\label{4.13}
\oint_{p^{B-kh}}|\mathfrak g_h(\boldsymbol \alpha,\boldsymbol \beta;\xi)|^{2s}\d\boldsymbol \beta \ll 
\bigl( p^{2(H-h)}\bigr)^{\lambda (s,k)+\varepsilon}I_1,
\end{equation}
where
\[
I_1=\rho_h(\xi)^{-2}\sum_{\eta\mmod{(p^{H-h},t^2)}}\rho_{H-h}(p^h\eta+\xi)^2
\oint_{p^{B-kh}}|\mathfrak g^*_{H-h}(\boldsymbol \alpha,\boldsymbol \beta ;\xi,\eta)|^{2s}\d\boldsymbol \beta ,
\]
and
\[
\mathfrak g^*_{H-h}(\boldsymbol \alpha,\boldsymbol \beta ;\xi,\eta)=\rho_{H-h}(p^h\eta+\xi)^{-1}
\sum_{y\equiv \eta\mmod{(p^{H-h},t^2)}}\mathfrak c_y(\xi;\boldsymbol \alpha)E(\psi^*(y;\boldsymbol \beta)).
\] 
The direct correspondence between residues $\zeta$ modulo $(p^H,t^2)$ and pairs $(\xi,\eta)$ of residues 
modulo $(p^h,t^2)$, and modulo $(p^{H-h},t^2)$, respectively, via the relation $\zeta=p^h\eta+\xi$, reveals that
\[
\mathfrak g^*_{H-h}(\boldsymbol \alpha,\boldsymbol \beta ;\xi,\eta)=\rho_{H}(\zeta)^{-1}
\sum_{\substack{y\in \mathfrak Z(t^2)\\ p^hy+\xi\equiv \zeta \mmod{(p^H,t^2)}}}\mathfrak a_{p^hy+\xi}
E\bigl( \psi(p^hy+\xi;\boldsymbol \alpha)+\psi^*(y;\boldsymbol \beta)\bigr) .
\]
In particular, we deduce from \eqref{4.12} and \eqref{4.13} that
\begin{align}
\sum_{\xi\mmod{(p^h,t^2)}}&\rho_h(\xi)^2\oint_{p^B}|\mathfrak f_h(\boldsymbol \alpha;\xi)|^{2s}
\d\boldsymbol \alpha \notag \\
&\ll \bigl( p^{2(H-h)}\bigr)^{\lambda(s,k)+\varepsilon}\sum_{\substack{\zeta\mmod{(p^H,t^2)}\\ 
\zeta=p^h\eta+\xi}}\rho_H(\zeta)^2\oint_{p^B}\oint_{p^{B-kh}}
|\mathfrak g^*_{H-h}(\boldsymbol \alpha,\boldsymbol \beta;\xi,\eta)|^{2s}\d\boldsymbol \beta \d\boldsymbol 
\alpha .\label{4.14}
\end{align}

\par By orthogonality, the mean value on the right hand side of \eqref{4.14} counts solutions of the system of 
congruences \eqref{4.4}, with each solution being counted with weight \eqref{4.5}, but subject to the additional 
condition that
\[
p^h\mathbf y+\xi\equiv p^h\mathbf z +\xi\equiv \zeta \mmod{(p^H,t^2)},
\]
and further, subject to the additional congruence conditions \eqref{4.9}. Note that the latter congruence conditions 
are generated by the integral over $\boldsymbol \beta$ in the mean value on the right hand side of \eqref{4.14}. The 
conditions \eqref{4.9} are implied by \eqref{4.4}, as we have shown above. Note also that when 
$\zeta=p^h\eta+\zeta$, one has
\[
\mathfrak g^*_{H-h}(\boldsymbol \alpha,\boldsymbol 0;\xi,\eta)=\rho_H(\zeta)^{-1}
\sum_{n\equiv \zeta \mmod{(p^H,t^2)}}\mathfrak a_ne(\psi(n;\boldsymbol \alpha))=
\mathfrak f_H(\boldsymbol \alpha;\zeta).
\]
Then, since the conditions \eqref{4.9} may be omitted, the definition \eqref{3.8} leads us to the conclusion that
\begin{align*}
U_{s,k}^{B,h}(\boldsymbol{\mathfrak a})&\ll \bigl( p^{2(H-h)}\bigr)^{\lambda(s,k)+\varepsilon}\rho_0(1)^{-2}
\sum_{\zeta\mmod{(p^H,t^2)}}\rho_H(\zeta)^2\oint_{p^B}|\mathfrak f_H(\boldsymbol \alpha;\zeta)|^{2s}
\d\boldsymbol \alpha \\
&=\bigl( p^{2(H-h)}\bigr)^{\lambda(s,k)+\varepsilon}U_{s,k}^{B,H}(\boldsymbol{\mathfrak a}).
\end{align*}
This completes the proof of the lemma.
\end{proof}

We require a crude estimate for $K_{a,b,c}^{r,\boldsymbol \varphi,\nu}(\boldsymbol{\mathfrak a})$. This we obtain 
by applying H\"older's inequality in combination with Lemma \ref{lemma4.4}.

\begin{lemma}\label{lemma4.5}
Suppose that $k\in \mathbb N$ and that $p$ is a prime number with $p>k$. Let $\tau$, $\varepsilon$ and $\delta$ 
be positive numbers with $\varepsilon<\tau<\delta<1$, and let $B$ be sufficiently large in terms of $s$, $k$ and 
$\varepsilon$. Write $H=\lceil B/k\rceil$ and suppose that $r$ and $\nu$ are non-negative integers with 
$1\le r\le k-1$. In addition, suppose that $0<\Lambda \le \lambda (s,k)$. Then for every 
$\boldsymbol \varphi\in \Phi_\tau(B)$, every sequence $(\boldsymbol{\mathfrak a})\in \mathbb D_0(B)$, and all 
non-negative integers $a$ and $b$ satisfying $\max \{a,b\}\le (1-\delta )H$, one has
\[
[[ K_{a,b,c}^{r,\boldsymbol \varphi,\nu}(\boldsymbol{\mathfrak a})]]_\Lambda 
\ll (p^{2H})^{\lambda(s,k)-\Lambda+\varepsilon}.
\]
\end{lemma}

\begin{proof}
A comparison with the statement and proof of \cite[Lemma 4.2]{Woo2019} reveals that the same argument as was 
applied to prove the latter remains valid, {\it mutatis mutandis}, in the present situation. It is important to note from 
the preamble to \eqref{3.9} that one has 
\[
\sum_{\zeta\mmod{(p^h,t^2)}}\rho_h(\zeta)^2\le \rho_0(1)^2,
\]
in alignment with the corresponding relation applied in the proof of \cite[Lemma 4.2]{Woo2019}. Moreover, in the 
conclusion of Lemma \ref{lemma4.4} above, we have a power of $p^{2(H-h)}$ in place of the corresponding power 
of $p^{H-h}$ in the conclusion of \cite[Lemma 4.1]{Woo2019}. Equipped with these observations, the reader will 
experience no difficulties in tracking through the argument of the proof of \cite[Lemma 4.2]{Woo2019} to obtain 
the conclusion of the present lemma.
\end{proof}

\section{The base case $k=1$: a variant of Hensel's lemma}
The base case for our nested induction is that in which $k=1$. The proof of Theorem \ref{theorem3.2} in this case is 
tantamount to a special case of Hensel's lemma in the ring $\mathbb Z[t]/(p^B,t^2)$, as we now explain.

\begin{lemma}\label{lemma5.1} For any prime number $p$, one has $\lambda(1,1)=0$.
\end{lemma}

\begin{proof} Let $\tau$ be a small positive number, and let $B$ be a positive integer sufficiently large in terms of 
$\tau$. Consider any sequence $(\boldsymbol{\mathfrak a})\in \mathbb D(B)$ and polynomial 
$\varphi\in \Phi_\tau(B)$. We suppose that $\varphi(u)$ is a $p^c$-spaced polynomial for some $c\ge \tau B$, so 
that 
\[
\varphi(u)=u+p^c\psi(u)+t\chi(u),
\]
for suitable polynomials $\psi(u),\chi(u)\in \mathbb Z[t,u]$. We observe that the mean value 
$U_{1,1}^B(\boldsymbol{\mathfrak a})$ counts the solutions of the congruence
\begin{equation}\label{5.1}
\varphi(x)\equiv \varphi(y)\mmod{(p^B,t^2)},
\end{equation}
with each solution $x,y\in \mathfrak Z(t^2)$ being counted with weight 
$\rho_0(1)^{-2}\mathfrak a_x{\overline{\mathfrak a}}_y$. By writing the congruence \eqref{5.1} in the form
\begin{equation}\label{5.2}
x-y+p^c(\psi(x)-\psi(y))+t(\chi(x)-\chi(y))\equiv 0\mmod{(p^B,t^2)},
\end{equation}
we see that
\begin{equation}\label{5.3}
x-y\equiv 0\mmod{(p^B,p^c,t)}.
\end{equation}
Note here that the ideal $(p^B,p^c,t)$ is, of course, equal to $(p^{\min\{B,c\}},t)$, and we have chosen the 
former as a less awkward avatar of the latter.\par

Since $x-y$ divides both $\psi(x)-\psi(y)$ and $\chi(x)-\chi(y)$, it follows by substituting \eqref{5.3} into the right 
hand side of \eqref{5.2} that
\begin{equation}\label{5.4}
x-y\equiv 0\mmod{(p^B,p^{2c},p^ct,t^2)}.
\end{equation}
By repeating this argument with \eqref{5.4} in place of \eqref{5.3}, and embarking on the inductive argument to 
hand, we find that for each non-negative integer $m$, one has
\[
x-y\equiv 0\mmod{(p^B,p^{(m+1)c},p^{mc}t,t^2)}.
\]
We take $m=\lceil 1/\tau\rceil$, and conclude that in fact $x-y\equiv 0\mmod{(p^B,t^2)}$. Thus, the solutions of 
\eqref{5.2} may be classified according to the common congruence class $\xi$ modulo $(p^B,t^2)$ of $x$ and $y$. 
Making use of the definitions \eqref{3.4}, \eqref{3.6} and \eqref{3.8}, we deduce via orthogonality that
\[
U_{1,1}^B(\boldsymbol{\mathfrak a})=\rho_0(1)^{-2}\sum_{\xi\mmod{(p^B,t^2)}}\rho_B(\xi)^2\oint_{p^B}
|\mathfrak f_B(\alpha;\xi)|^2\d\alpha =U_{1,1}^{B,B}(\boldsymbol {\mathfrak a}).
\]
Thus, we have
\[
\log \bigl( U_{1,1}^B(\boldsymbol{\mathfrak a})/U_{1,1}^{B,B}(\boldsymbol{\mathfrak a})\bigr) =0,
\]
and so it follows from the definitions \eqref{3.12} and \eqref{3.13} that one has $\lambda^*(1,1;\tau)=0$, and hence 
$\lambda(1,1)=0$. This completes the proof of the lemma.
\end{proof}

\section{Initial conditioning}
We continue our campaign on the proof of Theorem \ref{theorem3.2} in the case $k\ge 2$, assuming its validity for 
exponents smaller than $k$. Our initial step is to bound $U_{s,k}^B(\boldsymbol{\mathfrak a})$ in terms of a mean 
value of the shape $K_{a,b,c}^{1,\boldsymbol \varphi,\nu}(\boldsymbol{\mathfrak a})$. This section is devoted to the 
preliminary manoeuvres required for establishing such a bound. In this step, we encounter significant complications 
relative to the corresponding situation encountered in \cite{Woo2019} when handling efficient congruencing over 
the simpler ring $\mathbb Z$. We note, in particular, that it is in this initial conditioning process that we are forced 
to make use of the Hypothesis \ref{hypothesis3.2}.\par

Recall the definitions \eqref{3.12} and \eqref{3.13}. We fix a prime $p$ with $p>k$, also fix the parameter $s=k(k+1)/2$, 
and seek to establish that $\lambda(s,k)=0$. In order to confirm the latter, we proceed by deriving a contradiction to 
the assertion that $\lambda(s,k)=\Lambda>0$. In view of the upper bound \eqref{3.14}, we may suppose throughout 
that
\begin{equation}\label{6.1}
0<\Lambda\le s.
\end{equation}
Recall that we are implicitly assuming Hypothesis \ref{hypothesis3.2}, which asserts, via the implicit reference to 
Hypothesis \ref{hypothesis3.1}, that there exists a small positive number $\omega$ for which the conclusion of the 
hypothesis holds. The inductive argument calls for the assumption of this hypothesis in a nested form with a large, but 
finite, number of levels. We fix an admissible value for the smallest of the values $\omega$ that occur in this collection 
of levels, and introduce a heierarchy of sufficiently small positive numbers $\varepsilon$, $\tau$, $\delta$ and $\mu$ 
with
\begin{equation}\label{6.2}
\varepsilon <\tau<\delta<\mu<\omega <1.
\end{equation}
We suppose that each element in the heierarchy is sufficiently small in terms of $k$, $\Lambda$ and the 
larger elements of the heierarchy. Thus, for example, we are at liberty to assume that
\begin{equation}\label{6.3}
\varepsilon<k\Lambda \omega/8.
\end{equation}
By reference to the definition \eqref{3.13}, there is no loss of generality in supposing that 
$\lambda^*(s,k;\tau)\ge \Lambda-\varepsilon/2$. By virtue of \eqref{3.12}, there consequently exists a sequence 
$(B_m)_{m=1}^\infty$, with $B_m\rightarrow \infty$ as $m\rightarrow \infty$, having the property that whenever 
$m$ is sufficiently large, there exists $\boldsymbol \varphi_m\in \Phi_\tau(B_m)$ and 
$(\boldsymbol{\mathfrak a}^{(m)})\in \mathbb D(B_m)$ for which
\begin{equation}\label{6.4}
U_{s,k}^{B_m}(\boldsymbol{\mathfrak a}^{(m)})\ge (p^{2H_m})^{\Lambda-\varepsilon}
U_{s,k}^{B_m,H_m}(\boldsymbol{\mathfrak a}^{(m)}),
\end{equation}
in which $H_m=\lceil B_m/k\rceil$. We now fix any such integer $m$ having the property that $B=B_m$ is 
sufficiently large in terms of $k$, $\Lambda$, $\mu$, $\delta$, $\tau$, $\varepsilon$ and $\omega$. We suppress 
the index $m$ from this point on, so as to streamline our notation. We remark that since $B$ is supposed to be 
sufficiently large in terms of the aforementioned collection of parameters, then we may suppose that $\varepsilon B$ 
is also sufficiently large with respect to the same set of parameters.\par

On recalling Lemma \ref{lemma4.4}, we perceive that there is no loss of generality in supposing that, for every 
non-negative integer $h$ with $h\le (1-\delta)H$, uniformly in $\boldsymbol \varphi'\in \Phi_\tau(B)$ and 
$(\boldsymbol{\mathfrak a}')\in \mathbb D_0(B)$, one has
\begin{equation}\label{6.5}
U_{s,k}^{B,h}(\boldsymbol{\mathfrak a}')\le \bigl( p^{2(H-h)}\bigr)^{\lambda+\varepsilon}U_{s,k}^{B,H}(
\boldsymbol{\mathfrak a}').
\end{equation}

\par We are almost ready to announce our preliminary bound on $U_{s,k}^B(\boldsymbol{\mathfrak a})$ in terms of 
mean values of the shape $K_{a,b,c}^{1,\boldsymbol\varphi,\nu}(\boldsymbol{\mathfrak a})$. Here, we fix a choice 
for the parameter $\nu$ for future use by setting
\begin{equation}\label{6.6}
\nu = \lceil 8\varepsilon H/\Lambda \rceil .
\end{equation}
Notice that the assumption \eqref{6.3} ensures that one then has $\nu<\lfloor \omega B\rfloor$. Since we may 
suppose that $\boldsymbol \varphi\in \Phi_\tau(B)$, it follows that there exists an integer $c$ with $c\ge \tau B$ 
having the property that $\boldsymbol \varphi$ is $p^c$-spaced.\par

We first prepare an auxiliary estimate for a mean value associated with the exponential sum
\[
T_1(\boldsymbol \alpha)=\rho_0(1)^{-2s}\sum_{\zeta=1}^{p^\nu}\biggl| \sum_{\substack{\xi\mmod{(p^\nu,t^2)}\\ 
\xi\equiv \zeta\mmod{(p^\nu,t)}}}\rho_\nu(\xi)\mathfrak f_\nu(\boldsymbol \alpha;\xi)\biggr|^{2s}.
\]

\begin{lemma}\label{lemma6.1}
One has
\[
\oint_{p^B}T_1(\boldsymbol \alpha)\d\boldsymbol \alpha \ll p^{s\nu \varepsilon}
U_{s,k}^{B,\nu}(\boldsymbol{\mathfrak a}).
\]
\end{lemma}

\begin{proof} We begin by examining the mean value
\[
\oint_{p^B} \biggl| \sum_{\substack{\xi\mmod{(p^\nu,t^2)}\\ \xi\equiv \zeta \mmod{(p^\nu,t)}}}\rho_\nu (\xi)
\mathfrak f_\nu(\boldsymbol \alpha ;\xi)\biggr|^{2s}\d\boldsymbol \alpha =\sum_{\boldsymbol \xi}
I(\boldsymbol \xi)\prod_{i=1}^{2s}\rho_\nu (\xi_i),
\]
where
\[
I(\boldsymbol \xi)=\oint_{p^B}\prod_{i=1}^s\mathfrak f_\nu(\boldsymbol \alpha;\xi_i){\overline{\mathfrak f_\nu(
\boldsymbol \alpha ;\xi_{s+i})}}\d\boldsymbol \alpha .
\]
Here, the summation over $\boldsymbol \xi$ is taken over $\xi_i$ modulo $(p^\nu,t^2)$ with 
$\xi_i\equiv \zeta \mmod{(p^\nu,t)}$ for $1\le i\le 2s$. The mean value $I(\boldsymbol \xi)$ counts the solutions of 
the simultaneous congruences \eqref{3.7} with $\mathbf x,\mathbf y\in \mathfrak Z(t^2)$ satisfying the congruence 
conditions
\[
x_i\equiv \xi_i\mmod{(p^\nu,t^2)}\quad \text{and}\quad y_i\equiv \xi_{s+i}\mmod{(p^\nu,t^2)}\quad 
(1\le i\le s),
\]
and with each solution counted with weight
\[
\biggl( \prod_{i=1}^s\rho_\nu(\xi_i)^{-1}\mathfrak a_{x_i}\biggr) \biggl( \prod_{i=1}^s \rho_\nu(\xi_{s+i})^{-1}
{\overline{\mathfrak a}}_{y_i}\biggr) .
\]

\par Observe that $\boldsymbol \varphi$ is a $p^c$-spaced system, so it follows as in the argument of the proof of 
Lemma \ref{lemma4.4} leading to \eqref{4.6} that, in the first instance, there is no loss of generality in supposing that
\[
\varphi_j(u)=u^j+u^{k+1}\bigl( p^c\psi_j(u)+t\chi_j(u)\bigr)\quad (1\le j\le k),
\]
for appropriate polynomials $\psi_j(u),\chi_j(u)\in \mathbb Z[t,u]$. Our assumptions on $\nu$ ensure that $\nu$ is 
small compared to $c$. Thus, it follows from the congruence in \eqref{3.7} corresponding to $j=1$ that
\[
\sum_{i=1}^s \varphi_1^*(x_i)\equiv \sum_{i=1}^s \varphi_1^*(x_{s+i})\mmod{(p^\nu,t^2)},
\]
in which we have written
\[
\varphi_1^*(u)=u+tu^{k+1}\chi_1(u).
\]
Consequently, one has a solution of \eqref{3.7} to be counted by the mean value $I(\boldsymbol \xi)$ only when
\[
\sum_{i=1}^s\bigl( \varphi_1^*(\xi_i)-\varphi_1^*(\xi_{s+i})\bigr)\equiv 0\mmod{(p^\nu,t^2)}.
\]
Note here that for $1\le i\le 2s$, one has $\xi_i\equiv \zeta\mmod{(p^\nu,t)}$, and thus
\[
\varphi_1^*(\xi_i)\equiv \xi_i+t\zeta^{k+1}\chi_1(\zeta)\mmod{(p^\nu,t^2)}.
\]
Then we may conclude that one has a solution of \eqref{3.7} to be counted by $I(\boldsymbol \xi)$ only when
\begin{equation}\label{6.7}
\sum_{i=1}^s\bigl( \xi_i-\xi_{s+i}\bigr)\equiv 0\mmod{(p^\nu,t^2)}.
\end{equation}
Recalling again that for $1\le i\le 2s$, we have $\xi_i\equiv \zeta\mmod{(p^\nu,t)}$, we note that the congruence 
\eqref{6.7} serves only to constrain the congruence classes of the coefficients of $t$ in the variables $\xi_i$.\par

Working backwards from orthogonality, we may incorporate the condition \eqref{6.7} into the mean value to be 
estimated. Here, it is useful to observe that when $I(\boldsymbol \xi)$ is non-zero, it follows from the elementary 
inequality $|z_1\cdots z_n|\le |z_1|^n+\ldots +|z_n|^n$ that
\[
I(\boldsymbol \xi)\prod_{i=1}^{2s}\rho_\nu(\xi_i)\le \sum_{i=1}^{2s}\oint_{p^B}
|{\widetilde \rho}_{\nu,i}(\boldsymbol \xi)
{\widetilde{\mathfrak f}}_{\nu,i}(\boldsymbol \alpha;\boldsymbol \xi)|^{2s/(2s-1)},
\] 
where
\[
{\widetilde \rho}_{\nu,i}(\boldsymbol \xi)=\prod_{\substack{1\le j\le 2s\\ j\ne i}}\rho_\nu(\xi_j)\quad 
\text{and} \quad {\widetilde{\mathfrak f}}_{\nu,i}(\boldsymbol \alpha;\boldsymbol \xi)=
\prod_{\substack{1\le j\le 2s\\ j\ne i}}\mathfrak f_\nu (\boldsymbol \alpha ;\xi_j).
\]
Thus, we infer that
\begin{align*}
\oint_{p^B} &\biggl| \sum_{\substack{\xi\mmod{(p^\nu,t^2)}\\ \xi\equiv \zeta\mmod{(p^\nu,t)}}}\rho_\nu(\xi) 
\mathfrak f_\nu(\boldsymbol \alpha ;\xi)\biggr|^{2s}\d\boldsymbol \alpha \\
&\le \sum_{i=1}^{2s}
\sum_{\substack{\boldsymbol \xi\mmod{(p^\nu,t^2)}\\
\boldsymbol \xi\equiv \zeta\mmod{(p^\nu,t)}\\
\xi_1+\ldots +\xi_s\equiv \xi_{s+1}+\ldots +\xi_{2s}\mmod{(p^\nu,t^2)}}} \oint_{p^B}
|{\widetilde{\rho}}_{\nu,i}(\boldsymbol \xi)
{\widetilde{\mathfrak f}}_{\nu,i}(\boldsymbol \alpha;\boldsymbol \xi)|^{2s/(2s-1)}
\d\boldsymbol \alpha .
\end{align*}
By fixing the residue class of any $\xi_j$ modulo $(p^\nu,t)$ for $2s-1$ distinct choices of $j$ with $1\le j\le 2s$, the 
final class is fixed within the summation over $\boldsymbol \xi$. Thus, on employing symmetry amongst the indices, we 
see that
\begin{align}
\oint_{p^B} &\biggl| \sum_{\substack{\xi\mmod{(p^\nu,t^2)}\\ \xi\equiv \zeta\mmod{(p^\nu,t)}}}\rho_\nu(\xi) 
\mathfrak f_\nu(\boldsymbol \alpha ;\xi)\biggr|^{2s}\d\boldsymbol \alpha \notag \\
&\ll \oint_{p^B}\biggl( \sum_{\substack{\xi\mmod{(p^\nu,t^2)}\\ \xi\equiv \zeta\mmod{(p^\nu,t)}}}|\rho_\nu(\xi) 
\mathfrak f_\nu(\boldsymbol \alpha ;\xi)|^{2s/(2s-1)}\biggr)^{2s-1}\d\boldsymbol \alpha .\label{6.8}
\end{align}

\par By applying H\"older's inequality, we find that
\begin{align*}
\biggl( &\sum_{\substack{\xi\mmod{(p^\nu,t^2)}\\ \xi\equiv \zeta\mmod{(p^\nu,t)}}}|\rho_\nu(\xi) 
\mathfrak f_\nu(\boldsymbol \alpha ;\xi)|^{2s/(2s-1)}\biggr)^{2s-1}\\
&\le \biggl( \sum_{\substack{\xi\mmod{(p^\nu,t^2)}\\ \xi\equiv \zeta\mmod{(p^\nu,t)}}}\rho_\nu(\xi) \biggr)^{2s-2}
\biggl( \sum_{\substack{\xi\mmod{(p^\nu,t^2)}\\ \xi\equiv \zeta\mmod{(p^\nu,t)}}}\rho_\nu(\xi)^2
|\mathfrak f_\nu(\boldsymbol \alpha ;\xi)|^{2s}\biggr).
\end{align*}
It is at this point that we make use of Hypothesis \ref{hypothesis3.1}. By making use of the latter in combination 
with \eqref{3.2} and Cauchy's inequality, we see that
\begin{align*}
\biggl( \sum_{\substack{\xi\mmod{(p^\nu,t^2)}\\ \xi\equiv \zeta\mmod{(p^\nu,t)}}}\rho_\nu(\xi) \biggr)^2&
\le p^\nu \sum_{\substack{\xi\mmod{(p^\nu,t^2)}\\ \xi\equiv \zeta\mmod{(p^\nu,t)}}}\rho_\nu(\xi)^2\\
&=p^\nu \sum_{\eta_1=1}^{p^\nu}\| {\boldsymbol{\mathfrak a}}_\nu(\zeta+\eta_1t)\|_0\cdot \| 
{\boldsymbol{\mathfrak a}}_\nu(\zeta+\eta_1t)\|_\infty^2\\
&\ll p^{\nu\varepsilon}\|{\boldsymbol{\mathfrak a}}_0(1)\|_\infty^2\sum_{\eta_0=1}^{p^\nu}\sum_{\eta_1=1}^{p^\nu}
\| {\boldsymbol{\mathfrak a}}_\nu(\eta_0+\eta_1t)\|_0 .
\end{align*}
Thus, we infer that
\[
\biggl( \sum_{\substack{\xi\mmod{(p^\nu,t^2)}\\ \xi\equiv \zeta\mmod{(p^\nu,t)}}}\rho_\nu(\xi) \biggr)^2\ll 
 p^{\nu \varepsilon}\|{\boldsymbol{\mathfrak a}}_0(1)\|_\infty^2 \cdot \|{\boldsymbol{\mathfrak a}}_0(1)\|_0
 =p^{\nu \varepsilon}\rho_0(1)^2,
\]
whence
\begin{align*}
\biggl( \sum_{\substack{\xi\mmod{(p^\nu,t^2)}\\ \xi\equiv \zeta\mmod{(p^\nu,t)}}}&|\rho_\nu(\xi) 
\mathfrak f_\nu(\boldsymbol \alpha ;\xi)|^{2s/(2s-1)}\biggr)^{2s-1}\\
&\ll \bigl( p^{\nu\varepsilon}\rho_0(1)^2\bigr)^{s-1}\sum_{\substack{\xi\mmod{(p^\nu,t^2)}\\ \xi\equiv \zeta\mmod{(p^\nu,t)}}}\rho_\nu(\xi)^2
|\mathfrak f_\nu(\boldsymbol \alpha ;\xi)|^{2s}.
\end{align*}

\par By substituting this last upper bound into \eqref{6.8}, we arrive at the estimate
\begin{align*}
\sum_{\zeta=1}^{p^\nu}\oint_{p^B} &\biggl| \sum_{\substack{\xi\mmod{(p^\nu,t^2)}\\ 
\xi\equiv \zeta\mmod{(p^\nu,t)}}}\rho_\nu(\xi) \mathfrak f_\nu(\boldsymbol \alpha ;\xi)\biggr|^{2s}
\d\boldsymbol \alpha \\
&\ll p^{s\nu \varepsilon}\rho_0(1)^{2s-2}\sum_{\xi\mmod{(p^\nu,t^2)}}\rho_\nu(\xi)^2
\oint_{p^B}|\mathfrak f_\nu(\boldsymbol \alpha;\xi)|^{2s}\d\boldsymbol \alpha .
\end{align*}
On recalling the definition of $T_1(\boldsymbol \alpha)$ together with \eqref{3.8}, we therefore conclude that
\[
\oint_{p^B}T_1(\boldsymbol \alpha)\d\boldsymbol \alpha \ll 
p^{s\nu \varepsilon}\rho_0(1)^{-2}\sum_{\xi\mmod{(p^\nu,t^2)}}
\rho_\nu(\xi)^2\oint_{p^B}|\mathfrak f_\nu(\boldsymbol \alpha;\xi)|^{2s}\d\boldsymbol \alpha 
=p^{s\nu \varepsilon}U_{s,k}^{B,\nu}(\boldsymbol{\mathfrak a}).
\]
This completes the proof of the lemma.
\end{proof}

We are now equipped to establish the first of our conditioning lemmata.

\begin{lemma}\label{lemma6.2}
One has $U_{s,k}^B(\boldsymbol{\mathfrak a})\ll (p^{2\nu})^s 
K_{\nu,\nu,c}^{1,\boldsymbol \varphi,\nu}(\boldsymbol{\mathfrak a})$.
\end{lemma}

\begin{proof} Using the definition of $f_{\boldsymbol{\mathfrak a}}(\boldsymbol \alpha)$ and the definition 
\eqref{3.4}, one sees that
\begin{equation}\label{6.9}
\rho_0(1)^{2s}f_{\boldsymbol{\mathfrak a}}(\boldsymbol \alpha)^{2s}=\sum_{\xi\mmod{(p^\nu,t^2)}}\rho_\nu(\xi)
\mathfrak f_\nu(\boldsymbol \alpha;\xi)\bigl( \rho_0(1)\mathfrak f_{\boldsymbol{\mathfrak a}}(\boldsymbol \alpha)
\bigr)^{2s-1}.
\end{equation}
Moreover, for each residue $\xi$ modulo $(p^\nu,t^2)$, one has
\begin{equation}\label{6.10}
\rho_0(1)f_{\boldsymbol{\mathfrak a}}(\boldsymbol \alpha)=\sum_{\substack{\eta\mmod{(p^\nu,t^2)}\\ 
\eta\equiv \xi\mmod{(p^\nu,t)}}}\rho_\nu(\eta)\mathfrak f_\nu(\boldsymbol \alpha ;\eta)+
\sum_{\substack{\eta\mmod{(p^\nu,t^2)}\\ \eta\not\equiv \xi \mmod{(p^\nu,t)}}}\rho_\nu(\eta)\mathfrak f_\nu(
\boldsymbol \alpha ;\eta),
\end{equation}
whence
\begin{equation}\label{6.11}
|\rho_0(1)f_{\boldsymbol{\mathfrak a}}(\boldsymbol \alpha)|^{2s-1}\ll 
\biggl|\sum_{\substack{\eta\mmod{(p^\nu,t^2)}\\ 
\eta\equiv \xi\mmod{(p^\nu,t)}}}\rho_\nu(\eta)\mathfrak f_\nu(\boldsymbol \alpha ;\eta)\biggr|^{2s-1}+
\biggl| \sum_{\substack{\eta\mmod{(p^\nu,t^2)}\\ \eta\not\equiv \xi \mmod{(p^\nu,t)}}}
\rho_\nu(\eta)\mathfrak f_\nu(\boldsymbol \alpha ;\eta)\biggr|^{2s-1}.
\end{equation}
By substituting \eqref{6.11} into \eqref{6.9} and recalling the definition of $T_1(\boldsymbol \alpha)$ from the 
preamble to Lemma \ref{lemma6.1}, we arrive at the relation
\begin{equation}\label{6.12}
|f_{\boldsymbol{\mathfrak a}}(\boldsymbol \alpha)|^{2s}\ll T_1(\boldsymbol \alpha)+T_2(\boldsymbol \alpha),
\end{equation}
where
\begin{equation}\label{6.13}
T_2(\boldsymbol \alpha)=\rho_0(1)^{-2s}\sum_{\xi\mmod{(p^\nu,t^2)}}\rho_\nu(\xi)|\mathfrak f_\nu(\boldsymbol 
\alpha ;\xi)|\,\biggl| \sum_{\substack{\eta\mmod{(p^\nu,t^2)}\\ \eta\not\equiv \xi \mmod{(p^\nu,t)}}}\rho_\nu(\eta)
\mathfrak f_\nu(\boldsymbol \alpha ;\eta)\biggr|^{2s-1}.
\end{equation}

\par We must double up the number of generating functions $\mathfrak f_\nu(\boldsymbol \alpha;\xi)$ in the mean 
value augured by \eqref{6.13}. This we achieve with a swift application of Cauchy's inequality, obtaining the bound
\begin{equation}\label{6.14}
T_2(\boldsymbol \alpha)\le T_3(\boldsymbol \alpha)^{1/2}T_4(\boldsymbol \alpha)^{1/2},
\end{equation}
where
\begin{equation}\label{6.15}
T_3(\boldsymbol \alpha)=\rho_0(1)^{-2s}\sum_{\xi\mmod{(p^\nu,t^2)}}\rho_\nu(\xi)^2|\mathfrak f_\nu(\boldsymbol 
\alpha ;\xi)|^2\biggl| \sum_{\substack{\eta\mmod{(p^\nu,t^2)}\\ \eta\not\equiv \xi \mmod{(p^\nu,t)}}}\rho_\nu(\eta)
\mathfrak f_\nu(\boldsymbol \alpha ;\eta)\biggr|^{2s-2}
\end{equation}
and
\begin{equation}\label{6.16}
T_4(\boldsymbol \alpha)=\rho_0(1)^{-2s}\sum_{\xi\mmod{(p^\nu,t^2)}}\biggl| 
\sum_{\substack{\eta\mmod{(p^\nu,t^2)}\\ 
\eta\not\equiv \xi \mmod{(p^\nu,t)}}}\rho_\nu(\eta)\mathfrak f_\nu(\boldsymbol \alpha ;\eta)\biggr|^{2s}.
\end{equation}

\par We next analyse the expression $T_4(\boldsymbol \alpha)$, first removing the awkward summation restriction 
$\eta \not\equiv \xi \mmod{(p^\nu,t)}$. Thus, reversing the relation \eqref{6.10}, we see that
\[
\biggl| \sum_{\substack{\eta\mmod{(p^\nu,t^2)}\\ 
\eta\not\equiv \xi \mmod{(p^\nu,t)}}}\rho_\nu(\eta)\mathfrak f_\nu(\boldsymbol \alpha ;\eta)\biggr|^{2s}
\ll |\rho_0(1)f_{\boldsymbol{\mathfrak a}}(\boldsymbol \alpha)|^{2s}+
\biggl| \sum_{\substack{\eta\mmod{(p^\nu,t^2)}\\ 
\eta\equiv \xi \mmod{(p^\nu,t)}}}\rho_\nu(\eta)\mathfrak f_\nu(\boldsymbol \alpha ;\eta)\biggr|^{2s}.
\]
Hence, we see from \eqref{6.16} and the preamble to Lemma \ref{lemma6.1} that
\[
T_4(\boldsymbol \alpha)\ll \sum_{\xi\mmod{(p^\nu,t^2)}}|f_{\boldsymbol{\mathfrak a}}(\boldsymbol \alpha)|^{2s}
+p^\nu T_1(\boldsymbol \alpha),
\]
whence, as a consequence of \eqref{3.6} and Lemma \ref{lemma6.1},
\begin{equation}\label{6.17}
\oint_{p^B}T_4(\boldsymbol \alpha)\d\boldsymbol \alpha \ll p^{2\nu}U_{s,k}^B(\boldsymbol{\mathfrak a})+
p^{\nu+s\nu \varepsilon}U_{s,k}^{B,\nu}(\boldsymbol{\mathfrak a}).
\end{equation}

\par Next, we find through an application of H\"older's inequality that
\begin{align*}
\Biggl| &\sum_{\substack{\eta\mmod{(p^\nu.t^2)}\\ \eta\not\equiv \xi \mmod{(p^\nu,t)}}}\rho_\nu(\eta)
\mathfrak f_\nu(\boldsymbol \alpha ;\eta)\Biggr|^{2s-2}\\
&\le \Biggl( \sum_{\eta\mmod{(p^\nu,t^2)}}\rho_\nu(\eta)^2\Biggr)^{s-2}
\Biggl( \sum_{\eta\mmod{(p^\nu,t^2)}}1\Biggr)^{s-1}
\Biggl( \sum_{\substack{\eta\mmod{(p^\nu,t^2)}\\ \eta\not\equiv \xi \mmod{(p^\nu,t)}}}\rho_\nu(\eta)^2
|\mathfrak f_\nu(\boldsymbol \alpha ;\eta)|^{2s-2}\Biggr)\\
&\le \rho_0(1)^{2s-4}(p^{2\nu})^{s-1}
\sum_{\substack{\eta\mmod{(p^\nu,t^2)}\\ \eta\not\equiv \xi \mmod{(p^\nu,t)}}}
\rho_\nu(\eta)^2|\mathfrak f_\nu(\boldsymbol \alpha ;\eta)|^{2s-2}.
\end{align*}
Then we deduce from \eqref{6.15} that
\begin{align*}
T_3(\boldsymbol \alpha)\le \rho_0(1)^{-4}(p^{2\nu})^{s-1}\sum_{\xi\mmod{(p^\nu,t^2)}}
\sum_{\substack{\eta\mmod{(p^\nu,t^2)}\\ \eta\not\equiv \xi \mmod{(p^\nu,t)}}}\rho_\nu(\xi)^2\rho_\nu(\eta)^{2}
|\mathfrak f_\nu(\boldsymbol \alpha ;\xi)|^2|\mathfrak f_\nu(\boldsymbol \alpha ;\eta)|^{2s-2}.
\end{align*}
We therefore conclude from \eqref{3.16} that
\begin{equation}\label{6.18}
\oint_{p^B}T_3(\boldsymbol \alpha)\,{\rm d}\boldsymbol \alpha\le (p^{2\nu})^{s-1}
K_{\nu,\nu,c}^{1,\boldsymbol \varphi,\nu}(\boldsymbol {\mathfrak a}).
\end{equation}

\par On substituting \eqref{6.17} and \eqref{6.18} into \eqref{6.14}, we find via Schwarz's inequality that
\begin{align*}
\oint_{p^B}T_2(\boldsymbol \alpha)\d\boldsymbol \alpha &\le 
\biggl( \oint_{p^B}T_3(\boldsymbol \alpha)\d\boldsymbol \alpha\biggr)^{1/2}\biggl( \oint_{p^B}T_4(\boldsymbol 
\alpha)\d\boldsymbol \alpha \biggr)^{1/2}\\
&\ll \bigl( (p^{2\nu})^{s-1}K_{\nu,\nu,c}^{1,\boldsymbol \varphi,\nu}(\boldsymbol {\mathfrak a})\bigr)^{1/2}
\bigl( p^{2\nu}U_{s,k}^B(\boldsymbol {\mathfrak a})+
p^{\nu+s\nu\varepsilon}U_{s,k}^{B,\nu}(\boldsymbol {\mathfrak a})\bigr)^{1/2}.
\end{align*}
Consequently, on recalling Lemma \ref{lemma6.1} and \eqref{6.12}, we obtain the upper bound
\begin{align*}
\oint_{p^B}|f_{\boldsymbol{\mathfrak a}}(\boldsymbol \alpha)|^{2s}\d\boldsymbol \alpha &\ll \oint_{p^B}
T_1(\boldsymbol \alpha)\d\boldsymbol \alpha+\oint_{p^B}T_2(\boldsymbol \alpha)\d\boldsymbol \alpha \\
&\ll p^{s\nu\varepsilon}U_{s,k}^{B,\nu}(\boldsymbol {\mathfrak a})+(p^\nu)^{s-1}
\bigl(K_{\nu,\nu,c}^{1,\boldsymbol \varphi,\nu}
(\boldsymbol {\mathfrak a})\bigr)^{1/2}\bigl( p^{2\nu}U_{s,k}^B(\boldsymbol {\mathfrak a})+
p^{\nu+s\nu\varepsilon}U_{s,k}^{B,\nu}(\boldsymbol {\mathfrak a})\bigr)^{1/2}.
\end{align*}
Thus, from \eqref{3.6}, we deduce that
\[
U_{s,k}^B(\boldsymbol {\mathfrak a})\ll p^{s\nu\varepsilon}U_{s,k}^{B,\nu}(\boldsymbol {\mathfrak a})+
(p^\nu)^{s-1}\bigl(K_{\nu,\nu,c}^{1,\boldsymbol \varphi,\nu}(\boldsymbol {\mathfrak a})\bigr)^{1/2}
\bigl( p^{2\nu}U_{s,k}^B(\boldsymbol {\mathfrak a})+p^{\nu+s\nu\varepsilon}
U_{s,k}^{B,\nu}(\boldsymbol {\mathfrak a})\bigr)^{1/2},
\]
whence
\begin{align}
U_{s,k}^B(\boldsymbol {\mathfrak a})&\ll p^{s\nu\varepsilon}U_{s,k}^{B,\nu}(\boldsymbol {\mathfrak a})
+(p^{2\nu})^sK_{\nu,\nu,c}^{1,\boldsymbol \varphi,\nu}(\boldsymbol {\mathfrak a})+(p^\nu)^{s-1}
\bigl( K_{\nu,\nu,c}^{1,\boldsymbol \varphi,\nu}(\boldsymbol {\mathfrak a})\bigr)^{1/2}
\bigl(p^{\nu+s\nu \varepsilon}U_{s,k}^{B,\nu}(\boldsymbol {\mathfrak a})\bigr)^{1/2}\notag \\
&\ll p^{s\nu\varepsilon}U_{s,k}^{B,\nu}(\boldsymbol {\mathfrak a})+
(p^{2\nu})^sK_{\nu,\nu,c}^{1,\boldsymbol \varphi,\nu}(\boldsymbol {\mathfrak a}).\label{6.19}
\end{align}

\par We now make use of the estimate \eqref{6.5} to obtain the bound
\[
U_{s,k}^{B,\nu}(\boldsymbol{\mathfrak a})\le \bigl( p^{2(H-\nu)}\bigr)^{\Lambda+\varepsilon} U_{s,k}^{B,H}
(\boldsymbol{\mathfrak a}).
\]
Since the definition \eqref{6.6} ensures that
\[
s\nu \varepsilon+2(\Lambda+\varepsilon)(H-\nu)-2(\Lambda-\varepsilon)H<4\varepsilon H-
(2\Lambda-s\varepsilon)\nu <-4\varepsilon H,
\]
we discern via \eqref{6.4} that
\[
p^{s\nu \varepsilon}U_{s,k}^{B,\nu}(\boldsymbol{\mathfrak a})\le p^{-4\varepsilon H}(p^{2H})^{\Lambda -\varepsilon}
U_{s,k}^{B,H}(\boldsymbol{\mathfrak a})\le p^{-4\varepsilon H}U_{s,k}^B(\boldsymbol{\mathfrak a}).
\]
Thus, we infer from \eqref{6.19} that
\[
U_{s,k}^B(\boldsymbol{\mathfrak a})\ll p^{-4\varepsilon H}U_{s,k}^B(\boldsymbol{\mathfrak a})+
(p^{2\nu})^sK_{\nu,\nu,c}^{1,\boldsymbol \varphi,\nu}(\boldsymbol{\mathfrak a}),
\]
whence
\[
U_{s,k}^B(\boldsymbol{\mathfrak a})\ll (p^{2\nu})^s
K_{\nu,\nu,c}^{1,\boldsymbol \varphi,\nu}(\boldsymbol{\mathfrak a}).
\]
This completes the proof of the lemma.
\end{proof}

We amplify the congruence condition implicit in the mean value 
$K_{\nu,\nu,c}^{1,\boldsymbol \varphi,\nu}(\boldsymbol{\mathfrak a})$ in order to gain sufficient traction to 
pursue our main iteration. At this point, we introduce the last of our small quantities, namely the parameter 
$\theta$, defined by putting
\begin{equation}\label{6.20}
\theta =\lceil \mu H\rceil .
\end{equation}
An essentially trivial bound presented in the next lemma eases our discussion.

\begin{lemma}\label{lemma6.3}
Suppose that $a$ and $b$ are non-negative integers with $a\le b$. Then for any positive number $w$ and any 
$\xi\in \mathbb Z[t]$, one has
\[
\rho_a(\xi)^2|\mathfrak f_a(\boldsymbol \alpha ;\xi)|^{2w}\le (p^{2b-2a})^w
\sum_{\substack{\zeta\mmod{(p^b,t^2)}\\ \zeta\equiv \xi\mmod{(p^a,t^2)}}}\rho_b(\zeta)^2
|\mathfrak f_b(\boldsymbol \alpha ;\zeta)|^{2w}.
\]
\end{lemma}

\begin{proof} By applying H\"older's inequality together with \eqref{3.4}, we see that
\begin{align}
\rho_a(\xi)^{2w}|\mathfrak f_a(\boldsymbol \alpha;\xi)|^{2w}&=\biggl| \sum_{\substack{\zeta \mmod{(p^b,t^2)}\\
\zeta\equiv \xi\mmod{(p^a,t^2)}}}\rho_b(\zeta)\mathfrak f_b(\boldsymbol \alpha;\zeta)\biggr|^{2w}\notag \\
&\le V_1^wV_2^{w-1} \sum_{\substack{\zeta \mmod{(p^b,t^2)}\\
\zeta\equiv \xi\mmod{(p^a,t^2)}}}\rho_b(\zeta)^2|\mathfrak f_b(\boldsymbol \alpha;\zeta)|^{2w},\label{6.21}
\end{align}
where
\[
V_1=\sum_{\substack{\zeta \mmod{(p^b,t^2)}\\ \zeta\equiv \xi\mmod{(p^a,t^2)}}}1\quad \text{and}
\quad V_2=\sum_{\substack{\zeta \mmod{(p^b,t^2)}\\ \zeta\equiv \xi\mmod{(p^a,t^2)}}}\rho_b(\zeta)^2.
\]
However, we have $V_1=p^{2b-2a}$, and in view of the definition \eqref{3.2},
\[
V_2\le \| \boldsymbol{\mathfrak a}_a(\xi)\|_\infty^2
\sum_{\substack{\zeta \mmod{(p^b,t^2)}\\ \zeta\equiv \xi\mmod{(p^a,t^2)}}}
\| \boldsymbol{\mathfrak a}_b(\zeta)\|_0\le \| \boldsymbol{\mathfrak a}_a(\xi)\|_\infty^2\cdot 
\| \boldsymbol{\mathfrak a}_a(\xi)\|_0=\rho_a(\xi)^2.
\]
The conclusion of the lemma follows on substituting these relations into \eqref{6.21}.
\end{proof}

We at last reach the primary objective of this section.

\begin{lemma}\label{lemma6.4}
One has $U_{s,k}^B(\boldsymbol{\mathfrak a})\ll p^{2s\theta}K_{\theta,\theta,c}^{1,\boldsymbol \varphi,\nu}
(\boldsymbol{\mathfrak a})$.
\end{lemma}

\begin{proof} By appealing to \eqref{3.17}, we find via two applications of Lemma \ref{lemma6.3} that
\begin{align*}
 \rho_\nu(\xi)^2&\rho_\nu(\eta)^2K_{\nu,\nu,c}^{1,\boldsymbol \varphi,\nu}(\boldsymbol{\mathfrak a};\xi,\eta)\\
&\le (p^{2\theta-2\nu})^s\sum_{\xi',\eta'\mmod{(p^\theta,t^2)}}\oint_{p^B}
\rho_\theta(\xi')^2\rho_\theta(\eta')^2|\mathfrak f_\theta(\boldsymbol \alpha;\xi')^2\mathfrak 
f_\theta(\boldsymbol \alpha;\eta')^{2s-2}|\d\boldsymbol \alpha ,
\end{align*}
where the summation over $\xi'$ and $\eta'$ is subject to the conditions
\[
\xi'\equiv \xi\mmod{p^\nu}\quad \text{and}\quad \eta'\equiv \eta\mmod{p^\nu}.
\]
We therefore deduce from \eqref{3.16} that
\begin{align*}
K_{\nu,\nu,c}^{1,\boldsymbol \varphi,\nu}(\boldsymbol{\mathfrak a})&\le \rho_0(1)^{-4}(p^{2\theta -2\nu})^s
\sum_{\xi'\mmod{(p^\theta,t^2)}}
\sum_{\substack{\eta'\mmod{(p^\theta,t^2)}\\ \eta'\not\equiv \xi'\mmod{(p^\nu,t)}}}\rho_\theta(\xi')^2
\rho_\theta(\eta')^2K_{\theta,\theta,c}^{1,\boldsymbol \varphi,\nu}(\boldsymbol{\mathfrak a};\xi',\eta')
\\
&=(p^{2\theta-2\nu})^sK_{\theta,\theta,c}^{1,\boldsymbol \varphi,\nu}(\boldsymbol{\mathfrak a}).
\end{align*}
Therefore, by substituting this bound into the estimate supplied by Lemma \ref{lemma6.2}, we obtain
\[
U_{s,k}^B(\boldsymbol{\mathfrak a})\ll (p^{2\theta-2\nu})^s(p^{2\nu})^s
K_{\theta,\theta,c}^{1,\boldsymbol \varphi,\nu}(\boldsymbol{\mathfrak a}).
\]
The conclusion of the lemma follows.
\end{proof}

\section{Congruencing operations}
Throughout, we suppose that $k\ge 2$, and we assume the validity of Theorem \ref{theorem3.2} for exponents 
smaller than $k$. We continue to make use of the conventions of the previous section, and in particular we write 
$R=r(r+1)/2$.

\begin{lemma}\label{lemma7.1}
Suppose that $a$, $b$ and $r$ are integers with
\begin{equation}\label{7.1}
1\le r\le k-1,\quad a\ge \delta \theta,\quad b\ge \delta \theta\quad \text{and}\quad ra\le (k-r+1)b\le B.
\end{equation}
Put
\begin{equation}\label{7.2}
b'=\lceil (k-r+1)b/r\rceil .
\end{equation}
Then one has
\[
K_{a,b,c}^{r,\boldsymbol \varphi,\nu}(\boldsymbol{\mathfrak a})\ll \bigl( p^{2\nu}\bigr)^{2k^2}
K_{b',b,c}^{r,\boldsymbol \varphi,\nu}(\boldsymbol{\mathfrak a}).
\]
\end{lemma}

\begin{proof} The essential strategy of the proof of this lemma follows the plan adopted in the corresponding 
proof of \cite[Lemma 7.1]{Woo2019}. However, the fact that we are now working in $\mathbb Z[t]/(p^B,t^2)$ in 
place of $\mathbb Z$ leads to significant complications that warrant a detailed account. We begin by examining the 
definition \eqref{3.16} of the mean value $K_{a,b,c}^{r,\boldsymbol \varphi,\nu}(\boldsymbol{\mathfrak a};\xi,\eta)$, 
in which we suppose that $\eta\not\equiv \xi\mmod{(p^\nu,t)}$. The latter condition offers us permission to assume 
that for some non-negative integer $\gamma$ with $\gamma<\nu$, one has
\[
\xi\equiv \eta\mmod{(p^\gamma,t)}\quad \text{and}\quad \xi\not\equiv \eta\mmod{(p^{\gamma+1},t)}.
\]
Write $\xi=\xi_0+\xi_1t$ and $\eta=\eta_0+\eta_1t$ with $\xi_i,\eta_i\in \mathbb Z$. We define the integer 
$\omega$ by putting $\omega=(\xi_0-\eta_0)p^{-\gamma}$, so that $(\omega,p)=1$. It is convenient for future 
use to introduce the parameter
\begin{equation}\label{7.3}
B'=(k-r+1)b-ra-2(k-r)\gamma .
\end{equation}

\par We find it useful to provide sufficient room in which to work, and therefore explore the contrary situation in 
which $B'\le \nu$. In such circumstances, we find from \eqref{7.3} that
\[
(k-r+1)b-ra\le \nu+2(k-r)\gamma\le 2k\nu.
\]
In view of the definition \eqref{7.2}, we then have
\begin{equation}\label{7.4}
b'-a\le 1+2k\nu/r.
\end{equation}
Moreover, an application of Lemma \ref{lemma6.3} delivers the bound
\[
\rho_a(\xi)^2|\mathfrak f_a(\boldsymbol \alpha;\xi)|^{2R}\le (p^{2b'-2a})^R
\sum_{\substack{\zeta\mmod{(p^{b'},t^2)}\\ \zeta\equiv \xi\mmod{(p^a,t^2)}}}\rho_{b'}(\zeta)^2
|\mathfrak f_{b'}(\boldsymbol \alpha ;\zeta)|^{2R}.
\]
Since we have $R=r(r+1)/2<kr/2$, we deduce from \eqref{7.4}, \eqref{3.16} and \eqref{3.17} that
\[
K_{a,b,c}^{r,\boldsymbol \varphi,\nu}(\boldsymbol{\mathfrak a})\le (p^{2\nu})^{2k^2}
K_{b',b,c}^{r,\boldsymbol \varphi,\nu}(\boldsymbol{\mathfrak a}).
\]
This bound establishes the conclusion of the lemma when $B'\le \nu$.\par

We now turn to address the situation in which
\begin{equation}\label{7.5}
B'>\nu.
\end{equation}
Recall that, by orthogonality and the definition \eqref{3.17}, the mean value 
$K_{a,b,c}^{r,\boldsymbol \varphi,\nu}(\boldsymbol{\mathfrak a};\xi,\eta)$ counts the solutions in 
$\mathfrak Z(t^2)$ of the simultaneous congruences \eqref{3.18}, with the accompanying conditions \eqref{3.19}, 
and with each solution being counted with the weight \eqref{3.20}. In the first phase of our argument, we seek to 
extract from these congruences an underlying additional congruence condition on a subset of the variables that 
may subsequently be exploited.\par

We begin by examining the right hand side of the system \eqref{3.18}. Observe that $\boldsymbol \varphi$ is a 
$p^c$-spaced system, so it follows as in the argument of the proof of Lemma \ref{lemma4.4} leading to \eqref{4.6} 
that there is no loss of generality in supposing that
\[
\varphi_j(u)=u^j+u^{k+1}\bigl( p^c\psi_j(u)+t\chi_j(u)\bigr) \quad (1\le j\le k),
\]
for appropriate polynomials $\psi_j(u),\chi_j(u)\in \mathbb Z[t,u]$. In view of \eqref{3.19}, we may substitute
\[
v_l=p^bu_l+\eta\quad \text{and}\quad w_l=p^bz_l+\eta\quad (1\le l\le s-R),
\]
and then \eqref{3.18} is transformed into the new system
\begin{equation}\label{7.6}
\sum_{i=1}^R\bigl( \varphi_j(x_i)-\varphi_j(y_i)\bigr)\equiv \sum_{l=1}^{s-R}\bigl( \varphi_j(p^bu_l+\eta)-
\varphi_j(p^bz_l+\eta)\bigr) \mmod{(p^B,t^2)}\quad (1\le j\le k).
\end{equation}
Next, as in the argument of the proof of Lemma \ref{lemma4.4} leading from \eqref{4.6} to \eqref{4.7} and 
\eqref{4.8}, we may take appropriate integral linear combinations of the congruences comprising \eqref{7.6} to 
deduce that there exist polynomials $\Upsilon_j(u),\Theta_j(u)\in \mathbb Z[t,u]$ for which
\[
\sum_{i=1}^R\bigl( \Phi_j(x_i-\eta)-\Phi_j(y_i-\eta)\bigr)\equiv 
\sum_{l=1}^{s-R}\bigl( \Phi_j(p^bu_l)-\Phi_j(p^bz_l)\bigr) \mmod{(p^B,t^2)}\quad (1\le j\le k),
\]
in which
\begin{equation}\label{7.7}
\Phi_j(u)=u^j+u^{k+1}\bigl( p^c\Upsilon_j(u)+t\Theta_j(u)\bigr) .
\end{equation}
In particular, the system $\boldsymbol \Phi$ is $p^c$-spaced, and one has
\begin{equation}\label{7.8}
\sum_{i=1}^R \Phi_j(x_i-\eta)\equiv \sum_{i=1}^R\Phi_j(y_i-\eta)\mmod{(p^{jb},p^B,t^2)}\quad (1\le j\le k).
\end{equation}

\par Our next step is to extract from \eqref{7.8} a system related to the Vinogradov systems central to this memoir. 
Of necessity, we recycle variable names despite the risk of confusion. Recall from \eqref{3.20} that in the solutions 
of \eqref{3.19} of interest to us, we have
\[
\mathbf x\equiv \mathbf y\equiv \xi\mmod{(p^a,t^2)}.
\]
Then we may substitute
\begin{equation}\label{7.9}
x_i\equiv p^au_i+\xi\mmod{(p^B,t^2)}\quad \text{and}\quad y_i\equiv p^az_i+\xi\mmod{(p^B,t^2)}
\quad (1\le i\le R).
\end{equation}
We can say more concerning the element $\xi-\eta\in \mathbb Z[t]/(p^B,t^2)$ implicit in \eqref{7.8}. We recall that 
$\xi=\xi_0+\xi_1t$ and $\eta=\eta_0+\eta_1t$ with $\xi_0-\eta_0=\omega p^\gamma$, and that we may suppose 
that $B\ge (k-r+1)b$. Then by ignoring in \eqref{7.8} those congruences with index $j$ satisfying $1\le j\le k-r$, and 
substituting \eqref{7.9} into \eqref{7.8}, we find that for $k-r+1\le j\le k$, one has
\begin{equation}\label{7.10}
\sum_{i=1}^R \Phi_j\bigl( p^au_i+\omega p^\gamma +(\xi_1-\eta_1)t\bigr) \equiv 
\sum_{i=1}^R \Phi_j\bigl( p^az_i+\omega p^\gamma +(\xi_1-\eta_1)t\bigr) \mmod{(p^{(k-r+1)b},t^2)}.
\end{equation}

Just as in the corresponding argument in the proof of \cite[Lemma 7.1]{Woo2019}, the system \eqref{7.10} can be 
analysed by using an argument paralleling that leading from \eqref{4.6} to \eqref{4.7} and \eqref{4.8}. We apply the 
binomial theorem within \eqref{7.7}. Thus, when $1\le l\le r$, we find that for suitable polynomials 
$\Gamma_l(u)\in \mathbb Z[t,u]$, one has
\[
\Phi_{k-r+l}\bigl( p^ay+\omega p^\gamma +(\xi_1-\eta_1)t\bigr)-
\Phi_{k-r+l}\bigl( \omega p^\gamma +(\xi_1-\eta_1)t\bigr) ={\rm H}_l(p^ay),
\]
in which
\begin{equation}\label{7.11}
{\rm H}_l(u)=\sum_{i=1}^r\Omega_{il}\bigl( \omega p^\gamma+(\xi_1-\eta_1)t\bigr)^{k-r+l-i}u^i+
u^{r+1}\Gamma_l(u),
\end{equation}
and the $\mathbb Z[t]$-coefficients $\Omega_{il}$ satisfy the congruences
\[
\Omega_{il}\equiv \binom{k-r+l}{i}\mmod{(p^c,t)}\quad (1\le i,l\le r).
\]
As in the argument of the proof of \cite[Lemma 3.2]{Woo2013}, our hypothesis $p>k$ ensures that the matrix 
$A=(\Omega_{il})_{1\le i,l\le r}$ satisfies the property that $\det(A)\not\equiv 0\mmod{(p,t)}$. It follows that $A$ 
possesses an approximate multiplicative inverse $A_0^{-1}$ modulo $(p,t)$ having coefficients in $\mathbb Z$ with 
the property that $A_0^{-1}A\equiv I_r\mmod{(p,t)}$. We replace $\boldsymbol \Phi$ by 
$A_0^{-1}\boldsymbol \Phi$, and likewise $\boldsymbol \Gamma$ by $A_0^{-1}\boldsymbol \Gamma$. This 
transformation once again corresponds to taking linear combinations of the congruences comprising 
\eqref{7.10}. In this way, we find that there is no loss of generality in supposing that $\boldsymbol{\rm H}$ 
takes the shape defined by \eqref{7.11}, but with $A=(\Omega_{il})_{1\le i,l\le r}$ constrained to satisfy 
$A\equiv I_r\mmod{(p,t)}$. We therefore infer from Lemma \ref{lemma4.2} that this new matrix $A$ possesses a 
multiplicative inverse $A^{-1}$ modulo $(p^{(k-r+1)b},t^2)$ having coefficients in $\mathbb Z[t]/(p^B,t^2)$.\par

We now replace $\boldsymbol \Phi$ by $A^{-1}\boldsymbol \Phi$, and likewise $\boldsymbol \Gamma$ by 
$A^{-1}\boldsymbol \Gamma$. This, yet again, amounts to taking suitable linear combinations of the congruences 
comprising \eqref{7.10}. In this way, one sees that there is no loss of generality in supposing that the coefficient 
matrix $A$ is equal to $I_r$. Since $(\omega,p)=1$, moreover, there exists an integral multiplicative inverse 
$\omega^{-1}$ for $\omega$ modulo $p^{(k-r+1)b}$. This observation, in turn, implies that
\begin{align*}
\bigl( (\omega^{-1})^2(\omega p^\gamma-(\xi_1-\eta_1)t)\bigr) \bigl( \omega p^\gamma+(\xi_1-\eta_1)t\bigr) 
&=\omega^{-2}\bigl( \omega^2p^{2\gamma}-(\xi_1-\eta_1)^2 t^2\bigr) \\
&\equiv p^{2\gamma}\mmod{(p^{(k-r+1)b},t^2)}.
\end{align*}
Consequently, there exist polynomials $\Delta_l\in \mathbb Z[t,u]$ having the property that whenever the system 
\eqref{7.10} is satisfied, then
\[
(p^{2\gamma})^{k-r}\sum_{i=1}^R (p^a)^l\bigl( \Psi_l(u_i)-\Psi_l(z_i)\bigr) \equiv 0
\mmod{(p^{(k-r+1)b},t^2)}\quad (1\le l\le r),
\]
in which
\begin{equation}\label{7.12}
\Psi_l(u)=u^l+p^{a-2(k-r)\gamma}\Delta_l(u).
\end{equation}
Notice here that, since we may suppose from \eqref{7.1} that $a\ge \delta \theta$ and $\gamma \le \nu$, our 
hypotheses concerning the heierarchy \eqref{6.2} combines with \eqref{6.6} and \eqref{6.20} to give
\[
2k\gamma\le 2k\lceil 8\varepsilon H\Lambda^{-1}\rceil <\delta^2\lceil \mu H\rceil \le \delta a,
\]
whence
\[
a-2(k-r)\gamma>(1-\delta)a\ge \delta(1-\delta)\mu H>\tau B.
\]
We therefore perceive that the exponent of $p$ on the right hand side of \eqref{7.12} is positive, and moreover 
the system $\boldsymbol \Psi$ is $p^c$-spaced for some $c>\tau(k-r+1)b$.\par

The discussion thus far shows that the $\mathfrak Z(t^2)$-solutions $\mathbf x=p^a\mathbf u+\xi$, 
$\mathbf y=p^a\mathbf z+\xi$, $\mathbf v$, $\mathbf w$ of the original system of congruences \eqref{3.18}, 
counted with weight \eqref{3.20} associated with the definition \eqref{3.17} of 
$K_{a,b,c}^{r,\boldsymbol \varphi,\nu}(\boldsymbol{\mathfrak a};\xi,\eta)$, are constrained by the additional 
system of congruences
\begin{equation}\label{7.13}
\sum_{i=1}^R \Psi_l(u_i)\equiv \sum_{i=1}^R\Psi_l(z_i)\mmod{(p^{B'},t^2)}\quad (1\le l\le r),
\end{equation}
in which $B'$ is defined via \eqref{7.3}. In the final phase of our argument, our goal is to apply the grand inductive 
hypothesis, by which we mean the conclusion of Theorem \ref{theorem3.2} with $k$ replaced by $r$, so as to 
extract congruence information from the relations \eqref{7.13}.\par

Recall that in view of \eqref{7.5}, one has $B'>\nu$, and write $H'=\lceil B'/r\rceil$. In view of \eqref{7.2} and 
\eqref{7.3}, one has
\begin{equation}\label{7.14}
b'-H'\le a+1+2(k-r)\gamma/r.
\end{equation}
From \eqref{6.6}, moreover, we have $B'>\nu\ge 8\varepsilon H\Lambda^{-1}$. Thus, our hypotheses concerning 
$B$ and the heierarchy \eqref{6.2} ensure that $B'$ is sufficiently large in terms of $k$, $\Lambda$, $\mu$, 
$\tau$, $\delta$ and $\varepsilon$. We construe the definition \eqref{3.4} of $\mathfrak f_h(\boldsymbol \alpha;\xi)$ 
in the form \eqref{4.2}, so that
\begin{equation}\label{7.15}
\mathfrak f_a(\boldsymbol \alpha ;\xi)=\rho_a(\xi)^{-1}\sum_{y\in \mathfrak Z(t^2)}\mathfrak c_y(\boldsymbol \alpha),
\end{equation}
where
\[
\mathfrak c_y(\boldsymbol \alpha)=\mathfrak a_{p^ay+\xi}E(\psi(p^ay+\xi;\boldsymbol \alpha)),
\]
and $\psi(n;\boldsymbol \alpha)$ is defined by means of \eqref{3.3}. Note that, in view of \eqref{3.2}, the latter 
definition yields the relation
\begin{equation}\label{7.16}
\rho_0(1;\boldsymbol{\mathfrak c})=\|\boldsymbol{\mathfrak c}\|_0^{1/2}\cdot \| \boldsymbol{\mathfrak c}\|_\infty 
=\|\boldsymbol{\mathfrak a}_a(\xi)\|_0^{1/2}\cdot \|\boldsymbol{\mathfrak a}_a(\xi)\|_\infty=\rho_a(\xi).
\end{equation}
We also define the exponential sum
\begin{equation}\label{7.17}
\mathfrak g_{\boldsymbol{\mathfrak c}}(\boldsymbol \alpha,\boldsymbol \beta)
=\rho_0(1;\boldsymbol{\mathfrak c})^{-1}\sum_{y\in \mathfrak Z(t^2)}\mathfrak c_y(\boldsymbol \alpha)
E(\beta_1\Psi_1(y)+\ldots +\beta_r\Psi_r(y)),
\end{equation}
and the mean value
\[
J(\boldsymbol \alpha)=\oint_{p^{B'}}
|\mathfrak g_{\boldsymbol{\mathfrak c}}(\boldsymbol \alpha,\boldsymbol \beta)|^{2R}\d\boldsymbol \beta .
\]

\par By orthogonality, the mean value $J(\boldsymbol \alpha)$ counts the $\mathfrak Z(t^2)$-solutions of the 
system of congruences \eqref{7.13}, with each solution $\mathbf u,\mathbf z$ being counted with weight
\begin{align*}
\rho_0(1;\boldsymbol{\mathfrak c})^{-2R}&\prod_{i=1}^R \mathfrak c_{u_i}(\boldsymbol \alpha)
{\overline{\mathfrak c}}_{z_i}(\boldsymbol \alpha)\\
&=\rho_a(\xi)^{-2R}\biggl( \prod_{i=1}^R
\mathfrak a_{p^au_i+\xi}{\overline{\mathfrak a}}_{p^az_i+\xi}\biggr) E\biggl( \sum_{i=1}^R\bigl( 
\psi(p^au_i+\xi;\boldsymbol \alpha)-\psi(p^az_i+\xi;\boldsymbol \alpha)\bigr)\biggr) .
\end{align*}
The initial phase of our discussion revealed that the conditions \eqref{7.13} on $\mathbf x=p^a\mathbf u+\xi$ and 
$\mathbf y=p^a\mathbf z+\xi$ are implicitly implied by \eqref{3.18}. Consequently, since \eqref{7.15}, \eqref{7.16} 
and \eqref{7.17} yield the relation
\[
\mathfrak g_{\boldsymbol{\mathfrak c}}(\boldsymbol \alpha ,\mathbf 0)=\mathfrak f_a(\boldsymbol \alpha;\xi),
\]
we deduce that
\begin{equation}\label{7.18}
\oint_{p^B}|\mathfrak f_a(\boldsymbol \alpha ;\xi)^{2R}\mathfrak f_b(\boldsymbol \alpha;\eta)^{2s-2R}|
\d\boldsymbol \alpha =\oint_{p^B}\oint_{p^{B'}}|\mathfrak g_{\boldsymbol{\mathfrak c}}(\boldsymbol \alpha, 
\boldsymbol \beta )^{2R}\mathfrak f_b(\boldsymbol \alpha ;\eta)^{2s-2R}|\d \boldsymbol \beta \d\boldsymbol 
\alpha .
\end{equation}
Thus, we arrive at the relation
\begin{equation}\label{7.19}
K_{a,b,c}^{r,\boldsymbol \varphi,\nu}(\boldsymbol{\mathfrak a};\xi,\eta)
=\oint_{p^B}J(\boldsymbol \alpha)|\mathfrak f_b(\boldsymbol \alpha ;\eta)|^{2s-2R}\d\boldsymbol \alpha .
\end{equation}
We emphasise that the inner integral over $\boldsymbol \beta$ on the right hand side of \eqref{7.18} is, on a 
fundamental level, content-free, since the integral over $\boldsymbol \alpha$ restricts the variables underlying 
the factor $|\mathfrak g_{\boldsymbol{\mathfrak c}}(\boldsymbol \alpha,\boldsymbol \beta )|^{2R}$ so as to satisfy 
the system of congruences \eqref{7.13} automatically.\par

In the next phase of our argument, we finally apply our inductive hypothesis. The key observation is that the mean 
value $J(\boldsymbol \alpha)$ is amenable to Corollary \ref{corollary3.3} to Theorem \ref{theorem3.2}. The system 
$\boldsymbol \Psi$ is $p^c$-spaced for some $c\ge \tau B'$, and $B'$ is sufficiently large in terms of $k$, $\tau$ 
and $\varepsilon^2$.  In view of Hypothesis \ref{hypothesis3.2} for the collection of weights 
$(\boldsymbol{\mathfrak a})$, moreover, we may assume that the collection of weights 
$(\boldsymbol{\mathfrak c}_y(\boldsymbol \alpha))$ satisfies Hypothesis \ref{hypothesis3.2}.  
Hence, we have
\begin{align}
J(\boldsymbol \alpha)&=U_{R,r}^{B'}(\boldsymbol{\mathfrak c})\ll p^{B'\varepsilon^2}
U_{R,r}^{B',H'}(\boldsymbol{\mathfrak c}) \notag\\
&=p^{B'\varepsilon^2}\rho_0(1;\boldsymbol{\mathfrak c})^{-2}\sum_{\zeta \mmod{(p^{H'},t^2)}}
\rho_{H'}(\zeta;\boldsymbol{\mathfrak c})^2\oint_{p^{B'}}|\mathfrak g_{\boldsymbol{\mathfrak c}'}(\boldsymbol 
\alpha,\boldsymbol \beta)|^{2R}\d\boldsymbol \beta ,\label{7.20}
\end{align} 
 in which $\boldsymbol{\mathfrak c}'=\boldsymbol{\mathfrak c}'(\boldsymbol \alpha)$ is defined by putting
\[
\mathfrak c'_n(\boldsymbol \alpha)=\begin{cases} \mathfrak c_n(\boldsymbol \alpha),&\text{when $n\equiv \zeta 
\mmod{(p^{H'},t^2)}$,}\\
0,&\text{otherwise.}\end{cases}
\]
Note that $\rho_0(1;\boldsymbol{\mathfrak c})=\rho_a(\xi)$. Moreover, on putting $\kappa=p^a\zeta+\xi$, we 
find that
\begin{align*}
\rho_{H'}(\zeta;\boldsymbol{\mathfrak c})&=\| (\mathfrak a_{p^ay+\xi}: y\equiv \zeta\mmod{(p^{H'},t^2)}) \|_0^{1/2}
\cdot \|(\mathfrak a_{p^ay+\xi}: y\equiv \zeta\mmod{(p^{a+H'},t^2)})\|_\infty \\
&=\rho_{a+H'}(\kappa).
\end{align*}
Thus, on substituting \eqref{7.20} into \eqref{7.19}, we infer that
\begin{align}\label{7.21}
\rho_a(\xi)^2&K_{a,b,c}^{r,\boldsymbol \varphi,\nu}(\boldsymbol{\mathfrak a};\xi,\eta) \notag \\
&\ll p^{B'\varepsilon^2}
\sum_{\substack{\kappa\mmod{(p^{a+H'},t^2)}\\ \kappa\equiv \xi\mmod{(p^a,t^2)}}}\rho_{a+H'}(\kappa)^2
\oint_{p^B}\oint_{p^{B'}}|\mathfrak g_{\boldsymbol{\mathfrak c}'}(\boldsymbol \alpha, 
\boldsymbol \beta )^{2R}\mathfrak f_b(\boldsymbol \alpha ;\eta)^{2s-2R}|\d \boldsymbol \beta \d\boldsymbol 
\alpha .
\end{align}

\par Having extracted our strong congruence constraint on a subset of the implicit variables, our next task is to 
reverse course so as to relate the mean value on the right hand side of \eqref{7.21} to simpler mean values closer in 
spirit to the target $K_{b',b,c}^{r,\boldsymbol \varphi,\nu}(\boldsymbol{\mathfrak a})$. We begin by observing that 
the mean value
\[
\oint_{p^{B'}}|\mathfrak g_{\boldsymbol{\mathfrak c}'}(\boldsymbol \alpha, \boldsymbol \beta )|^{2R}
\d \boldsymbol \beta
\]
counts the $\mathfrak Z(t^2)$-solutions of the system of equations \eqref{7.13}, with each solution $\mathbf u$, 
$\mathbf z$ being counted with weight
\begin{align*}
\rho_0(1;\boldsymbol{\mathfrak c}')^{-2R}&\prod_{i=1}^R \mathfrak c'_{u_i}(\boldsymbol \alpha)
{\overline{\mathfrak c}}'_{z_i}(\boldsymbol \alpha)\\
&=\rho_{a+H'}(\kappa )^{-2R}\biggl( \prod_{i=1}^R
\mathfrak a_{p^au_i+\xi}{\overline{\mathfrak a}}_{p^az_i+\xi}\biggr) E\biggl( \sum_{i=1}^R\bigl( 
\psi(p^au_i+\xi;\boldsymbol \alpha)-\psi(p^az_i+\xi;\boldsymbol \alpha)\bigr)\biggr) ,
\end{align*}
but now subject to the constraint $\mathbf u\equiv \mathbf z\equiv \zeta \mmod{p^{H'}}$. Note that the 
conditions \eqref{7.13} on $\mathbf x=p^a\mathbf u+\xi$ and $\mathbf y=p^a\mathbf z+\xi$ are again 
automatically implied by \eqref{3.18}. Consequently, since we have 
$\mathfrak g_{\boldsymbol{\mathfrak c}'}(\boldsymbol \alpha, \mathbf 0)=
\mathfrak f_{a+H'}(\boldsymbol \alpha;\kappa)$, we perceive that
\begin{equation}\label{7.22}
\oint_{p^B}\oint_{p^{B'}}|\mathfrak g_{\boldsymbol{\mathfrak c}'}(\boldsymbol \alpha, 
\boldsymbol \beta )^{2R}\mathfrak f_b(\boldsymbol \alpha ;\eta)^{2s-2R}|
\d \boldsymbol \beta \d\boldsymbol \alpha
=\oint_{p^B}|\mathfrak f_{a+H'}(\boldsymbol \alpha;\kappa)^{2R}\mathfrak f_b(\boldsymbol \alpha ;\eta)^{2s-2R}|
\d\boldsymbol \alpha .
\end{equation}
By reference to Lemma \ref{lemma6.3}, moreover, we obtain the upper bound
\begin{equation}\label{7.23}
\rho_{a+H'}(\kappa)^2|\mathfrak f_{a+H'}(\boldsymbol \alpha;\kappa)|^{2R}\le (p^{2b'-2(a+H')})^R
\sum_{\substack{\xi'\mmod{(p^{b'},t^2)}\\ \xi'\equiv \kappa\mmod{(p^{a+H'},t^2)}}}\rho_{b'}(\xi')
|\mathfrak f_{b'}(\boldsymbol \alpha;\xi')|^{2R}.
\end{equation}
Therefore, by substituting \eqref{7.23} into \eqref{7.22}, and the resulting bound into \eqref{7.21}, we arrive at the 
estimate
\[
\rho_a(\xi)^2K_{a,b,c}^{r,\boldsymbol \varphi,\nu}(\boldsymbol{\mathfrak a};\xi,\eta)\ll 
p^{B'\varepsilon^2+2R(b'-(a+H'))}\sum_{\substack{\xi'\mmod{(p^{b'},t^2)}\\ \xi'\equiv \xi\mmod{(p^a,t^2)}}}
\rho_{b'}(\xi')^2K_{b',b,c}^{r,\boldsymbol \varphi,\nu}(\boldsymbol{\mathfrak a};\xi',\eta). 
\]
This upper bound may now, in turn, be substituted into the definition \eqref{3.16}, and thus we find that
\begin{align*}
\rho_0(1)^4&K_{a,b,c}^{r,\boldsymbol \varphi,\nu}(\boldsymbol{\mathfrak a})\\
&\ll p^{B'\varepsilon^2+2R(b'-(a+H'))}
\sum_{\xi'\mmod{(p^{b'},t^2)}}\sum_{\substack{\eta\mmod{(p^{b'},t^2)}\\ \eta\not\equiv \xi'\mmod{(p^\nu,t)}}}
\rho_{b'}(\xi')^2\rho_b(\eta)^2K_{b',b,c}^{r,\boldsymbol \varphi,\nu}(\boldsymbol{\mathfrak a};\xi',\eta),
\end{align*}
and hence
\begin{equation}\label{7.24}
K_{a,b,c}^{r,\boldsymbol \varphi,\nu}(\boldsymbol{\mathfrak a})\ll p^{B'\varepsilon^2+2R(b'-(a+H'))}
K_{b',b,c}^{r,\boldsymbol \varphi,\nu}(\boldsymbol{\mathfrak a}).
\end{equation}

\par The proof of the lemma is now swiftly made complete. We observe that, as a consequence of \eqref{7.14}, one 
has
\[
2R(b'-a-H')=r(r+1)(b'-a-H')\le (r+1)(r+2(k-r)\gamma)<2k^2\nu.
\]
Consequently, on recalling the definition \eqref{6.6} and the heierarchy of small quantities \eqref{6.2}, we find that
\[
B'\varepsilon^2+2R(b'-a-H')<\nu+2k^2\nu<3k^2\nu.
\]
The conclusion of the lemma therefore follows from \eqref{7.24} in this situation in which we assume that $B'>\nu$.
\end{proof}

\section{The iterative process}
The work of the previous section, in common with \cite[\S7]{Woo2019} in its own setting, shows how to generate 
powerful congruence constraints on the variables underlying the mean value 
$K_{a,b,c}^{r,\boldsymbol \varphi,\nu}(\boldsymbol{\mathfrak a})$. In all essentials, the conclusion of our own 
Lemma \ref{lemma7.1} provides a relation more or less identical to that of \cite[Lemma 7.1]{Woo2019}, save that we 
work in $\mathbb Z[t]/t^2\mathbb Z[t]$ instead of in $\mathbb Z$. The subtle differences in the definitions of 
$K_{a,b,c}^{r,\boldsymbol \varphi,\nu}(\boldsymbol{\mathfrak a})$ in the respective situations do not obstruct us 
from drawing the same conclusions. We therefore provide only a sketch of the argument, noting the very few 
differences encountered as the argument is pursued. At this point, the parameters $\boldsymbol \varphi$, $\nu$ 
and $c$ have ceased to possess any particular significance, and henceforth we may safely omit their mention in our 
various notations. Also, throughout this section, we suppose that $a$, $b$ and $r$ are integers satisfying \eqref{7.1}, 
and we define $b'$ via \eqref{7.2}. Our first step is to interchange the roles of $b'$ and $b$ in the conclusion of 
Lemma \ref{lemma7.1}

\begin{lemma}\label{lemma8.1}
When $r\ge 2$, one has
\[
K_{a,b}^r(\boldsymbol{\mathfrak a})\ll (p^{2\nu})^{2k^2}K_{b,b'}^{k-r}(\boldsymbol{\mathfrak a})^{\frac{1}{k-r+1}}
K_{b',b}^{r-1}(\boldsymbol{\mathfrak a})^{\frac {k-r}{k-r+1}}.
\]
Meanwhile, when $r=1$, one instead has
\[
K_{a,b}^1(\boldsymbol{\mathfrak a})\ll (p^{2\nu})^{2k^2}K_{b,kb}^{k-1}(\boldsymbol{\mathfrak a})^{1/k}
U_{s,k}^{B,b}(\boldsymbol{\mathfrak a})^{1-1/k}.
\]
\end{lemma}

\begin{proof} Both conclusions follow by pursuing the argument of the proof of \cite[Lemma 8.1]{Woo2019}, 
mutatis mutandis. The factor $p^{k^2\nu}$ in the latter is here replaced by the factor $(p^{2\nu})^{2k^2}$ from 
the conclusion of  Lemma \ref{lemma7.1} throughout. Moreover, the condition $\xi\not \equiv \eta \mmod{p^\nu}$ 
in \cite{Woo2019} is replaced in the present context by $\xi\not \equiv \eta \mmod{(p^\nu,t)}$. In all other respects, 
the argument follows almost word for word as in the proof of \cite[Lemma 8.1]{Woo2019}.
\end{proof}

By employing the definition \eqref{3.21}, the conclusion of Lemma \ref{lemma8.1} may be converted into a 
normalised analogue.

\begin{lemma}\label{lemma8.2} Suppose that $b\le (1-\delta)B/k$. Then, when $r\ge 2$, one has
\[
[[ K_{a,b}^r(\boldsymbol{\mathfrak a})]]_\Lambda \ll (p^{2\nu})^{2k^2}
[[ K_{b,b'}^{k-r}(\boldsymbol{\mathfrak a})]]_\Lambda^{\frac{1}{k-r+1}}
[[ K_{b',b}^{r-1}(\boldsymbol{\mathfrak a})]]_\Lambda^{1-1/r}.
\]
Meanwhile, when $r=1$, one instead has
\begin{equation}\label{8.1}
[[ K_{a,b}^1(\boldsymbol{\mathfrak a})]]_\Lambda \ll (p^{2\nu})^{4k^2}
[[ K_{b,kb}^{k-1}(\boldsymbol{\mathfrak a})]]_\Lambda^{1/k}\bigl( p^{-2b}\bigr)^{\Lambda (1-1/k)}.
\end{equation}
\end{lemma}

\begin{proof} The proof follows that of \cite[Lemma 8.2]{Woo2019}, mutatis mutandis. In our setting, the upper 
bound \eqref{6.5} contains an extra factor of $2$ in the exponent relative to its analogue in \cite{Woo2019}, and 
hence one obtains an extra factor $2$ in the exponent of $p$ on the right hand side of \eqref{8.1} in the final factor.
\end{proof}

At this point, we streamline our notation further by writing
\[
{\widetilde K}_{a,b}^r=[[ K_{a,b}^r(\boldsymbol{\mathfrak a})]]_\Lambda .
\]
Also, when $1\le j\le k-1$, we write
\[
\rho_j=\frac{j}{k-j+1}\quad \text{and}\quad b_j=\biggl\lceil \frac{b}{\rho_j}\biggr\rceil .
\]

\begin{lemma}\label{lemma8.3}
Suppose that
\begin{equation}\label{8.2}
1\le r\le k-1,\quad a\ge \delta \theta,\quad b\ge k^2\delta \theta\quad \text{and}\quad ra\le (k-r+1)b.
\end{equation}
Then, whenever $k^2b\le (1-\delta )B$, there exist integers $a'$, $b'$, $r'$, and there exists a positive number 
$\rho$, having the property that
\[
{\widetilde K}_{a,b}^r\ll \bigl( {\widetilde K}_{a',b'}^{r'}\bigr)^{\rho}\bigl( p^{-2b}\bigr)^{\Lambda /(2k)},
\]
with
\[
a'\ge \delta \theta,\quad b'\ge k^2\delta \theta,\quad r'a'\le (k-r'+1)b',\quad 1\le r'\le k-1,
\]
\[
0<\rho<(1-1/k)^2,\quad (1+2/k)b\le b'\le k^2b,\quad b'=\biggl\lceil \frac{r'+1}{k-r'}\biggr\rceil ,\quad 
\rho b'\ge b.
\]
\end{lemma}

\begin{proof} Under the hypothesis \eqref{8.2}, the argument of the proof of \cite[Lemma 9.1]{Woo2019} applies, 
essentially verbatim, to show that whenever $kb\le (1-\delta)B$, then
\[
{\widetilde K}_{a,b}^r\ll p^{4(r+1)k^2\nu}\bigl( p^{-2b}\bigr)^{(1-1/k)\Lambda/r} \prod_{j=1}^r \bigl( 
{\widetilde K}_{b,b_j}^{k-j}\bigr)^{\rho_j/r}.
\]
Under the same hypotheses, it follows from the argument of the proof of \cite[Lemma 9.2]{Woo2019}, again applied 
almost verbatim, that there exists an integer $r'$ with $1\le r'\le r$ having the property that
\[
{\widetilde K}_{a,b}^r\ll \bigl( {\widetilde K}_{b,b_{r'}}^{k-r'}\bigr)^{\rho_{r'}}\bigl( p^{-2b}\bigr)^{\Lambda/(2k)}.
\]
The conclusion of the present lemma then follows from this estimate by applying the argument of the proof of 
\cite[Lemma 9.3]{Woo2019}, without any essential alteration.
\end{proof}

\section{The proof of Theorem \ref{theorem3.2}}
The conclusion of Lemma \ref{lemma8.3} facilitates the same kind of $p$-adic concentration argument as was made 
available by \cite[Lemma 9.3]{Woo2019}. Indeed, our proof of Theorem \ref{theorem3.2} may be engineered 
essentially verbatim as in the corresponding argument of \cite[\S10]{Woo2019}. For this reason, we shall be 
extremely economical on details, since it suffices to sketch out differences only.

\begin{proof}[The proof of Theorem \ref{theorem3.2}] In view of the conclusion of Lemma \ref{lemma5.1}, we may 
suppose throughout this section that $k$ is a natural number with $k\ge 2$. In addition, we put $s=k(k+1)/2$. We 
may suppose that the conclusion of Theorem \ref{theorem3.2} has already been established for exponents smaller 
than $k$. By way of seeking a contradiction, we suppose that $\lambda(s,k)=\Lambda$ with $\Lambda>0$. We 
initiate our iterative method by appealing to Lemma \ref{lemma6.4}, which shows that
\[
U_{s,k}^B(\boldsymbol{\mathfrak a})\ll p^{2s\theta}K_{\theta,\theta}^1(\boldsymbol{\mathfrak a}).
\]
In view of the assumption \eqref{6.4}, we deduce from \eqref{3.21} that
\[
[[ K_{\theta,\theta}^1(\boldsymbol{\mathfrak a})]]_\Lambda \gg \frac{p^{-2s\theta}
U_{s,k}^B(\boldsymbol{\mathfrak a})}{(p^{2H})^\Lambda U_{s,k}^{B,H}(\boldsymbol{\mathfrak a})}\ge 
p^{-2s\theta -2H\varepsilon}.
\]
The heierarchy \eqref{6.2} consequently combines with \eqref{6.20} to reveal that
\begin{equation}\label{9.1}
{\widetilde K}_{\theta,\theta}^1\gg p^{-4s\theta}.
\end{equation}

\par Next, we apply Lemma \ref{lemma8.3} repeatedly. Put
\begin{equation}\label{9.2}
N=\lceil 16sk/\Lambda\rceil .
\end{equation}
Our assumption that $B$, and hence also $H$, is sufficiently large in terms of our heierarchy of parameters 
\eqref{6.2} ensures that $2k^{2N+2}\theta<H$. We cliam that the sequences $(a_n)$, $(b_n)$, $(r_n)$, $(\rho_n)$ may 
be defined for $0\le n\le N$ in such a manner that
\begin{align*}
1\le r_n\le k-1,&\quad k^2\delta \theta\le b_n\le k^{2n+2}\theta ,\\
\delta \theta\le a_n\le &\,(k-r_n+1)b_n/r_n,\\
0<\rho_n<(1-1/k)^2,&\quad \rho_nb_n\ge b_{n-1}\quad (n\ge 1),
\end{align*}
and so that
\begin{equation}\label{9.3}
{\widetilde K}_{\theta,\theta}^1\ll \bigl( {\widetilde K}_{a_n,b_n}^{r_n}\bigr)^{\rho_1\cdots \rho_n}
\bigl( p^{-\Lambda /(2k)}\bigr)^{2nb_0}.
\end{equation}
We initiate these sequences by putting $a_0=b_0=\theta$ and $\rho_0=r_0=1$, so that \eqref{9.3} is immediate in 
the case $n=0$.\par

The confirmation of this claim is essentially a numerical check that follows via an argument identical to that bridging 
\cite[equation (10.6)]{Woo2019} and \cite[equation (10.9)]{Woo2019}, mutatis mutandis. From here, by combining 
\eqref{9.1} and \eqref{9.3}, we find that
\begin{equation}\label{9.4}
p^{-4s\theta}\ll {\widetilde K}_{\theta,\theta}^1\ll \bigl( {\widetilde K}_{a_N,b_N}^{r_N}\bigr)^\rho 
\bigl( p^{-\Lambda/(2k)}\bigr)^{2N\theta},
\end{equation}
where $\rho=\rho_1\cdots \rho_N<1$. On recalling the definition \eqref{3.21} and the assumption that 
$\lambda(s,k)=\Lambda$, the conclusion of Lemma \ref{lemma4.5} shows that
\[
[[ K_{a_N,b_N}^{r_N}]]_\Lambda \ll p^{2H\varepsilon}.
\]
Then, again employing the properties of the heierarchy \eqref{6.2}, we may suppose that 
${\widetilde K}_{a_N,b_N}^{r_N}\ll p^\theta$, and thus \eqref{9.4} delivers the bound
\[
p^{-4s\theta}\ll p^{\theta (1-N\Lambda /k)},
\]
and hence
\begin{equation}\label{9.5}
(p^\theta)^{8s}\gg (p^\theta)^{N\Lambda/k}.
\end{equation}
The definition \eqref{6.20} of $\theta$ shows that $p^\theta$ is sufficiently large in terms of $s$, $k$ and 
$\Lambda$. Hence, the upper bound \eqref{9.5} can hold only when $8s\ge N\Lambda/k$. In view of \eqref{9.2}, we 
therefore obtain the bound
\[
\Lambda \le 8sk/N\le \Lambda/2,
\]
yielding a contradiction to the assumption \eqref{6.1} that $\Lambda>0$.\par

This discussion forces us to conclude that $\Lambda$ cannot be positive, and thus $\lambda(s,k)\le 0$. This 
completes the proof of Theorem \ref{theorem3.2} for the exponent $k$, and the theorem therefore follows in full by 
induction on $k$.
\end{proof}

\begin{proof}[The proof of Corollary \ref{corollary3.3}] The conclusion of Corollary \ref{corollary3.3} follows from 
Theorem \ref{theorem3.2} by simply interpreting the definitions \eqref{3.12} and \eqref{3.13} of 
$\lambda^*(s,\theta;\tau)$ and $\lambda (s,\theta)$. Thus, since $\lambda (s,k)\le 0$, we may suppose that for every 
$\varepsilon>0$, whenever $\tau>0$ is sufficiently small, one has $\lambda^*(s,k;\tau)<\varepsilon$. But then, for 
each $\varepsilon>0$, one has
\[
U_{s,k}^B(\boldsymbol{\mathfrak a})\ll p^{(\lambda^*(s,k;\tau)+\varepsilon)H}
U_{s,k}^{B,H}(\boldsymbol{\mathfrak a})
\]
for all coefficient weights $(\boldsymbol{\mathfrak a})\in \mathbb D$, all $\boldsymbol \varphi\in \Phi_\tau(B)$ and 
all large enough values of $B$. The corresponding conclusion for the situation with coefficient weights 
$({\mathbf 0})$ all equal to zero is of course, trivial.
\end{proof}

\section{Mean values of exponential sums}
Our goal in this section is to establish the main conclusions from the introduction by application of Theorem 
\ref{theorem3.2} and its corollary. In pursuit of this objective, we follow closely the path laid down in the classical 
case associated with $\mathbb Z$ described in \cite[\S12]{Woo2019}.

\begin{proof}[The proof of Theorem \ref{theorem1.1}] Let $\varphi_j\in \mathbb Z[t,u]$ $(1\le j\le k)$ be 
polynomials with the Wronskian $W(u;\boldsymbol \varphi)$ non-vanishing modulo $t$, and let $\varepsilon>0$ be 
a small positive number. Let $\mathcal Z$ denote the set of zeros of the congruence 
\[
W(u;\boldsymbol \varphi)\equiv 0\mmod{t},
\]
with $u\in \mathbb Z$, and let $X$ be a positive integer sufficiently large in terms of $\boldsymbol \varphi$, $s$, 
$k$ and $\varepsilon$. We note that the non-vanishing of $W(u;\boldsymbol \varphi)$ modulo $t$ implies that the 
expression $W(u;\boldsymbol \varphi)\mmod{t}$ is a non-zero polynomial having degree bounded in terms of 
$\boldsymbol \varphi$ and $k$. Hence
\begin{equation}\label{10.1}
\text{card}(\mathcal Z)\le \text{deg}\bigl( W(u;\boldsymbol \varphi)\mmod{t}\bigr) \ll 1.
\end{equation}
We suppose that $s=k(k+1)/2$. Finally, when $(\boldsymbol{\mathfrak a})$ is a set of complex weights 
with $|\mathfrak a_{n_0+n_1t}|=1$ for $|n_0|\le X$ and $|n_1|\le X$, and $\mathfrak a_{n_0+n_1t}=0$ otherwise, we 
write
\[
F(\boldsymbol \alpha;X)=(2X+1)^{-1}\sum_{|n_0|\le X}\sum_{|n_1|\le X}\mathfrak 
a_{n_0+n_1t}E\bigl( \psi(n_0+n_1t;\boldsymbol \alpha)\bigr) 
\]
and
\[
F_0(\boldsymbol \alpha;X)=(2X+1)^{-1}\sum_{\substack{|n_0|\le X\\ n_0\not \in \mathcal Z}}\sum_{|n_1|\le X}
\mathfrak a_{n_0+n_1t}E\bigl( \psi(n_0+n_1t;\boldsymbol \alpha)\bigr) ,
\]
where $\psi(n;\boldsymbol \alpha)$ is defined as in \eqref{3.3}. We frequently find it convenient to abbreviate 
$n_0+n_1t$ simply to $n$, with similar conventions applied with related notation.\par

As a consequence of the triangle inequality, one has
\[
\biggl| \sum_{\substack{|n_0|\le X\\ n_0\in \mathcal Z}}\sum_{|n_1|\le X}\mathfrak a_{n_0+n_1t}
E\bigl( \psi(n_0+n_1t; \boldsymbol \alpha)\bigr)\biggr| \le (2X+1)\text{card}(\mathcal Z),
\]
and thus, it follows from \eqref{10.1} that
\[
|F(\boldsymbol \alpha;X)|\ll 1+|F_0(\boldsymbol \alpha ;X)|.
\]
Thus, one has
\begin{equation}\label{10.2}
\oint |F(\boldsymbol \alpha;X)|^{2s}\d\boldsymbol \alpha \ll 1+\oint |F_0(\boldsymbol \alpha ;X)|^{2s}
\d\boldsymbol \alpha .
\end{equation}
In this way, it is apparent that it suffices to restrict attention to mean values in which the underlying variables 
possess non-vanishing Wronskians modulo $t$. By orthogonality, the mean value on the right hand side of 
\eqref{10.2} counts the number of integral solutions of the system of equations \eqref{1.3} with 
$|u_i|,|v_i|\le X$ and $u_i\not\in \mathcal Z\mmod{t}$ $(1\le i\le 2s)$, and with each solution being counted with 
weight
\begin{equation}\label{10.3}
(2X+1)^{-2s}\prod_{i=1}^s\mathfrak a_{u_i+v_it}{\overline{\mathfrak a}}_{u_{s+i}+v_{s+i}t}.
\end{equation}

\par Our next step is to shift from non-vanishing of Wronskians modulo $t$ to corresponding non-vanishing modulo 
$(p,t)$, for a suitable prime number $p$. When $n\in \mathbb Z[t]$, denote by $W_t(n;\boldsymbol \varphi)$ the 
integer congruent modulo $t$ to $W(n;\boldsymbol \varphi)$. Given a solution $\mathbf u,\mathbf v$ of the system 
\eqref{1.3} of the type just described, the integer
\[
\Xi (\mathbf u)=\prod_{i=1}^{2s}W_t(u_i;\boldsymbol \varphi)
\]
is non-zero. Moreover, one has $|\Xi(\mathbf u)|\le CX^D$, for some $C>0$ depending at most on $s$, $k$ and the 
coefficients of $\boldsymbol \varphi$, and $D$ a positive integer with
\[
D\le 2s\sum_{j=1}^k \text{deg}(\varphi_j).
\] 
Let $\mathcal P$ denote the set of primes $p$ with $(\log X)^2<p\le 3(\log X)^2$. Thus, when $X$ is large, it is a 
consequence of the prime number theorem that
\[
\prod_{p\in \mathcal P}p>(\log X)^{(\log X)^2/\log \log X}>CX^D.
\]
We therefore deduce that for each solution $\mathbf u,\mathbf v$ of the system \eqref{1.3} counted by the integral on 
the right hand side of \eqref{10.2}, there exists a prime $p\in \mathcal P$ with
\[
\prod_{i=1}^{2s}W_t(u_i;\boldsymbol \varphi)\not\equiv 0\mmod{p}.
\]
In particular, one has
\[
\oint |F_0(\boldsymbol \alpha;X)|^{2s}\d\boldsymbol \alpha \le \sum_{p\in \mathcal P}
\oint |F_p(\boldsymbol \alpha;X)|^{2s}\d\boldsymbol \alpha ,
\]
where
\begin{equation}\label{10.4}
F_p(\boldsymbol \alpha;X)=(2X+1)^{-1}\sum_{W(n;\boldsymbol \varphi)\not\equiv 0\mmod{(p,t)}}
|\mathfrak a_n|E\bigl( \psi(n;\boldsymbol \alpha)\bigr) .
\end{equation}
Thus, we conclude from \eqref{10.2} that
\begin{equation}\label{10.5}
\oint |F(\boldsymbol \alpha;X)|^{2s}\d\boldsymbol \alpha \ll 1+(\log X)^2\max_{p\in \mathcal P}\oint 
|F_p(\boldsymbol \alpha;X)|^{2s}\d\boldsymbol \alpha.
\end{equation}

\par We seek to establish that for each fixed prime $p\in \mathcal P$, one has
\begin{equation}\label{10.6}
\oint |F_p(\boldsymbol \alpha;X)|^{2s}\d\boldsymbol \alpha \ll X^{\varepsilon/2}.
\end{equation}
Granted this estimate, it follows from \eqref{10.5} that
\[
\oint |F(\boldsymbol \alpha;X)|^{2s}\d\boldsymbol \alpha \ll 1+X^{\varepsilon/2}(\log X)^2\ll X^\varepsilon ,
\]
whence
\begin{equation}\label{10.7}
\oint \biggl| \sum_{|n_0|\le X}\sum_{|n_1|\le X}\mathfrak a_{n_0+n_1t}E\bigl( \psi(n_0+n_1t;\boldsymbol \alpha)\bigr) 
\biggr|^{2s}\d\boldsymbol \alpha \ll X^{2s+\varepsilon}.
\end{equation}
By specialising the coefficient weights $\mathfrak a_{n_0+n_1t}$ to be $1$ for $|n_0|\le X$ and $|n_1|\le X$, we find 
via orthogonality that this establishes the upper bound $\Upsilon_{s,k}(X)\ll X^{2s+\varepsilon}$ claimed in equation 
\eqref{1.4} of the statement of Theorem \ref{theorem1.1} in the special case $s=k(k+1)/2$. The desired conclusion for 
smaller values of $s$ follows by applying H\"older's inequality.\par

We henceforth focus on the proof of the estimate \eqref{10.6}. We put
\[
B=\biggl\lceil \frac{k\log X}{\log p}\biggr\rceil .
\]
Thus, we have $X^k\le p^B\le pX^k$. Next, by orthogonality, we observe that the mean value on the left hand side 
of \eqref{10.6} is equal to a weighted count of the solutions to the system of equations \eqref{1.3} over 
$\mathfrak Z(t^2)$, with each weight non-negative. An upper bound for this count is given by replacing the 
equations by corresponding congruences modulo $p^B$. Thus, by reversing our earlier observation utilising 
orthogonality, one obtains the upper bound
\begin{equation}\label{10.8}
\oint |F_p(\boldsymbol \alpha;X)|^{2s}\d\boldsymbol \alpha \le 
\oint_{p^B} |F_p(\boldsymbol \alpha;X)|^{2s}\d\boldsymbol \alpha .
\end{equation}

\par The next step towards the proof of \eqref{10.6} is the preparation of a $p^c$-spaced system associated with 
$\boldsymbol \varphi$. In view of the definition \eqref{10.4}, we may suppose that $\mathfrak a_n=1$ for each $n$. 
The mean value on the left hand side of \eqref{10.6} counts the $\mathfrak Z(t^2)$-solutions of the system of 
equations \eqref{1.3} with $|\mathbf u|,|\mathbf v|\le X$ satisfying
\[
\prod_{i=1}^{2s}W_t(u_i;\boldsymbol \varphi)\not\equiv 0\mmod{p},
\]
in which each solution is counted with weight \eqref{10.3}, namely $(2X+1)^{-2s}$. Suppose that $\tau>0$ is 
sufficiently small in terms of $s$ and $k$. We suppose also that $B$ is sufficiently large in terms of $\tau$, as well 
as $s$, $k$ and $\varepsilon$. Such is permissible, since we may suppose that $X$ is likewise sufficiently large. We 
put $c=\lceil \tau B\rceil$. We next sort the implicit summation defining $F_p(\boldsymbol \alpha;X)$ into congruence 
classes modulo $(p^c,t^2)$. In view of our assumptions concerning the sequence $(\boldsymbol{\mathfrak a})$ 
implied by \eqref{10.4}, we may suppose that $\mathfrak a_n=0$ whenever $n\equiv \xi\mmod{(p^c,t^2)}$ and
\[
W_t(\xi;\boldsymbol \varphi)\equiv 0\mmod{p}.
\]
For such values of $\xi$, one then has $\rho_c(\xi)=0$. With this assumption, one finds that
\[
F_p(\boldsymbol \alpha;X)=(2X+1)^{-1}\sum_{\xi\mmod{(p^c,t^2)}}\rho_c(\xi)\mathfrak f_c(\boldsymbol \alpha;\xi),
\]
and hence Lemma \ref{lemma6.3} delivers the bound
\[
\rho_0(1)^2|F_p(\boldsymbol \alpha;X)|^{2s}\le (p^{2c})^s\sum_{\xi\mmod{(p^c,t^2)}}\rho_c(\xi)^2
|\mathfrak f_c(\boldsymbol \alpha;\xi)|^{2s}.
\]
Thus, we obtain the estimate
\begin{equation}\label{10.9}
\oint_{p^B}|F_p(\boldsymbol \alpha;X)|^{2s}\d\boldsymbol \alpha \le (p^{2c})^s\rho_0(1)^{-2}
\sum_{\xi\mmod{(p^c,t^2)}}\rho_c(\xi)^2I_p(\xi),
\end{equation}
where
\[
I_p(\xi)=\oint_{p^B}|\mathfrak f_c(\boldsymbol \alpha;\xi)|^{2s}\d\boldsymbol \alpha .
\]

\par The mean value $I_p(\xi)$ counts the solutions $\mathbf y,\mathbf z\in \mathfrak Z(t^2)$ of the system of 
congruences
\begin{equation}\label{10.10}
\sum_{i=1}^s \bigl( \varphi_j(p^cy_i+\xi)-\varphi_j(p^cz_i+\xi)\bigr) \equiv 0\mmod{(p^B,t^2)}\quad (1\le j\le k),
\end{equation}
with each solution being counted with weight
\begin{equation}\label{10.11}
\rho_c(\xi)^{-2s}\prod_{i=1}^s \mathfrak a_{p^cy_i+\xi}{\overline{\mathfrak a}}_{p^cz_i+\xi}.
\end{equation}
Here, we note that we may restrict attention to the situation in which
\[
W_t(\xi;\boldsymbol \varphi)\not\equiv 0\mmod{p}\quad \text{and}\quad \rho_c(\xi)>0.
\]
In particular, it follows from \eqref{3.8} and \eqref{10.9} that
\[
\oint_{p^B}|F_p(\boldsymbol \alpha;X)|^{2s}\d\boldsymbol \alpha \le (p^{2c})^s
U_{s,k}^{B,c}(\boldsymbol{\mathfrak a}).
\]

\par We now reinterpret the system of congruences \eqref{10.10} by applying Taylor's theorem to expand the 
expression $\varphi_j(p^cv+\xi)$. When $1\le j\le k$, we find that for suitable polynomials 
$\Phi_j(u)\in \mathbb Z[t,u]$, one has
\[
\varphi_j(p^cv+\xi)-\varphi_j(\xi)\equiv \sum_{l=1}^k\omega_{lj}(p^cv)^l+(p^cv)^{k+1}\Phi_j(p^cv)
\mmod{(p^B,t^2)},
\]
in which the $\mathbb Z[t]$-coefficients $\omega_{lj}$ are defined by taking
\[
\omega_{lj}=(l!)^{-1}\varphi_j^{(l)}(\xi)\quad (1\le l,j\le k).
\]
Notice that these coefficients lie in $\mathbb Z[t]$, because $l!$ divides the $l$-th differential of a polynomial 
with coefficients lying in $\mathbb Z[t]$.\par

The determinant of the coefficient matrix $\Omega=(\omega_{lj})_{1\le l,j\le k}$ is given by the formula
\[
\det(\Omega)=W(\xi;\boldsymbol \varphi)\biggl( \prod_{l=1}^k l!\biggr)^{-1}.
\]
Since we may suppose that $p>k$, the hypothesis that $W(\xi;\boldsymbol \varphi)\not \equiv 0\mmod{(p,t)}$ 
permits us to conclude that $\det(\Omega)\not\equiv 0\mmod{(p,t)}$, and thus it follows that a multiplicative 
inverse $\Omega^{-1}$ modulo $(p,t)$ exists having coefficients in $\mathbb Z$. We now replace 
$\boldsymbol \varphi$ by $\Omega^{-1}\boldsymbol \varphi$ and $\boldsymbol \Phi$ by 
$\Omega^{-1}\boldsymbol \Phi$. This amounts to taking appropriate integral linear combinations of the 
congruences comprising \eqref{10.10}. Proceeding in this manner, we see that there is no loss of generality in 
supposing that the original coefficient matrix $\Omega$ is congruent to $I_k$ modulo $(p,t)$. Assuming now that 
this assumption applies ab initio, it follows from Lemma \ref{lemma4.2} that the matrix of coefficients $\Omega$ 
has an inverse matrix $\Omega^{-1}$ in $\mathbb Z[t]/(p^B,t^2)$ satisfying $\Omega^{-1}\equiv I_k\mmod{(p,t)}$.  
In this new scenario, we again replace $\boldsymbol \varphi$ by $\Omega^{-1}\boldsymbol \varphi$ and 
$\boldsymbol \Phi$ by $\Omega^{-1}\boldsymbol \Phi$. This process permits us the assumption that the coefficient 
matrix $\Omega$ is in fact equal to $I_k$. Hence, there exist polynomials $\Xi_j(u)\in \mathbb Z[t,u]$ having the 
property that whenever the system of congruences \eqref{10.10} is satisfied, then
\[
\sum_{i=1}^s (p^c)^j\bigl(\Psi_j(y_i)-\Psi_j(z_i)\bigr) \equiv 0\mmod{(p^B,t^2)}\quad (1\le j\le k),
\]
in which $\Psi_j(u)=u^j+p^c\Xi_j(u)$. In particular, the system of polynomials $\boldsymbol \Psi$ is $p^c$-spaced.

\par The discussion of the previous paragraph shows that in any solution $\mathbf y,\mathbf z$ of the system 
\eqref{10.10}, counted with weight \eqref{10.11}, one has the additional constraints
\[
\sum_{i=1}^s\bigl( \Psi_j(y_i)-\Psi_j(z_i)\bigr) \equiv 0 \mmod{(p^{B-kc},t^2)}\quad (1\le j\le k).
\]
Define the coefficients
\[
\mathfrak c_y(\boldsymbol \alpha)=\mathfrak a_{p^cy+\xi}E(\psi(p^cy+\xi;\boldsymbol \alpha)\bigr) .
\]
Also, write
\[
\mathfrak g_{\boldsymbol{\mathfrak c}}(\boldsymbol \alpha,\boldsymbol \beta)=
\rho_0(1;\boldsymbol{\mathfrak c})^{-1}\sum_{y\in \mathbb Z}
\mathfrak c_y(\boldsymbol \alpha)E\bigl(\beta_1\Psi_1(y)+\ldots +\beta_k\Psi_k(y)\bigr),
\]
and define the mean value
\[
J(\boldsymbol \alpha)=\oint_{p^{B-kc}}|\mathfrak g_{\boldsymbol{\mathfrak c}}
(\boldsymbol \alpha,\boldsymbol \beta)|^{2s}\d\boldsymbol \beta .
\]
Note that $\rho_0(1;{\boldsymbol{\mathfrak c}})=\rho_c(\xi;\boldsymbol{\mathfrak a})$. Then, just as in the 
argument leading to \eqref{7.19} above, one sees that
\begin{equation}\label{10.12}
I_p(\xi)=\oint_{p^B}|\mathfrak g_{\boldsymbol{\mathfrak c}}(\boldsymbol \alpha;\mathbf 0)|^{2s}
\d\boldsymbol \alpha=\oint_{p^B}J(\boldsymbol \alpha)\d\boldsymbol \alpha .
\end{equation}

Observe next that, when written in the notation defined in \eqref{3.6}, one has
\[
J(\boldsymbol \alpha)=U_{s,k}^{B-kc,\boldsymbol \Psi}(\boldsymbol{\mathfrak c}).
\]
Moreover, the system $\boldsymbol \Psi$ is $p^c$-spaced and $B-kc$ is sufficiently large in terms of $k$, $\tau$ 
and $\varepsilon$. Notice that, as essentially noted in the discussion following the statement of Hypothesis 
\ref{hypothesis3.2}, the coefficient weights $(\boldsymbol{\mathfrak a_n})$ with which we are working automatically 
satisfy both Hypotheses \ref{hypothesis3.1} and \ref{hypothesis3.2} with $\varepsilon=0$, and so likewise do the 
coefficients $(\boldsymbol{\mathfrak c}_y(\boldsymbol \alpha))$. Then Corollary \ref{corollary3.3} supplies the 
bound
\[
J(\boldsymbol \alpha)\ll p^{B\varepsilon}U_{s,k}^{B-kc,H,\boldsymbol \Psi}(\boldsymbol{\mathfrak c}),
\]
where $H=\lceil B/k\rceil -c$. On substituting this bound into \eqref{10.12}, and from there into \eqref{10.9}, we infer 
that
\begin{equation}\label{10.13}
\oint_{p^B}|F_p(\boldsymbol \alpha;X)|^{2s}\d\boldsymbol \alpha\ll (p^{2c})^sp^{B\varepsilon}
\rho_0(1;\boldsymbol{\mathfrak a})^{-2}\sum_{\xi\mmod{(p^c,t^2)}}\rho_c(\xi;\boldsymbol{\mathfrak a})^2
U_{s,k}^{B-kc,H,\boldsymbol \Psi}(\boldsymbol{\mathfrak c}).
\end{equation}

\par For ease of exposition, we temporarily abbreviate 
$\mathfrak f_H(\boldsymbol \beta ;\eta;\boldsymbol{\mathfrak c};\boldsymbol \Psi)$ to 
$\mathfrak f^*_H(\boldsymbol \beta ;\eta)$. Thus, we have
\[
\mathfrak f^*_H(\boldsymbol \beta ;\eta)=\rho_H(\eta;\boldsymbol{\mathfrak c})^{-1}
\sum_{y\equiv \eta \mmod{(p^H,t^2)}}\mathfrak c_y(\boldsymbol \alpha)
E\bigl( \beta_1\Psi_1(y)+\ldots +\beta_k\Psi_k(y)\bigr).
\]
Note here that the coefficients $\mathfrak c_y(\boldsymbol \alpha)$ are $0$ unless $p^cy+\xi=n_0+n_1t$, with 
the integers $n_0$ and $n_1$ satisfying $|n_i|\le X$. This condition we abbreviate by writing $|p^cy+\xi|\le X$. As a 
consequence of the triangle inequality, one consequently has
\[
\rho_H(\eta;\boldsymbol{\mathfrak c})|\mathfrak f^*_H(\boldsymbol \beta ;\eta)|\le 
\sum_{\substack{y\equiv \eta \mmod{(p^H,t^2)}\\ |p^cy+\xi|\le X}}|\mathfrak c_y(\boldsymbol \alpha)| \le 
\biggl( \frac{2X+1}{p^{c+H}}+1\biggr)^2.
\]
Moreover, a similar computation reveals that
\[
\rho_H(\eta;\boldsymbol{\mathfrak c})=\frac{2X+1}{p^{c+H}}+O(1).
\]
Since $c+H=\lceil B/k\rceil$, the definition \eqref{3.8} yields
\begin{align*}
U_{s,k}^{B-kc,H,\boldsymbol \Psi}(\boldsymbol{\mathfrak c})&\ll \rho_0(1;\boldsymbol{\mathfrak c})^{-2}
\sum_{\eta\mmod{(p^c,t^2)}}\rho_c(\eta;\boldsymbol{\mathfrak c})^2\bigl( 1+X/p^{c+H}\bigr)^{2s}\\
&\ll (1+X/p^{B/k})^{2s}.
\end{align*}
Hence, we deduce from \eqref{10.13} that
\begin{align*}
\oint_{p^B}|F_p(\boldsymbol \alpha;X)|^{2s}\d\boldsymbol \alpha &\ll p^{2(cs+B\varepsilon)} (1+X/p^{B/k})^{2s}
\rho_0(1;\boldsymbol{\mathfrak a})^{-2}\sum_{\xi\mmod{(p^c,t^2)}}\rho_c(\xi;\boldsymbol{\mathfrak a})^2\\
&\ll p^{(4s\tau +2\varepsilon)B}\bigl( 1+X/p^{B/k}\bigr)^{2s}.
\end{align*}

\par Recall that $\tau$ was chosen sufficiently small in terms of $s$ and $k$. Moreover, we are at liberty to choose 
$B$ to satisfy the bounds $X^k\le p^B<pX^k$. Then we conclude that, for each positive number $\delta$, one has
\[
\oint_{p^B}|F_p(\boldsymbol \alpha;X)|^{2s}\d\boldsymbol \alpha \ll p^{B\delta}\ll X^{2k\delta}.
\]
In view of \eqref{10.8}, this establishes the upper bound \eqref{10.6} for each positive number $\varepsilon$, and thus 
we may conclude that the estimate \eqref{10.7} does indeed hold.
\end{proof}

With Theorem \ref{theorem1.1} now established, it is a relatively simple matter to derive the corollaries announced 
in \S1, as we now show.

\begin{proof}[The proof of Corollary \ref{corollary1.2}]
We consider a $k$-tuple of polynomials $\varphi_j\in \mathbb Z[u]$ having the property that 
$W(u;\boldsymbol \varphi)$ is non-zero as a polynomial in $u$. Then one has
\[
\varphi_j(x+yt)\equiv \varphi_j(x)\mmod{t}\quad (1\le j\le k),
\]
so that $W_t(u;\boldsymbol \varphi)$ is non-zero as a polynomial in $u$. Thus, as a consequence of orthogonality in 
combination with Theorem \ref{theorem1.1}, we have
\[
\oint \biggl| \sum_{|n_0|\le X}\sum_{|n_1|\le X}E\bigl( \psi(n_0+n_1t;\boldsymbol \alpha)\bigr) \biggr|^{2s}
\d\boldsymbol \alpha =\Upsilon_{s,k}(X;\boldsymbol \varphi)\ll X^{2s+\varepsilon}.
\]
However, an application of Taylor's theorem reveals that
\[
\varphi_j(x+yt)\equiv \varphi_j(x)+y\varphi_j'(x)t\mmod{t^2}\quad (1\le j\le t).
\]
Thus, the definitions of \S2 ensure that with
\[
\boldsymbol \alpha=\biggl( \frac{\beta_1}{t}+\frac{\gamma_1}{t^2},\ldots ,\frac{\beta_k}{t}+\frac{\gamma_k}{t^2}
\biggr) ,
\]
one has
\[
\sum_{|n_0|\le X}\sum_{|n_1|\le X}E(\psi(n_0+n_1t;\boldsymbol \alpha))=\sum_{|n_0|\le X}\sum_{|n_1|\le X}
e(\Omega(\mathbf n;\boldsymbol \beta,\boldsymbol \gamma)),
\]
where
\[
\Omega(\mathbf n;\boldsymbol \beta,\boldsymbol \gamma)=\sum_{i=1}^k
\bigl( \beta_i\varphi_i(n_0)+\gamma_in_1\varphi_i'(n_0)\bigr) .
\]
In particular, it follows via orthogonality that
\begin{align*}
\Upsilon_{s,k}(X;\boldsymbol \varphi)&=\oint \biggl| \sum_{|n_0|\le X}\sum_{|n_1|\le X}
E(\psi(n_0+n_1t;\boldsymbol \alpha))\biggr|^{2s}\d\boldsymbol \alpha \\
&=\oint \oint \biggl| \sum_{|n_0|\le X}\sum_{|n_1|\le X}
e(\Omega(\mathbf n;\boldsymbol \beta,\boldsymbol \gamma))\biggr|^{2s}\d\boldsymbol \beta \d\boldsymbol \gamma .
\end{align*}
By orthogonality, the mean value on the right hand side counts the integral solutions of the system of equations
\[
\left. \begin{aligned} \sum_{i=1}^s\bigl( \varphi_j(z_i)-\varphi_j(z_{i+s})\bigr) &=0\\
\sum_{i=1}^s \bigl( h_i\varphi_j'(z_i)-h_{i+s}\varphi_j'(z_{i+s}\bigr) &=0\end{aligned} \right\} \quad 
(1\le j\le k),
\]
with $|\mathbf z|,|\mathbf h|\le X$. On invoking the conclusion of Theorem \ref{theorem1.1}, we have therefore shown 
that 
\[
\Upsilon_{s,k}^*(X;\boldsymbol \varphi)=\Upsilon_{s,k}(X;\boldsymbol \varphi)\ll X^{2s+\varepsilon},
\]
completing the proof of the corollary.
\end{proof}

A specialisation of the polynomials $\varphi_j(u)$ occurring in Corollary \ref{corollary1.2} yields the conclusion of 
Corollary \ref{corollary1.3}.

\begin{proof}[The proof of Corollary \ref{corollary1.3}]
We may suppose that $d_1,\ldots ,d_k$ are integers with
\[
1\le d_1<d_2<\ldots <d_k.
\]
A comparison of the systems \eqref{1.5} and \eqref{1.7} shows that the quantity $\Upsilon_{s,k}^+(X;\mathbf d)$ 
central to Corollary \ref{corollary1.3} is equal to the quantity $\Upsilon_{s,k}^*(X;\boldsymbol \varphi)$ estimated 
in Corollary \ref{corollary1.2}, provided that one sets $\varphi_j(u)=u^{d_j}$ $(1\le j\le k)$. In order to apply 
Corollary \ref{corollary1.2}, all that must be verified is that $W(u;\boldsymbol \varphi)$ is non-zero as a 
polynomial in $u$. However, one has
\[
W(u;\boldsymbol \varphi)=\det\bigl( d_j(d_j-1)\cdots (d_j-i+1)u^{d_j-i}\bigr)_{1\le i,j\le k}.
\]
Every monomial in this polynomial is a multiple of $u^D$, where
\[
D=\biggl( \sum_{j=1}^kd_j\biggr)-k(k+1)/2=\sum_{j=1}^k(d_j-j).
\]
Thus $W(u;\boldsymbol \varphi)=u^D\det(\Omega)$, where
\[
\Omega =\bigl( d_j(d_j-1)\cdots (d_j-i+1)\bigr)_{1\le i,j\le k}.
\]
However, as was shown in the proof of Corollary 1.2 of \cite[\S12]{Woo2019}, we have
\[
\det(\Omega)=d_1\cdots d_k\prod_{1\le i<j\le k}(d_i-d_j)\ne 0.
\]
Therefore, one has $W(u;\boldsymbol \varphi)\ne 0$. Thus, we may conclude from Corollary \ref{corollary1.2} that
\[
\Upsilon_{s,k}^+(X;\mathbf d)=\Upsilon_{s,k}^*(X;\boldsymbol \varphi)\ll X^{2s+\varepsilon},
\]
completing the proof of the corollary.
\end{proof}

Finally, the proof of Corollary \ref{corollary1.4} is essentially a specialisation of Corollary \ref{corollary1.3}.

\begin{proof}[The proof of Corollary \ref{corollary1.4}]
A comparison of the systems \eqref{1.7} and \eqref{1.8} reveals that on putting $\mathbf d=(1,2,\ldots ,k)$, one has 
$J_{s,k}^*(X)=\Upsilon_{s,k}^*(X;\mathbf d)$. Thus we deduce from Corollary \ref{corollary1.3} that with 
$s_0=k(k+1)/2$, one has
\begin{equation}\label{10.14}
J_{s_0,k}^*(X)\ll X^{2s_0+\varepsilon}=X^{k(k+1)+\varepsilon}.
\end{equation}
Put
\[
\psi_k(u,v;\boldsymbol \alpha,\boldsymbol \beta)=\sum_{j=1}^k\bigl( u^j\alpha_j+vu^{j-1}\beta_j\bigr),
\]
and
\[
w_k(\boldsymbol \alpha,\boldsymbol \beta;X)=\sum_{|u|\le X}\sum_{|v|\le X}e(\psi_k(u,v;\boldsymbol \alpha,
\boldsymbol \beta)).
\]
Then it follows via orthogonality that
\[
J_{s,k}^*(X)=\oint\oint |w_k(\boldsymbol \alpha,\boldsymbol \beta ;X)|^{2s}
\d\boldsymbol \alpha \d \boldsymbol \beta .
\]

\par When $s>s_0$, we may make use of the trivial estimate
\[
|w_k(\boldsymbol \alpha,\boldsymbol \beta ;X)|\le (2X+1)^2
\]
in combination with \eqref{10.14} to deduce that
\begin{align*}
J_{s,k}^*(X)&\ll (X^2)^{2s-2s_0}\oint\oint |w_k(\boldsymbol \alpha,\boldsymbol \beta ;X)|^{2s_0}
\d\boldsymbol \alpha \d \boldsymbol \beta \\
&=(X^2)^{2s-k(k+1)}J_{s_0,k}(X)\ll X^{4s-k(k+1)+\varepsilon}.
\end{align*}
This confirms the conclusion of Corollary \ref{corollary1.4} when $s>s_0$. Meanwhile, when $0<s\le s_0$, we may 
apply H\"older's inequality to deduce in like manner that
\begin{align*}
J_{s,k}^*(X)&\ll \biggl( \oint\oint |w_k(\boldsymbol \alpha,\boldsymbol \beta ;X)|^{2s_0}
\d\boldsymbol \alpha \d \boldsymbol \beta \biggr)^{s/s_0} \\
&=\bigl( J_{s_0,k}(X)\bigr)^{s/s_0}\ll X^{2s+\varepsilon}.
\end{align*}
The conclusion of Corollary \ref{corollary1.4} consequently holds for all positive values of $s$.
\end{proof}

\section{A number field variant}
The main conclusions of \S1 all have their natural analogues in algebraic number fields. A key advantage of our 
efficient congruencing approach, with regard to such variants, is that the proofs of such conclusions proceed without 
any significant alteration of the accompanying argument, once the associated infrastructure has been adjusted.\par

We begin by modifying the infrastructure necessary to describe our number field variants. We consider an algebraic 
extension $K$ of degree $n$ over $K$. We write $\mathfrak O_K$ for the ring of integers of $K$, and we fix a basis 
$\Omega=\{ \omega_1,\ldots ,\omega_n\}$ for $\mathfrak O_K$ over $\mathbb Z$. When $X$ is a large positive 
number, we define the box $\mathcal B(X)\subset \mathfrak O_K$ given by
\[
\mathcal B(X)=\{ r_1\omega_1+\ldots +r_n\omega_n: 
\text{$r_i\in \bigl[-\tfrac{1}{2}X^{1/n},\tfrac{1}{2}X^{1/n}\bigr)\cap \mathbb Z$ $(1\le i\le n)$}\} .
\]
One then has $\text{card}(\mathcal B(X))\asymp X$. We consider the system of congruences \eqref{1.1} with 
\[
x_i=u_i+v_it\quad (1\le i\le 2s),
\]
and denote by $\Upsilon_{s,k}(X;\boldsymbol \varphi;K)$ the number of solutions of \eqref{1.1} with 
$\mathbf u,\mathbf v\in \mathcal B(X)$.\par

The analogue of Theorem \ref{theorem1.1} for number fields may now be announced.

\begin{theorem}\label{theorem11.1} Suppose that $\varphi_j(u)\in \mathfrak O_K[t,u]$ $(1\le j\le k)$ is a system of 
polynomials with $W(u;\boldsymbol \varphi)\not\equiv 0\mmod{t}$. Let $s$ be a positive integer with 
$s\le k(k+1)/2$. Then for each $\varepsilon>0$, one has 
$\Upsilon_{s,k}(X;\boldsymbol \varphi;K)\ll X^{2s+\varepsilon}$.
\end{theorem}

Analogues of Corollaries \ref{corollary1.2} and \ref{corollary1.3} relevant for number fields follow so directly that the 
need for elaboration is rendered void. We content ourselves by recording an analogue of Corollary \ref{corollary1.4} 
that counts solutions of the system of equations \eqref{1.8} with $u_i,v_i\in \mathcal B(X)$ $(1\le i\le 2s)$, a number 
we denote by $J_{s,k}^*(X;K)$.

\begin{corollary}\label{corollary11.2}
Let $s$ be a positive integer. Then for each $\varepsilon>0$, one has
\[
J_{s,k}^*(X;K)\ll X^\varepsilon (X^{2s}+X^{4s-k(k+1)}).
\]
\end{corollary}

The conclusion of Theorem \ref{theorem11.1} follows by applying an analogue of Theorem \ref{theorem3.2}. The 
latter requires the introduction of further infrastructure. In the setting described above, we denote by $K^{(l)}$ 
$(1\le l\le n_1)$ the real conjugate fields associated with $K$, and by $K^{(m)}$ and $K^{(m+n_2)}$ 
$(n_1+1\le m\le n_1+n_2)$ the pairs of complex conjugate fields associated with $K$. Thus, one has $n=n_1+2n_2$. 
When $\gamma\in K$, we denote by $\gamma^{(i)}$ the conjugates of $\gamma$, where 
$\gamma^{(i)}\in K^{(i)}$ $(1\le i\le n)$. Then, the trace map of relevance to us is given by
\[
{\rm{Tr}}_{K/\mathbb Q}(\gamma)=\gamma^{(1)}+\ldots +\gamma^{(n)}.
\]
When $\gamma_j\in K$ and $\theta_j\in \mathbb R$ for $1\le j\le n$, and
\[
\lambda(\boldsymbol \theta)=\theta_1\gamma_1+\ldots +\theta_n\gamma_n,
\]
we define
\[
\lambda^{(i)}(\boldsymbol \theta)=\theta_1\gamma_1^{(i)}+\ldots +\theta_n\gamma_n^{(i)}\quad (1\le i\le n).
\]
The trace map on $K$ then extends to the map defined by
\[
{\rm{Tr}}_{K/\mathbb Q}(\lambda(\boldsymbol \theta))=\lambda^{(1)}(\boldsymbol \theta)+\ldots +\lambda^{(n)}
(\boldsymbol \theta),
\]
and our basic exponential function is then given by
\[
\mathcal E(\lambda(\theta))=e\bigl( {\rm{Tr}}_{K/\mathbb Q}(\lambda(\boldsymbol \theta))\bigr) .
\]

\par Next, we recall the inverse different $\mathfrak d^{-1}$, so that
\[
\mathfrak d^{-1}=\{ \eta\in K:\text{${\rm{Tr}}_{K/\mathbb Q}(\eta\xi)\in \mathbb Z$ for all $\xi\in \mathfrak O_K$}\}.
\]
There is a basis ${\rm P}=\{\rho_1,\ldots ,\rho_n\}$ of $\mathfrak d^{-1}$ dual to $\Omega$ having the property that
\[
{\rm{Tr}}_{K/\mathbb Q}(\rho_i\omega_j)=\begin{cases} 1,&\text{when $i=j$,}\\
0,&\text{when $i\ne j$.}\end{cases}
\]
We next introduce a special case of the linear form $\lambda(\boldsymbol \theta)$ introduced above, defining 
$\alpha=\alpha(\boldsymbol \theta)$ by putting
\[
\alpha(\boldsymbol \theta)=\theta_1\rho_1+\ldots +\theta_n\rho_n.
\]
When this and allied forms appear in an $n$-fold integral, we abbreviate the $n$-fold differential 
${\rm d}\theta_1\cdots \d\theta_n$ to ${\rm d}\alpha$. We also write $\mathbb T$ for $[0,1)^n$. In this notation, 
we have the orthogonality relation
\[
\int_{\mathbb T}\mathcal E(\alpha\gamma)\d\alpha =\begin{cases}1,&\text{when $\gamma=0$,}\\
0,&\text{when $\gamma\in \mathfrak O_K\setminus \{0\}$.}\end{cases}
\]

\par The introduction of all of this notation enables us to reformulate Theorem \ref{theorem11.1} in the form of a 
mean-value estimate. Define the exponential sum
\[
F_k(\boldsymbol \alpha,\boldsymbol \beta;X)=\sum_{u\in \mathcal B(X)}\sum_{v\in \mathcal B(X)}
\mathcal E(\psi_k(u,v;\boldsymbol \alpha,\boldsymbol \beta)),
\]
in which $\psi_k(u,v;\boldsymbol \alpha,\boldsymbol \beta)$ is defined precisely as in the preamble to equation 
\eqref{1.6}. Then, as a consequence of the analogue in number fields of Corollary \ref{corollary1.2}, one has the 
following conclusion.

\begin{corollary}\label{corollary11.3}
Suppose that $\varphi_j\in \mathfrak O_K[u]$ $(1\le j\le k)$ is a system of polynomials with 
$W(u;\boldsymbol \varphi)\ne 0$. Let $s$ be a positive integer with $s\le k(k+1)/2$. Then for each $\varepsilon>0$, 
\[
\int_{\mathbb T^k}\int_{\mathbb T^k}|F_k(\boldsymbol \alpha,\boldsymbol \beta ;X)|^{2s}\d\boldsymbol \alpha 
\d \boldsymbol \beta \ll X^{2s+\varepsilon}.
\]
\end{corollary}

We provide a proof of Theorem \ref{theorem11.1} in outline only, since the interested reader will have no difficulty 
in following the plan laid out in \cite[\S15]{Woo2019} in order to modify our treatment in \S\S2--10 so as to 
substitute a number field setting in place of $\mathbb Z$.\par

Our methods involve congruences, and this requires some rapid discussion. We again consider an integer $k\ge 1$, 
and polynomials $\varphi_1,\ldots ,\varphi_k\in \mathfrak O_K[t,u]$. Throughout the argument, the ring 
$\mathfrak O_K$ replaces $\mathbb Z$. Let $\mathfrak p=(\pi)$ be a prime ideal of $\mathfrak O_K$, and write 
$p=N_{K/\mathbb Q}(\pi)$. We work implicitly in the $\mathfrak p$-adic field completing $K$ at the place 
$\mathfrak p$, and we assume throughout that $p>(k!)^n$ (for safety). Congruences modulo $p^h$ are now 
replaced throughout by congruences modulo $\mathfrak p^h$, and the function $e(z)$ is replaced by 
$\mathcal E(z)$. When $F:\mathbb T^k\rightarrow \mathbb C$ is integrable, we write
\[
\oint F(\boldsymbol \alpha)\d\boldsymbol \alpha =\int_{\mathbb T^k}F(\boldsymbol \alpha)\d\alpha .
\]
Next, when $B$ is a positive integer, we define
\begin{equation}\label{11.1}
\oint_{\mathfrak p^B}F(\boldsymbol \alpha)\d\boldsymbol \alpha =p^{-kB}\sum_{u_1\mmod{\mathfrak p^B}}\ldots 
\sum_{u_k\mmod{\mathfrak p^B}}F({\mathbf u}\pi^{-B}).
\end{equation}
Here, the summations on the right hand side of \eqref{11.1} are taken over complete sets of residues modulo 
$\mathfrak p^B$, with the argument $\mathbf u \pi^{-B}$ interpreted as follows. We interpret each coordinate 
$u_i\pi^{-B}$ in the argument on the right hand side of \eqref{11.1} in the shape
\[
u_i\pi^{-B}=\beta_{i1}\rho_1+\ldots +\beta_{in}\rho_n,
\]
with $\beta_{ij}\in \mathbb R$, and then reduce each coefficient $\beta_{ij}$ modulo $1$. With this convention, the 
expression $\mathbf u\pi^{-B}$ may be identified as belonging to $\mathbb T^k$. Integrals with subscripts $p^B$ 
in \S\S2--10 are now replaced by this newly defined integral with subscript $\mathfrak p^B$, in obvious fashion. 
Now, for $\nu\in \mathfrak O_K[t]$, one has the orthogonality relation
\[
\oint_{\mathfrak p^B}E(\nu \alpha )\d\alpha=\begin{cases} 1,&\text{when $\nu \equiv 0
\mmod{(\mathfrak p^B,t^2)}$,}\\
0,&\text{when $\nu \not\equiv 0\mmod{(\mathfrak p^B,t^2)}$.}\end{cases}
\]

\par Modified versions of Theorem \ref{theorem3.2} and its corollary may now be stated in the number field setting. 
We recall the definitions \eqref{3.6} and \eqref{3.8}, though of course now in their modified apparitions.

\begin{theorem}\label{theorem11.4}
Let $K$ be an algebraic extension of $\mathbb Q$ with $[K:\mathbb Q]=n$. Suppose that $k\in \mathbb N$, 
and that $\mathfrak p=(\pi)$ is a prime ideal of $\mathfrak O_K$ with $p=N_{K/\mathbb Q}(\pi)>(k!)^n$. Then one 
has $\lambda (k(k+1)/2,k)=0$.
\end{theorem}

\begin{corollary}\label{corollary11.5}
Let $K$ be an algebraic extension of $\mathbb Q$ with $[K:\mathbb Q]=n$. Suppose that $k\in \mathbb N$, and 
that $\mathfrak p=(\pi)$ is a prime ideal of $\mathfrak O_K$ with $p=N_{K/\mathbb Q}(\pi)>(k!)^n$. Suppose that 
$\tau>0$ and $\varepsilon>0$. Let $B$ be sufficiently large in terms of $n$, $k$, $\tau$ and $\varepsilon$. Put 
$s=k(k+1)/2$ and $H=\lceil B/k\rceil$. Then, for every $\boldsymbol \varphi\in \Phi_\tau(B)$, and every sequence 
$(\boldsymbol{\mathfrak a})\in \mathbb D_0(B)$, one has
\[
U_{s,k}^B(\boldsymbol{\mathfrak a})\ll p^{B\varepsilon}U_{s,k}^{B,H}(\boldsymbol{\mathfrak a}).
\]
\end{corollary}

With Corollary \ref{corollary11.5} in hand, the proofs of Theorem \ref{theorem11.1} and its corollaries follow just as 
in the discussion of \S10, as modified in the proof of \cite[Theorem 15.1]{Woo2019} described in 
\cite[\S15]{Woo2019}. The reader should experience no difficulty in following this pedestrian pathway.

\section{A function field variant}
The conclusions of \S1 also have their natural analogues in function fields, and again proofs are supplied relatively 
easily owing to the flexibility of efficient congruencing methods. We illustrate our ideas by considering the function 
field $\mathbb F_q(T)$. We once again begin by describing the infrastructure necessary to describe our function 
field variants.\par

Consider a finite field $\mathbb F_q$ having characteristic $p>k$, and with $q=p^l$. Associated with the polynomial 
ring $\mathbb O=\mathbb F_q[T]$ is the field of fractions $\mathbb K =\mathbb F_q(T)$. We now take as our basic 
parameter the integer $X$, a sufficiently large natural number, and we put
\[
\mathbb O_X=\{ \nu\in \mathbb F_q[T]:\text{deg}(\nu)\le X\}.
\]
In the present setting, we consider the system of congruences \eqref{1.1} with
\[
x_i=u_i+v_it\quad (1\le i\le 2s),
\]
and denote by $\Upsilon_{s,k}(X;\boldsymbol \varphi;\mathbb K)$ the number of solutions of \eqref{1.1} with 
$\mathbf u,\mathbf v\in \mathbb O_X$. The analogue of Theorem \ref{theorem1.1} for the class of function fields 
$\mathbb F_q(T)$ at hand may now be recorded.

\begin{theorem}\label{theorem12.1} Suppose that $\varphi_j(u)\in \mathbb O[t,u]$ $(1\le j\le k)$ is a system of 
polynomials with $W(u;\boldsymbol \varphi)\not \equiv 0\mmod{t}$. Let $s$ be a positive integer with 
$s\le k(k+1)/2$. Then, whenever $p>k$ and $\varepsilon>0$, one has
\[
\Upsilon_{s,k}(X;\boldsymbol \varphi;\mathbb K)\ll (q^X)^{2s+\varepsilon}.
\]
\end{theorem}

As in the corresponding situation in \S11, analogues of Corollaries \ref{corollary1.2} and \ref{corollary1.3} 
appropriate for function fields follow directly from Theorem \ref{theorem12.1}. We again record an analogue of 
Vinogradov's mean value theorem, though now in function fields. In this context, we denote by 
$J_{s,k}^*(X;\mathbb K)$ the number of solutions of the system \eqref{1.8} with $\mathbf u,\mathbf v\in \mathbb O_X$.

\begin{corollary}\label{corollary12.2}
Let $s$ be a positive integer. Then, whenever $p>k$ and $\varepsilon>0$, one has
\[
J_{s,k}^*(X;\mathbb K)\ll (q^X)^\varepsilon \bigl( (q^X)^{2s}+(q^X)^{4s-k(k+1)}\bigr) .
\]
\end{corollary}

The conclusion of Theorem \ref{theorem12.1} follows by applying an analogue of Theorem \ref{theorem3.2}. In 
order to describe this analogue, we require some additional infrastructure. We denote by 
$\mathbb K_\infty=\mathbb F_q((1/T))$ the completion of $\mathbb F_q(T)$ at $\infty$. One may write each 
element $\alpha \in \mathbb K_\infty$ in the shape $\alpha=\sum_{i\le n}a_iT^i$ for some $n\in \mathbb Z$ and 
coefficients $a_i=a_i(\alpha)$ in $\mathbb F_q$ for $-\infty <i\le n$. We define $\text{ord}\, \alpha$ to the largest 
integer $i$ for which $a_i(\alpha)\ne 0$. We then write $\langle \alpha\rangle$ for $q^{\text{ord}\,\alpha}$. Here, 
we adopt the convention that $\text{ord}\, 0=-\infty$ and $\langle 0\rangle =0$.\par

Consider the compact additive subgroup $\mathbb T$ of $\mathbb K_\infty$ defined by $\mathbb T=\{ \alpha \in 
\mathbb K_\infty: \langle \alpha\rangle<1\}$. Every element $\alpha$ of $\mathbb K_\infty$ can be written uniquely 
in the shape $\alpha =\lfloor \alpha\rfloor +\|\alpha\|$, where $\lfloor \alpha\rfloor\in \mathbb F_q[T]$ and 
$\|\alpha \|\in \mathbb T$, and we may normalise any Haar measure ${\rm d}\alpha$ on $\mathbb K_\infty$ in such 
a manner that $\int_\mathbb T\d\alpha =1$.\par

In order to define the analogue of the exponential function, we begin by observing that there is a non-trivial 
character $e_q:\mathbb F_q\rightarrow \mathbb C^\times$ defined for each $a\in \mathbb F_q$ by taking 
$e_q(a)=e({\rm tr}(a)/p)$, where ${\rm tr}(a)=a+a^p+\ldots +a^{p^{l-1}}$ denotes the familiar trace map. This 
character induces a map $e:\mathbb K_\infty \rightarrow \mathbb C^\times$ by defining, for each element 
$\alpha\in \mathbb K_\infty$, the value of $e(\alpha)$ to be $e_q(\text{res}(\alpha))$. One then has
\[
\int_\mathbb T e(h\alpha)\d\alpha =\begin{cases}1,&\text{when $h=0$,}\\
0,&\text{when $h\in \mathbb F_q[T]\setminus \{0\}$.}\end{cases}
\]

\par Just as in the corresponding situation in \S11, the introduction of this notation allows us to reformulate 
Theorem \ref{theorem12.1} in the form of a mean value estimate. Define now
\[
F_k(\boldsymbol \alpha,\boldsymbol \beta;X)=\sum_{u\in \mathbb O_X}\sum_{v\in \mathbb O_X}
e(\psi_k(u,v;\boldsymbol \alpha,\boldsymbol \beta)),
\]
in which $\psi_k(u,v;\boldsymbol \alpha,\boldsymbol \beta)$ is again defined as in the preamble to equation 
\eqref{1.6}. Then, as a consequence of the analogue in function fields of Corollary \ref{corollary1.2}, one has the 
following conclusion.

\begin{corollary}\label{corollary12.3}
Suppose that $\varphi_j\in \mathbb O[u]$ $(1\le j\le k)$ is a system of polynomials with 
$W(u;\boldsymbol \varphi)\ne 0$. Let $s$ be a positive integer with $s\le k(k+1)/2$. Then for each $\varepsilon>0$, 
one has
\[
\int_{\mathbb T^k}\int_{\mathbb T^k}|F_k(\boldsymbol \alpha,\boldsymbol \beta;X)|^{2s}
\d\boldsymbol \alpha \d\boldsymbol \beta \ll (q^X)^{2s+\varepsilon}.
\]
\end{corollary}

We again provide a proof of Theorem \ref{theorem12.1} in brief only. The reader will experience no difficulty in 
following the plan laid out in \cite[\S17]{Woo2019} in order to modify our treatment in \S\S2--10 so as to replace the 
ring of integers $\mathbb Z$ in the former setting by its function field analogue $\mathbb O=\mathbb F_q[T]$ 
presently the focus of attention.\par

We again consider an integer $k\ge 1$, and polynomials $\varphi_1(u),\ldots ,\varphi_k(u)\in \mathbb O[t,u]$. 
Throughout the argument, the ring $\mathbb O$ replaces $\mathbb Z$. Let $\pi\in \mathbb O$ be a monic 
irreducible polynomial, having positive degree of course. We assume throughout, as previously noted, that 
$\text{char}(\mathbb F_q)=p>k$. Congruences modulo $p^h$ are now replaced throughout by congruences 
modulo $\pi^h$, and the function $e(z)$ of \S\S2--10 is replaced by its surrogate from this section. When 
$F:\mathbb T^k\rightarrow \mathbb C$ is integrable, we write
\[
\oint F(\boldsymbol \alpha)\d\boldsymbol \alpha =\int_{\mathbb T^k}F(\boldsymbol \alpha)\d\boldsymbol \alpha .
\]
Also, when $B$ is a positive integer, we define
\[
\oint_{\pi^B}F(\boldsymbol \alpha)\d\boldsymbol \alpha =\langle \pi\rangle^{-kB}\sum_{u_1\mmod{\pi^B}}
\cdots \sum_{u_k\mmod{\pi^B}}F(\mathbf u \pi^{-B}),
\]
where the summations are taken over complete sets of residues modulo $\pi^B$. Integrals with subscripts $p^B$ in 
\S\S2--10 are now replaced by this newly defined integral with subscript $\pi^B$, in obvious fashion. Now, for 
$\nu\in \mathbb O[t]$, one has the orthogonality relation
\[
\oint_{\pi^B}E(\alpha \nu)\d\alpha =\begin{cases}1,&\text{when $\nu \equiv 0\mmod{(\pi^B,t^2)}$,}\\
0,&\text{when $\nu\not\equiv 0\mmod{(\pi^B,t^2)}$.}\end{cases}
\]
Modified versions of Theorem \ref{theorem3.2} and its corollary may now be stated in the function field setting. We 
recall the definitions \eqref{3.6} and \eqref{3.8}, now in their modified attire.

\begin{theorem}\label{theorem12.4}
Let $\mathbb O=\mathbb F_q[T]$. Suppose that $k\in \mathbb N$, and that $\pi\in \mathbb O$ is a monic 
irreducible polynomial. Then, under the assumption that ${\rm ch}(\mathbb F_q)>k$, one has 
\[
\lambda(k(k+1)/2,k)=0.
\]
\end{theorem}

\begin{corollary}\label{corollary12.5}
Let $\mathbb O=\mathbb F_q[T]$. Suppose that $k\in \mathbb N$, and that $\pi\in \mathbb O$ is a monic 
irreducible polynomial. Suppose also that $\tau>0$ and $\varepsilon>0$. Let $B$ be sufficiently large in terms of 
$k$, $\tau$ and $\varepsilon$. Put $s=k(k+1)/2$ and $H=\lceil B/k\rceil$. Finally, suppose that 
${\rm ch}(\mathbb F_q)>k$. Then for every $\boldsymbol \varphi\in \Phi_\tau(B)$, and every sequence 
$(\boldsymbol{\mathfrak a})\in \mathbb D_0(B)$, one has
\[
U_{s,k}^B(\boldsymbol{\mathfrak a})\ll \langle \pi\rangle^{B\varepsilon}U_{s,k}^{B,H}(\boldsymbol{\mathfrak a}).
\]
\end{corollary}

Granted Corollary \ref{corollary12.5}, the proofs of Theorem \ref{theorem12.1} and its corollaries follow just as in the 
discussion of \S10, as modified in the proof of \cite[Theorem 17.1]{Woo2019} described in \S17 of the latter source. 
The well-travelled reader will experience no discomfort in pursuing this pedestrian route.

\section{Appendix: minimal TDI systems in dimension $2$}
Our goal in this section is to provide an informal classification of translation-dilation invariant (TDI) systems in two 
dimensions having the minimal number of equations. Such systems as concern us are generated by polynomials 
$\varphi_j(x,y),\psi_j(x,y)\in \mathbb C[x,y]$ $(1\le j\le k)$, with $k\ge 2$ and
\[
\text{deg}\,\varphi_j=\text{deg}\,\psi_j\quad (1\le j\le k),
\]
with
\[
1\le \text{deg}\, \varphi_1<\ldots <\text{deg}\, \varphi_k.
\]
We shall suppose that these polynomials satisfy the condition that whenever $s\in \mathbb N$, then the system of 
equations
\begin{equation}\label{13.1}
\left. \begin{aligned}
\sum_{i=1}^s\bigl( \varphi_j(x_i,y_i)-\varphi_j(x_{s+i},y_{s+i})\bigr)&=0\\
\sum_{i=1}^s\bigl( \psi_j(x_i,y_i)-\psi_j(x_{s+i},y_{s+i})\bigr)&=0
\end{aligned}\right\}\quad (1\le j\le k)
\end{equation}
is translation-dilation invariant. By this we mean that whenever $\mathbf x,\mathbf y$ is a solution of the system 
\eqref{13.1}, then for any $q\in \mathbb C^\times$ and $a,b\in \mathbb C$, the tuple 
$q\mathbf x+a,q\mathbf y+b$ is also a solution of the system \eqref{13.1}. We show that, up to reasonable notions 
of equivalence and non-degeneracy of systems, all systems of the shape \eqref{13.1} fall into one of three families 
simple to describe.\par

We observe first that, by taking appropriate linear combinations of the polynomials $\varphi_j(x,y)$ and $\psi_j(x,y)$, 
for $1\le j\le k$, we may suppose without loss of generality that these polynomials are both homogeneous. This is an 
immediate consequence of the dilation invariance condition associated with the system \eqref{13.1}, specifically that 
$q\mathbf x,q\mathbf y$ should be a solution of \eqref{13.1} whenever $\mathbf x,\mathbf y$ is a solution.\par

The TDI property also ensures that for $1\le j\le k$, the polynomials
\[
\frac{\partial \varphi_j}{\partial x}(x,y),\quad \frac{\partial \varphi_j}{\partial y}(x,y),\quad 
\frac{\partial \psi_j}{\partial x}(x,y),\quad \frac{\partial \psi_j}{\partial y}(x,y),
\]
must all lie in the set
\[
\bigl\{ \lambda \varphi_{j-1}(x,y)+\mu \psi_{j-1}(x,y): \lambda,\mu\in \mathbb C\bigr\} .
\]
In order to confirm this observation, one has merely to observe that a Taylor series expansion of 
$\varphi_j(x+a,y+b)$ takes the shape
\[
\varphi_j(x,y)+\Bigl( a\frac{\partial \varphi_j}{\partial x}(x,y)+b\frac{\partial \varphi_j}{\partial y}(x,y)\Bigr) +\ldots ,
\]
with a similar expansion holding for $\psi_j(x+a,y+b)$. Then, on substituting this expansion into the system 
\eqref{13.1}, we see that the TDI property of this system requires that
\[
\sum_{i=1}^s\Bigl( \frac{\partial \varphi_j}{\partial x}(x_i,y_i)-\frac{\partial \varphi_j}{\partial x}(x_{s+i},y_{s+i})\Bigr)=0
\]
and
\[
\sum_{i=1}^s\Bigl( \frac{\partial \varphi_j}{\partial y}(x_i,y_i)-\frac{\partial \varphi_j}{\partial y}(x_{s+i},y_{s+i})\Bigr)=0,
\]
with a similar property associated to $\psi_j(x,y)$. It follows that
\begin{equation}\label{13.2}
\text{span}_\mathbb C\biggl\{ \frac{\partial \varphi_j}{\partial x},\frac{\partial \varphi_j}{\partial y}, 
\frac{\partial \psi_j}{\partial x},\frac{\partial \psi_j}{\partial y}\biggr\}\subseteq 
\text{span}_\mathbb C\{\varphi_{j-1},\psi_{j-1}\}\quad (2\le j\le k),
\end{equation}
whence, in particular,
\begin{equation}\label{13.3}
\dim \text{span}_\mathbb C\biggl\{ \frac{\partial \varphi_j}{\partial x},\frac{\partial \varphi_j}{\partial y}, 
\frac{\partial \psi_j}{\partial x},\frac{\partial \psi_j}{\partial y}\biggr\}\le 2\quad (2\le j\le k).
\end{equation}

\par We concentrate initially on the consequences of the relations \eqref{13.2} and \eqref{13.3} relating to the 
polynomials $\varphi_j$. By a second application of \eqref{13.2}, we find that
\[
\dim \text{span}_\mathbb C\biggl\{ \frac{\partial^2\varphi_k}{\partial x^2},
\frac{\partial^2\varphi_k}{\partial x\,\partial y},\frac{\partial^2\varphi_k}{\partial y^2}\biggr\}\le 2.
\]
Thus, there exist $a,b,c\in \mathbb C$, not all zero, such that
\begin{equation}\label{13.4}
a\frac{\partial^2\varphi_k}{\partial x^2}+b\frac{\partial^2\varphi_k}{\partial x\,\partial y}+
c\frac{\partial^2\varphi_k}{\partial y^2}=0.
\end{equation}
For ease of exposition, write
\[
\varphi_k^{(r)}=\frac{\partial^k\varphi_k}{\partial x^r\, \partial y^{k-r}}\quad (0\le r\le k),
\]
and note that each quantity $\varphi_k^{(r)}$ is a constant. Then we deduce from \eqref{13.4} that
\begin{equation}\label{13.5}
a\varphi_k^{(r+2)}+b\varphi_k^{(r+1)}+c\varphi_k^{(r)}=0\quad (0\le r\le k-2).
\end{equation}

\par We divide our further discussion into three cases:\vskip.1cm

\noindent (a) Suppose first that the relation \eqref{13.4} holds with $a\ne 0$. In this case we examine the 
characteristic equation
\begin{equation}\label{13.6}
a\lambda^2+b\lambda+c=0
\end{equation}
of the recurrence relation \eqref{13.5}. Let $\lambda_1$ and $\lambda_2$ be the roots of the equation \eqref{13.6}. 
Then, for some $A,B\in \mathbb C$, one has
\[
\varphi_k^{(r)}=\begin{cases} A\lambda_1^r+B\lambda_2^r,&\text{when $\lambda_1\ne \lambda_2$,}\\
(Ar+B)\lambda_1^{r},&\text{when $\lambda_1=\lambda_2$.}\end{cases}
\]

\par When $\lambda_1\ne \lambda_2$, this observation reveals that the coefficient $c_r$ of $x^ry^{k-r}$ in 
$\varphi_k(x,y)$ satisfies
\[
r!(k-r)!c_r=\varphi_k^{(r)}=A\lambda_1^r+B\lambda_2^r,
\]
and hence
\begin{equation}\label{13.7}
k!\varphi_k(x,y)=\sum_{r=0}^k\binom{k}{r}\varphi_k^{(r)}x^ry^{k-r}=A(\lambda_1x+y)^k+B(\lambda_2x+y)^k.
\end{equation}
When $\lambda_1=\lambda_2$, meanwhile, the recurrence relation has a singular solution. One finds in this scenario 
that the coefficient $c_r$ of $x^ry^{k-r}$ in $\varphi_k(x,y)$ satisfies
\[
r!(k-r)!c_r=\varphi_k^{(r)}=(Ar+B)\lambda_1^r,
\]
whence
\begin{equation}\label{13.8}
k!\varphi_k(x,y)=\sum_{r=0}^k\binom{k}{r}\varphi_k^{(r)}x^ry^{k-r}=
Ak\lambda_1x(\lambda_1x+y)^{k-1}+B(\lambda_1x+y)^k.
\end{equation}
\vskip.1cm

\noindent (b) Next we suppose that the relation \eqref{13.4} holds with $a=c=0$, whence $b\ne 0$. Then we have
\[
\frac{\partial^2\varphi_k}{\partial x\,\partial y}=0,
\]
so that for some $A,B\in \mathbb C$, one has $\varphi_k(x,y)=Ax^k+By^k$.\vskip.1cm

\vskip.1cm\noindent (c) Finally, we suppose that the relation \eqref{13.4} holds with $a=0$ and $c\ne 0$. In this 
situation we may interchange the roles of $x$ and $y$ to obtain formulae of the shapes \eqref{13.7} and \eqref{13.8}, 
save that $x$ and $y$ have been interchanged.\par

A similar argument applies for the polynomial $\psi_k(x,y)$. Thus, we see that for suitable complex numbers $A'$ and 
$B'$, and for certain characteristic equation roots $\lambda_1'$ and $\lambda_2'$, the polynomial $\psi_k(x,y)$ 
takes one of the forms
\[
A'x^k+B'y^k,\quad A'(\lambda_1'x+y)^k+B'(\lambda_2'x+y)^k,\quad 
A'k\lambda_1'x(\lambda_1'x+y)^{k-1}+B'(\lambda_1'x+y)^k,
\]
or similar expressions with the roles of $x$ and $y$ interchanged.\par

We now return to examine the consequences of our classification for the relation
\[
\dim \text{span}_\mathbb C\biggl\{ \frac{\partial \varphi_k}{\partial x},\frac{\partial \varphi_k}{\partial y}, 
\frac{\partial \psi_k}{\partial x},\frac{\partial \psi_k}{\partial y}\biggr\}\le 2
\]
implied by \eqref{13.3}. Note that we also require that
\[
\dim \text{span}_\mathbb C\{ \varphi_k,\psi_k\}=2.
\]
There is no loss of generality in replacing $(\varphi_k,\psi_k)$ by the pair $(\Phi_k,\Psi_k)$, where 
\[
\Phi_k=\alpha \varphi_k+\beta \psi_k\quad \text{and}\quad \Psi_k=\gamma \varphi_k+\delta \psi_k,
\]
for any
\[
\left( \begin{matrix} \alpha&\beta\\ \gamma&\delta\end{matrix}\right) \in {\rm{GL}}_2(\mathbb C).
\]
Likewise, we may replace $\varphi_k(x,y)$ and $\psi_k(x,y)$ by
\[
\varphi_k(\mu x+\nu y,\rho x+\sigma y)\quad \text{and}\quad \psi_k(\mu x+\nu y,\rho x+\sigma y)
\]
for any
\[
\left( \begin{matrix} \mu&\nu\\ \rho&\sigma\end{matrix}\right) \in {\rm{GL}}_2(\mathbb C).
\]
With such transformations in mind, a modicum of contemplation reveals that our classification discussed thus far 
implies that, up to the equivalence implicit in these ${\rm GL_2}(\mathbb C)$ transformations just considered, the 
polynomial pairs $(\varphi_k(x,y),\psi_k(x,y))$ must take one of the forms
\[
(x^k,y^k)\quad \text{and}\quad (x^k,x^{k-1}y).
\]
It then follows that the systems of polynomials defining the equations \eqref{13.1} are equivalent to ones taking the 
form
\[
(\varphi_j(x,y),\psi_j(x,y))=(x^j,y^j)\quad (1\le j\le k),
\]
or
\[
(\varphi_j(x,y),\psi_j(x,y))=(x^j,x^{j-1}y)\quad (1\le j\le k).
\]

\par One additional detail now briefly detains us, for we have considered polynomials defined over $\mathbb C[x,y]$, 
and ones equivalent via ${\rm GL}_2(\mathbb C)$ transformations. If we seek to restrict attention to polynomials 
defined over $\mathbb Z[x,y]$, then similar arguments apply, although our transformations must apply in such a 
manner that the resulting polynomials remain defined over $\mathbb Z[x,y]$. Since the characteristic roots 
$\lambda_1,\lambda_2$ introduced above must lie in a quadratic extension of $\mathbb Q$, when 
$\lambda_1\ne \lambda_2$, and must lie in $\mathbb Q$, when $\lambda_1=\lambda_2$, we arrive at the 
following conclusion. Thus, the polynomials $(\varphi_j(x,y),\psi_j(x,y))$ must be equivalent via 
${\rm{GL}}_2(\mathbb Q)$ transformations to pairs taking one of the following forms, for $1\le j\le k$:
\begin{enumerate}
\item[(a)] When $\lambda_1$ and $\lambda_2$ are distinct rational numbers, then
\[
\varphi_j(x,y)=x^j\quad \text{and}\quad \psi_j(x,y)=y^j.
\]
\item[(b)] When $\lambda_1=\lambda_2$, a rational number, then
\[
\varphi_j(x,y)=x^j\quad \text{and}\quad \psi_j(x,y)=x^{j-1}y.
\]
\item[(c)] When $\lambda_1$ and $\lambda_2$ are distinct, but neither are rational, then necessarily 
$\lambda_1$ and $\lambda_2$ are conjugate quadratic irrational numbers in a number field of the shape 
$\mathbb Q(\sqrt{d})$, for some non-square integer $d$. One then has
\[
\varphi_j(x,y)=(x+y\sqrt{d})^j+(x-y\sqrt{d})^j
\]
and
\[
\psi_j(x,y)=\frac{1}{\sqrt{d}}\bigl( (x+y\sqrt{d})^j-(x-y\sqrt{d})^j\bigr),
\]
so that
\[
\varphi_j(x,y)+\sqrt{d}\psi_j(x,y)=(x+y\sqrt{d})^j.
\]
\end{enumerate}

\noindent So far as mean values of the corresponding exponential sums are concerned, the scenario (a) is 
accommodated via the classical theory resolved in work of the author \cite{Woo2016} and \cite{Woo2019}, and also 
\cite{BDG2016}. Meanwhile, the scenario (c) is handled via the extension of Vinogradov's mean value theorem to 
quadratic number fields, a matter handled for quite general number fields in the author's work \cite[\S15]{Woo2019}. 
Finally, the scenario (b) is the subject of this memoir.\par

We finish by observing that the system \eqref{13.1}, in the scenarios (a), (b) and (c), correspond to systems of 
Vinogradov type over $\mathbb Z[t]$ modulo $\omega(t)$, in the following sense. Thus, scenario (a) corresponds to 
the Vinogradov system of congruences
\begin{equation}\label{13.9}
\sum_{i=1}^s\bigl( (x_i+y_it)^j-(x_{s+i}+y_{s+i}t)^j\bigr)\equiv 0\mmod{\omega(t)}\quad (1\le j\le k),
\end{equation}
wherein $\omega(t)$ is a quadratic polynomial in $\mathbb Z[t]$ having two distinct rational roots. Meanwhile, 
scenario (b) corresponds to the Vinogradov system of congruences \eqref{13.9} in which $\omega(t)$ is a quadratic 
polynomial in $\mathbb Z[t]$ having repeated (rational) roots. Finally, scenario (c) corresponds to the Vinogradov 
system of congruences \eqref{13.9} in which $\omega(t)$ is a quadratic polynomial in $\mathbb Z[t]$ irreducible 
over $\mathbb Q$, but splitting in $\mathbb Q(\sqrt{d})$.

\bibliographystyle{amsbracket}
\providecommand{\bysame}{\leavevmode\hbox to3em{\hrulefill}\thinspace}

\end{document}